\documentclass[11pt,letter]{amsart}
\usepackage[latin1]{inputenc} 
\usepackage{multicol}

\newtheorem{theo}{Theorem}[section]
\newtheorem{prop}[theo]{Proposition}
\newtheorem{cor}[theo]{Corollary}
\newtheorem{lemma}[theo]{Lemma}

\theoremstyle{definition}
\newtheorem{rque}[theo]{Remark}
\newtheorem{defin}[theo]{Definition}
\newtheorem{ex}[theo]{Example}

\usepackage[all]{xy}
\usepackage{amsfonts}
\usepackage{amsmath}
\usepackage{amssymb}
\usepackage{graphics}

\DeclareMathOperator{\A}{\mathbb{A}}
\DeclareMathOperator{\C}{\mathbb{C}}

\DeclareMathOperator{\F}{\mathbb{F}}
\DeclareMathOperator{\Q}{\mathbb{Q}}
\DeclareMathOperator{\R}{\mathbb{R}}
\DeclareMathOperator{\Z}{\mathbb{Z}}
\DeclareMathOperator{\n}{\mathbb{N}}
\DeclareMathOperator{\p}{\mathbb{P}}

\DeclareMathOperator{\T}{\mathbb{T}}
\DeclareMathOperator{\bV}{\mathbb{V}}

\DeclareMathOperator{\cA}{\mathcal{A}}
\DeclareMathOperator{\cF}{\mathcal{F}}
\DeclareMathOperator{\cH}{\mathcal{H}}
\DeclareMathOperator{\cM}{\mathcal{M}}

\DeclareMathOperator{\cO}{\mathcal{O}}
\DeclareMathOperator{\cP}{\mathcal{P}}

\DeclareMathOperator{\cT}{\mathcal{T}}
\DeclareMathOperator{\cU}{\mathcal{U}}
\DeclareMathOperator{\cV}{\mathcal{V}}
\DeclareMathOperator{\cW}{\mathcal{W}}

\DeclareMathOperator{\ga}{\mathfrak{a}}
\DeclareMathOperator{\gb}{\mathfrak{b}}
\DeclareMathOperator{\gc}{\mathfrak{c}}
\DeclareMathOperator{\gd}{\mathfrak{d}}
\DeclareMathOperator{\gm}{\mathfrak{m}}
\DeclareMathOperator{\gM}{\mathfrak{M}}
\DeclareMathOperator{\gn}{\mathfrak{n}}
\DeclareMathOperator{\go}{\mathfrak{o}}
\DeclareMathOperator{\gO}{\mathfrak{O}}
\DeclareMathOperator{\gp}{\mathfrak{p}}
\DeclareMathOperator{\gt}{\mathfrak{t}}
\DeclareMathOperator{\gu}{\mathfrak{u}}

\DeclareMathOperator{\Ad}{Ad}

\DeclareMathOperator{\Ass}{Ass}
\DeclareMathOperator{\Aut}{Aut}

\DeclareMathOperator{\Cl}{Cl}

\DeclareMathOperator{\End}{End}

\DeclareMathOperator{\Fil}{Fil}

\DeclareMathOperator{\Frob}{Frob}
\DeclareMathOperator{\GL}{GL}

\DeclareMathOperator{\Gal}{Gal}
\DeclareMathOperator{\Gm}{\mathbb{G}_m}
\DeclareMathOperator{\G}{\mathcal{G}}
\DeclareMathOperator{\Hom}{Hom}
\DeclareMathOperator{\Ind}{Ind}
\DeclareMathOperator{\Isom}{Isom}

\DeclareMathOperator{\MF}{MF}

\DeclareMathOperator{\Nm}{\N_{F\!/\!\Q}}
\DeclareMathOperator{\N}{N}
\DeclareMathOperator{\PGL}{PGL}
\DeclareMathOperator{\PSL}{PSL}

\DeclareMathOperator{\Rep}{Rep}
\DeclareMathOperator{\Res}{Res}

\DeclareMathOperator{\SL}{SL}
\DeclareMathOperator{\SO}{SO}

\DeclareMathOperator{\Spec}{Spec}

\DeclareMathOperator{\St}{St}
\DeclareMathOperator{\Stab}{Stab}
\DeclareMathOperator{\Supp}{Supp}
\DeclareMathOperator{\Sym}{Sym}

\DeclareMathOperator{\Tor}{Tor}

\DeclareMathOperator{\an}{an}
\DeclareMathOperator{\bs}{\!\backslash\!}

\DeclareMathOperator{\crys}{crys}

\DeclareMathOperator{\disc}{disc}
\DeclareMathOperator{\dlog}{dlog}
\DeclareMathOperator{\dR}{dR}

\DeclareMathOperator{\fil}{fil}
\DeclareMathOperator{\gal}{gal}

\DeclareMathOperator{\gr}{gr}
\DeclareMathOperator{\id}{id}
\DeclareMathOperator{\im}{im}
\DeclareMathOperator{\limind}{\underset{\rightarrow}{lim}} 
\DeclareMathOperator{\limproj}{\underset{\leftarrow}{lim}}
\DeclareMathOperator{\ldR}{log-dR}
\DeclareMathOperator{\ord}{ord}
\DeclareMathOperator{\pr}{pr}

\DeclareMathOperator{\rE}{\mathrm{E}}
\DeclareMathOperator{\rH}{\mathrm{H}}

\DeclareMathOperator{\tr}{tr}

\begin{document}

\title[{\tiny   Mod $p$ Galois representations  and
cohomology of Hilbert  modular  varieties}]{Galois 
representations modulo $p$ and
cohomology of Hilbert modular varieties}

\author{Mladen Dimitrov}
\date{\today}

\begin{abstract} 
The aim of this paper is to extend some  arithmetic results 
on elliptic modular forms to the case of Hilbert modular  forms.
Among these results let's mention :

$\relbar$ the control of the image of the Galois representation modulo $p$  
\cite{serre2}\cite{ribet3},

$\relbar$ Hida's congruence criterion outside an explicit set of primes $p$  
\cite{hida63}, 

$\relbar$ the freeness of the integral cohomology of the Hilbert modular 
variety  over certain local components of the Hecke algebra and 
the Gorenstein property of these local algebras \cite{Ma}\cite{FaJo}.

\medskip 
We study the arithmetic of the  Hilbert modular  forms  by studying 
their modulo $p$ Galois representations and our main tool is 
the action of the inertia groups at the primes above $p$.
In order to determine  this action,  we compute the 
Hodge-Tate  (resp. the Fontaine-Laffaille) weights of the 
$p$-adic (resp. the modulo $p$) {\'e}tale  cohomology of the 
 Hilbert modular variety. The cohomological part of our paper is
inspired  by the work of Mokrane,  Polo and  Tilouine \cite{MoTi,PoTi} 
on the cohomology of the Siegel modular varieties and  builds upon 
the geometric constructions  of \cite{dimdg,dimtildg}.

\end{abstract}

\maketitle

\setcounter{tocdepth}{1}
\tableofcontents
\setcounter{tocdepth}{3}

\newpage
\section*{Introduction}

Let $F$ be a totally real number field of degree $d$,  ring of integers 
$\go$ and different $\gd$. Denote by   $\widetilde{F}$   
the Galois  closure of  $F$ in $\overline{\Q}$ and by  $J_F$ 
 the set of   all embeddings of $F$ into  $\overline{\Q}\subset \C$.
 
We fix an ideal $\gn\!\subset\! \go$ and we put  $\Delta= \Nm(\gn\gd)$. 

For a weight  $k=\sum_{\tau\in J_F}k_\tau\tau\in\Z[J_F]$ as in Def.\ref{weight} 
we put  $k_0=\max\{k_\tau |\tau \in J_F\}$. If $\psi$ is 
a Hecke character  of $F$  of conductor dividing $\gn$ and  type 
$2-k_0$ at infinity, we denote  by  $ S_k(\gn,\psi)$ the corresponding 
space of  Hilbert modular cuspforms (see  Def.\ref{space-hmf}).

Let  $f\!\in \!S_k(\gn,\psi)$ be  a  newform (that is a primitive
 normalized eigenform). For all  ideals $\ga\!\subset\! \go$,
 we denote by  $c(f,\ga)$ the eigenvalue of  the standard Hecke
operator $T_{\ga}$ on $f$. 

Let $p$ be a prime number and let 
$\iota_p:\overline{\Q}\hookrightarrow \overline{\Q}_p$ be an embedding.

Denote by   $E$  a sufficiently large  
$p$-adic field,   of ring of integers $\cO$, maximal ideal $\cP$ and
residue field $\kappa$.

\subsection{Galois image results.} 
The  absolute Galois  group of a field $L$ is denoted  by $\G_L$.
  By results of  Taylor \cite{Ta} and   Blasius-Rogawski  \cite{BR} 
there exists a continuous representation 
$\rho=\rho_{f,p}:\G_F \rightarrow \GL_2(E)$
 which is absolutely irreducible, totally odd, 
unramified outside $\gn p$ and 
such that for each prime ideal $v$ of $\go$,
not dividing $p\gn$, we have:
\begin{equation*}
\tr(\rho(\Frob_v))=\iota_p(c(f,v)), \enspace 
\det(\rho(\Frob_v))=\iota_p(\psi(v))\Nm(v),
\end{equation*}

\noindent where  $\Frob_v$ denotes a geometric Frobenius at $v$.

 By taking a Galois stable $\cO$-lattice, we define 
$\overline{\rho}=\rho\mod{\cP}:\G_F
\rightarrow \GL_2(\kappa)$,  whose  semi-simplification is   
independent of the particular choice of a lattice.

The following proposition is a generalization  to the Hilbert modular case,
of results of  Serre \cite{serre2} and  Ribet \cite{ribet3}
on elliptic modular forms (see Prop.\ref{Irr}, Prop.\ref{LI} and 
Prop.\ref{LI-Ind}).

\begin{prop} 
{\rm(i)} For all but finitely many primes $p$,

\medskip
$\mathrm{\bf (Irr_{\overline{\rho}})}$ $\overline{\rho}=\overline{\rho}_{f,p}$ is absolutely
 irreducible.

\smallskip
{\rm(ii)} If $f$ is not a theta series, then for 
all but finitely many primes $p$,

\medskip
$\mathrm{\bf (LI_{\overline{\rho}})}$ there exists a power $q$ of $p$ 
such that $\SL_2(\F_q)\!\subset\! 
\im(\overline{\rho})
\!\subset\! \kappa^\times \GL_2(\F_q)$.

\smallskip 
{\rm(iii)} Assume that $f$ is not a twist by a character of any of 
 its $d$ internal conjugates   and  is not a theta series. Then  
for  all but finitely many primes $p$, 
there exist a power $q$ of $p$,  a partition $J_F=\coprod_{i\in I} J_F^i $ 
and  $\sigma_{i,\tau}\!\in \!\Gal(\F_q/\F_p)$, $\tau\!\in \!J_F^i $,
such that $(\tau\neq \tau'\Rightarrow \sigma_{i,\tau}\neq \sigma_{i,\tau'})$
 and

$\mathrm{{\bf (LI}_{Ind \overline{\rho}}{\bf )}}$ 
$\Ind_F^{\Q}\overline{\rho}:\G_{\widehat{F}''}\rightarrow 
\SL_2(\F_q)^{J_F}$ factors as  a surjection 
$\G_{\widehat{F}''}\twoheadrightarrow 
\SL_2(\F_q)^{I}$ followed by the  map $(M_i)_{i\in I}\mapsto 
(M_i^{\sigma_{i,\tau}})_{i\in I,\tau\in J_F^i }$,
where  $\widehat{F}''$ denotes the compositum of   $\widetilde{F}$ and 
the fixed field of $(\Ind_F^{\Q}\overline{\rho})^{-1}(\SL_2(\F_q)^{J_F})$.
\end{prop}

\subsection{Cohomological results.} \label{coh-results}
Let  $Y_{/\Z[\frac{1}{\Delta}]}$ be  the Hilbert modular  variety
of level $K_1(\gn)$ (see \S\ref{hmv}). Consider the $p$-adic  {\'e}tale  cohomology
$\rH^\bullet(Y_{\overline{\Q}},\bV_n(\overline{\Q}_p))$, where 
 $\bV_n(\overline{\Q}_p)$  denotes the local system of weight 
$n=\sum_{\tau\in J_F}(k_\tau-2)\tau\in\n[J_F]$ (see \S\ref{MWC}). 
By a result of Brylinski and Labesse  \cite{BL} the subspace 
$W_f:=\underset{\ga\subset \go}{\bigcap} \ker(T_{\ga}-c(f,\ga))$ of 
 $\rH^d(Y_{\overline{\Q}},\bV_n(\overline{\Q}_p))$  is 
isomorphic,  as $\G_{\widetilde{F}}$-module and after
semi-simplification, to the tensor induced representation 
$\otimes\Ind_F^{\Q}\rho$.  

\newpage 

Assume that 

{\bf (I) }$p$ does not divide $\Delta$. 

\medskip 
Then $Y$ has  smooth toroidal compactifications  
over $\Z_p$ (see \cite{dimdg}). For each   $J\!\subset\! J_F$,  we put 
$|p(J)|=\sum_{\tau\in J}(k_0\!-\!m_\tau\!-\!1)+
\sum_{\tau\in J_F\bs J}m_\tau$, where $m_\tau=(k_0-k_\tau)/2\in\n$.
 By applying a  method 
of Chai and Faltings \cite{FaCh} Chap.VI,  one can  prove  
(see \cite{dimtildg} Thm.7.8, Cor.7.9)

\begin{theo} Assume that $p$ does not divide $\Delta$. Then 

{\rm(i)} the Galois representation 
 $\rH^j(Y_{\overline{\Q}},\bV_n(\overline{\Q}_p))$ is 
crystalline at $p$ and its Hodge-Tate weights 
belong to the set   $\{|p(J)|, J\!\subset\! J_F\enspace |J|\leq j \}$, and

\smallskip
{\rm(ii)} the   Hodge-Tate weights of  $W_f$ are given by the 
multiset  $\{|p(J)|, J\!\subset\! J_F\}$.
\end{theo}

For our main  arithmetic applications we need to establish a modulo $p$ version
 of the above theorem. This is achieved under the  following
 additional assumption  : 
 
\medskip 
{\bf (II) } $p\!-\!1>\underset{\tau\in J_F}{\sum}(k_\tau\!-\!1)$.

\smallskip
\noindent The integer  $\sum_{\tau\in J_F}(k_\tau\!-\!1)$ is equal to the difference  
$|p(J_F)|-|p(\varnothing)|$ 
between the largest and the smallest  Hodge-Tate weights of the 
cohomology of the Hilbert modular  variety.
We use  {\bf (I) } and {\bf (II) } in order to apply 
Fontaine-Laffaille's Theory \cite{FoLa},
as well as  Faltings' Comparison Theorem modulo $p$ \cite{Fa-jami}.
By  adapting to the case of Hilbert  modular varieties  some techniques  developed  by 
Mokrane, Polo and Tilouine \cite{MoTi, PoTi} for  Siegel modular
varieties, such as the construction of an  integral Bernstein-Gelfand-Gelfand 
complex for distribution algebras, we compute the  Fontaine-Laffaille 
weights of $\rH^\bullet(Y_{\overline{\Q}},\bV_n(\kappa))$
(see Thm.\ref{bgg-modp}).

\subsection{Arithmetic results.}
Consider the $\cO$-module of interior   cohomology $\rH_{!}^d(Y,\bV_n(\cO))'$, 
defined as the image of $\rH_{c}^d(Y,\bV_n(\cO))$ in $\rH^d(Y,\bV_n(E))$. 
Let $\T=\cO[T_{\ga},\ga\!\subset\! \go]$ 
be the full Hecke  algebra acting on it, and let $\T'\subset \T$ be the subalgebra generated by the Hecke operators
outside a finite set of places containing those dividing $\gn p$. 
Denote by $\gm$ the maximal ideal of $\T$  corresponding to  $f$ and $\iota_p$
and  put $\gm'=\gm\cap \T'$.

\begin{theo}\label{cohom}
Assume that the conditions  $\mathrm{{\bf (I)}}$ and 
$\mathrm{{\bf (II)}}$ from \S\ref{coh-results}  hold.

 {\rm (i)} If $\mathrm{\bf (Irr_{\overline{\rho}})}$ holds, 
 $d(p\!-\!1)>5\sum_{\tau\in J_F}(k_\tau\!-\!1)$ and 

\medskip
 $\mathrm{\bf (MW)}$  the middle weight  
$\frac{|p(J_F)|+|p(\varnothing)|}{2}=\frac{d(k_0-1)}{2}$ 
does not belong to   $\{|p(J)|,J\!\subset\! J_F\}$,

\vspace{2mm}
then  the local component  $\rH^\bullet_{\partial}(Y,\bV_n(\cO))_{\gm'}$ 
of the boundary  cohomology   vanishes, and  the  Poincar{\'e}  pairing  
$\rH_{!}^d(Y,\bV_n(\cO))'_{\gm'} \times
\rH_{!}^d(Y,\bV_n(\cO))'_{\gm'}\rightarrow \cO$ is a perfect duality.

{\rm (ii)} If $\mathrm{\bf (LI_{\Ind\overline{\rho}})}$ holds, then 
$\rH^\bullet(Y,\bV_n(\cO))_{\gm'}=\rH^d(Y,\bV_n(\cO))_{\gm'}$ is a free  
$\cO$-module  of  finite rank  and its Pontryagin dual is isomorphic to 
 $\rH^d(Y,\bV_n(E/\!\cO))_{\gm'}$.
\end{theo}

 The proof involves a ``local-global'' Galois argument. For the (i), 
we use lemma \ref{key2}(ii) and a  theorem of Pink \cite{Pi} on the 
{\'e}tale cohomology   of a local system  restricted to the boundary of
$Y$ (see Thm.\ref{boundary}). For the (ii),  we use lemma \ref{key}
and the  computation of the Fontaine-Laffaille weights of the
cohomology of Thm.\ref{bgg-modp} (see Thm.\ref{coh-loc}).

\medskip
Let $\Lambda^*(\Ad^0(f),s)$  be the imprimitive adjoint $L$-function  of $f$, 
completed  by its Euler factors at infinity and let 
$W(f)$ be the  complex constant from the functional equation 
of the  standard $L$-function  of $f$ (see {\S}\ref{L-ad}).
We denote  by  $\Omega_f^{\pm}\in \C^\times\!\!/\cO^\times$ any 
two complementary  periods
defined  by the  Eichler-Shimura-Harder isomorphism  (see {\S}\ref{periodes}).

\medskip
{\bf Theorem A} (Thm.\ref{theo-A})
{\it  Let  $f$ and  $p$ be such that {\bf (I)}, $\mathrm{\bf (Irr_{\overline{\rho}})}$ and {\bf (MW)} hold,
and $p-1>\max(1,\frac{5}{d})\sum_{\tau\in J_F}(k_\tau-1)$.  Assume
that  $\iota_p(\frac{W(f)\Lambda^*(\Ad^0(f),1)}{\Omega_f^+\Omega_f^-})\in\cP$.
Then there  exists another normalized eigenform  $g\in  S_k(\gn,\psi)$
such that $f\equiv g\!\pmod{\cP}$, in the sense that 
$c(f,\ga)\equiv c(g,\ga)\!\pmod{\cP}$ for  each  ideal $\ga\!\subset\! \go$.}

\medskip
The proof follows closely the original one given  by  Hida  \cite{hida63}
in the elliptic modular case, and uses Thm.\ref{cohom}(i)
as well as a formula of Shimura  relating 
$\Lambda^*(\Ad^0(f),1)$ to the Petersson inner product of $f$ (see 
(\ref{adjoint})). Let us note that Ghate  \cite{ghate1} has obtained a
very similar  result when the weight $k$ is parallel. A converse for 
Thm A is provided by the (ii) of the following

\medskip

{\bf Theorem B} (Thm.\ref{theo-B})
{\it Let  $f$ and  $p$ be such that   $\mathrm{{\bf (I)}}$, 
$\mathrm{{\bf (II)}}$ and $\mathrm{\bf (LI_{\Ind\overline{\rho}})}$  hold. Then

{\rm(i)} $\rH^\bullet(Y,\bV_n(\kappa))[\gm]=
\rH^d(Y,\bV_n(\kappa))[\gm]$ is a  $\kappa$-vector 
space  of dimension $2^d$. 

{\rm(ii)} $\rH^\bullet(Y,\bV_n(\cO))_{\gm}=
\rH^d(Y,\bV_n(\cO))_{\gm}$ 
is free of rank  $2^d$ over  $\T_{\gm}$.  

{\rm(iii)} $\T_{\gm}$ is Gorenstein.}

\medskip

By \cite{Ma} it is enough to prove (i), which is a consequence of 
  Thm.\ref{cohom}(ii) and the  $q$-expansion principle \S\ref{qdev}.

\medskip
The last theorem is due,  under milder assumptions, 
to Mazur \cite{Ma} for $F=\Q$ and $k=2$, and to 
Faltings and Jordan \cite{FaJo} for $F=\Q$. The Gorenstein
property is proved  by Diamond \cite{diamond2} when  $F$ is  quadratic and 
$k=(2,2)$, under the assumptions {\bf (I)}, {\bf (II)} and $\mathrm{\bf (Irr_{\overline{\rho}})}$.
We expect that Diamond's approach via intersection cohomology 
could be generalized in order to  prove the Gorenstein property 
of $\T_{\gm}$ under the assumptions 
{\bf (I)}, {\bf (II)} and $\mathrm{\bf (LI_{\overline{\rho}})}$ 
(see lemma \ref{key2}(i) and  remark \ref{key-diamond}).

\medskip
When $f$ is ordinary at $p$ (see Def.\ref{ord-autom}) we can replace 
 the assumptions {\bf (I)} and  {\bf (II)} of theorems A and B by the weaker 
assumptions that $p$ does not divide $\Nm(\gd)$ and that $k\pmod{p-1}$
satisfies {\bf (II)} (see Cor.\ref{thm-ord}). The proof uses  Hida's 
 families of $p$-adic ordinary   Hilbert modular forms. We prove an exact 
control theorem for the ordinary part of the cohomology of the 
 Hilbert modular variety, and give a new proof of Hida's  exact 
control theorem for the ordinary Hecke algebra (see Prop.\ref{control}). 

\medskip
 Theorems A and B relate  $\cO/\iota_p
(\frac{W(f)\Lambda^*(\Ad^0(f),1)}{\Omega_f^+\Omega_f^-})$ to the 
congruence module associated to the $\cO$-algebra homomorphism 
$\T\rightarrow \cO$, $T_{\ga}\mapsto \iota_p(c(f,\ga))$. In a subsequent paper
 \cite{dim-ihara} we relate these two
with the cardinality of the Selmer group of $\Ad^0(\rho)\otimes E/\cO$. 
An interesting question is whether $\Omega_f^{\pm}$ are the periods 
involved in the Bloch-Kato conjecture for $\Ad^0(f)$ (see
the work of Diamond, Flach and Guo \cite{dfg} for the elliptic modular case).

\subsection{Explicit results.} 
By a classical theorem of Dickson, the image in  $\PGL_2(\kappa)$
 of an irreducible subgroup of $\GL_2(\kappa)$  not 
satisfying $\mathrm{\bf (LI_{\overline{\rho}})}$, is either a  
 dihedral, a  tetrahedral, an octahedral or  an icosahedral group.  
In the next proposition we consider the later exceptional cases
for the image of $\overline{\rho}$ in $\PGL_2(\kappa)$.

\medskip
Denote by $\go^{\times}_+$ (resp. by $\go^{\times}_{\gn,1}$) be the group of 
 totally positive (resp. congruent to $1$ modulo $\gn$) units of $\go$. 

\newpage 

\begin{prop}
Assume   $p$ does not divide $\Delta$ and   $p>k_0$. 

{\rm(i)} Assume that $k$  is non-parallel.  If  for all  $J\!\subset\! J_F$, 
there exists an unit  $\epsilon\in\go^{\times}_+\cap\go^{\times}_{\gn,1}$,  such that  $p$ 
  does not divide  $\N_{\widetilde{F}\!/\!\Q}
\left(\underset{\tau\in J}{\prod}\tau(\epsilon)^{k_0-m_\tau-1}-
\underset{\tau\in J_F\!\bs\!J}{\prod}\tau(\epsilon)^{m_\tau}\right)
\neq 0$,  then  $\mathrm{\bf (Irr_{\overline{\rho}})}$ holds.

\smallskip
{\rm(ii)} If $d(p\!-\!1) >  5\underset{\tau\in J_F}{\sum}(k_\tau\!-\!1)$,
then the image of $\overline{\rho}$ in $\PGL_2(\kappa)$ is not 
a tetrahedral, an  octahedral nor  an  icosahedral group.

\smallskip
{\rm(iii)}  Assume  that for all $\tau\in J_F$, $p\neq 2k_\tau-1$ and the   following condition :

\medskip
\noindent {\bf (non-CM)}  for  each  quadratic extension 
 $K$  of  $F$ of discriminant dividing  $\gn$ and splitting  
all the primes  $\gp$ of $F$ above  $p$, one of the following
holds :  
\begin{enumerate}
\item 
$  K \text{ is CM and   there does not exist a 
Hecke character }   \varphi \text{ of }   K \text{ of   conductor dividing}\\
\gn\Delta_{K/F}^{-1} \text{ and   infinity type }
(m_\tau,k_0\!-\!1\!-\!m_\tau)_{\tau\in J_F}\text{, such that }
\rho\equiv \Ind_{K}^{F} \varphi \pmod{\cP}, $

\item 
$  K\text{  is not  CM and for all extensions }
 \widetilde{\tau}  \text{   of } \tau\in J_F \text{ to } K\text{,
 there exists a unit }
 \epsilon\in\gO^{\times},\\  \epsilon-1\in\gn  \text{ such that } 
 p \text{  does not divide }
 \displaystyle \N_{K/\Q}\left(\prod_{\tau\in J_F}
\widetilde{\tau}(\epsilon)^{m_{\tau}}\widetilde{\tau}(c(\epsilon))^{k_0-m_{\tau}-1}\!-\!1\right). $
\end{enumerate}

Then   the image   of  $\overline{\rho}$ in $\PGL_2(\kappa)$ 
 is not a  dihedral group. 

\smallskip
{\rm(iv)} Assume $\mathrm{\bf (LI_{\overline{\rho}})}$ and 
that $k$ is not induced from a weight for a strict subfield 
$F'$ of $F$. Assume moreover that for all $\tau,\tau'\in J_F$,
 $p\neq k_\tau+k_{\tau'}-1$. Then  $\mathrm{\bf (LI_{\Ind\overline{\rho}})}$ holds. 
\end{prop}

By the last proposition, we obtain the following  corollary to theorems A and B.

\begin{cor} Let   $\epsilon$ be any element of 
$\go^{\times}_+\cap\go^{\times}_{\gn,1}$.

{\rm(i)} Assume $d\!=\!2$ and $k\!=\!(k_0,k_0\!-\!2m_1)$,  with  $m_1\neq 0$. 
If $p\nmid \Delta 
\Nm((\epsilon^{m_1}\!-\!1)(\epsilon^{k_0\!-\!m_1\!-\!1}\!-\!1))$ and 
$p\!-\!1>4(k_0\!-\!m_1\!-\!1)$, then theorem A holds.  If additionally  
we have the  {\bf (non-CM)}  condition then   theorem B  also holds.

{\rm(ii)} Assume $d\!=\!3$, $\id\!\neq\!\tau\in J_F$ and 
$k\!=\!(k_0,k_0-2m_1,k_0-2m_2)$, 
with $0\leq m_1\leq m_2\neq 0$. 
If $p\nmid \Delta \N_{\widetilde{F}\!/\!\Q}(
(\tau(\epsilon)^{m_1}-\epsilon^{-m_2})
(\tau(\epsilon)^{m_1}-\epsilon^{m_2+1-k_0})
(\tau(\epsilon)^{m_1+1-k_0}-\epsilon^{m_2})
(\tau(\epsilon)^{k_0-m_1-1}-\epsilon^{m_2+1-k_0})
)$ and  $p-1>\frac{5}{3}(3k_0-2m_1-2m_2-3)$, 
then theorem A holds.  If additionally  we
have $p\neq 2k_0-1$ and the  {\bf (non-CM)}  condition then   theorem B  also holds.
\end{cor}

\subsection{Acknowledgements.}
 I would like to thank J. Bella{\"\i}che, D. Blasius, G. Chenevier,
 L. Dieulefait, E. Ghate, M. Kisin, 
V. Lafforgue, A. Mokrane, E. Urban, J. Wildeshaus and J.-P. Wintenberger for helpful 
conversations.
This article was completed during my visit at UCLA on an invitation by  H. Hida. 
I would like to thank him heartily for his  hospitality and for many inspiring 
discussions. 
I am grateful to G. Kings and R. Taylor for their interesting comments 
on an earlier version of this paper that have much improved it.
Finally, I would like to thank J. Tilouine who suggested me to study this problem 
and supported me during the  preparation of this article.

\newpage

\section{Hilbert modular forms and varieties.}\label{hmfv}

We define the algebraic groups  $D_{\!/\!\Q}=\Res^{F}_{\Q} \Gm$,
 $G_{\!/\!\Q}=\Res^{F}_{\Q} \GL_2$ and   $G^*_{/\Q}=G\times_{D}\Gm$,  
where  the fiber product is relative to  the reduced norm map $\nu: G \rightarrow D$. The standard 
Borel subgroup  of $G$, its unipotent  radical  and  its standard maximal 
torus    are denoted  by $B$, $U$ and  $T$, respectively. 
We identify $D \times D$ with $T$, by $(u,\epsilon)\mapsto 
\begin{pmatrix} u\epsilon & 0 \\ 0 & u^{-1} \end{pmatrix}$.

\subsection{Analytic Hilbert modular varieties.}

Let $D(\R)_+$ (resp. $G(\R)_+$) be the identity component  of  
$D(\R)=(F\otimes \R)^\times$  (resp. of $G(\R)$). The group $G(\R)_+$ 
acts  by linear fractional transformations  on     the space  
$\mathfrak{H}_F=\{z\in F\otimes\C \enspace |\enspace \im(z)\in 
D(\R)_+\}$. We have $\mathfrak{H}_F\cong \mathfrak{H}^{J_F}$, where 
$\mathfrak{H}=\{z\in \C\enspace |\enspace \im(z)>0\} $ is  the 
 Poincar{\'e}'s upper half-plane (the isomorphism being  given by 
$\xi\otimes z \mapsto (\tau(\xi)z)_{\tau\in J_F}$, for $\xi\in F$, $z\in \C$).
We consider the unique group action  of $G(\R)$ on  the space  
$\mathfrak{H}_F$ extending  the action of $G(\R)_+$  and such that, 
on each copy of $\mathfrak{H}$ 
the element  $\begin{pmatrix} -1 & 0\\ 0 & 1\end{pmatrix}$ acts 
 by $z\mapsto -\overline{z}$.
 We put  $\underline{i}=(\sqrt{-1},...,\sqrt{-1})\in \mathfrak{H}_F$,  
$K_{\infty}^+=\Stab_{G(\R)_+}(\underline{i})=\SO_2(F\otimes\R)D(\R)$ and 
$K_{\infty}=\Stab_{G(\R)}(\underline{i})=\mathrm{O}_2(F\otimes \R)D(\R)$.

We denote by $\widehat{\Z}=\prod_l \Z_l$ the profinite completion  of $\Z$ 
and we put  $\widehat{\go}=\widehat{\Z}\otimes \go=
\prod_v \go_v$, where $v$ runs over all the finite places  of $F$.
Let  $\A$ (resp. $\A_f$) be  the ring of the  ad{\`e}les 
(resp. of the finite ad{\`e}les) of  $\Q$. We 
consider the following open compact subgroup of $G(\A_f)$ :
$$K_1(\gn)=\left\{\begin{pmatrix}a& b \\ c & d\end{pmatrix}\in
G(\widehat{\Z}) | d-1\in \gn, c\in \gn \right\}.$$

 The ad{\'e}lic Hilbert modular variety of level 
$K_1(\gn)$ is defined as 
$$Y^{\an}=Y_1(\gn)^{\an}=
G(\Q)\bs G(\A)/K_1(\gn)K_{\infty}^+.$$

By the Strong Approximation Theorem, the 
 connected components of $Y^{\an}$ are indexed by the narrow ideal class group 
$\Cl_F^{+}\!=\!D(\A)/D(\Q)D(\widehat{\Z})D(\R)_+$
of $F$. For each fractional ideal $\gc$ of $F$, we put $\gc^*=\gc^{-1}\gd^{-1}$.
We define the following congruence subgroup of $G(\Q)$ :
$$\Gamma_1(\gc,\gn)=\Big{\{}
\begin{pmatrix}a &b \\c &d\end{pmatrix} \in G(\Q)\cap
\begin{pmatrix} \go & \gc^* \\ 
\mathfrak{cdn} & \go\end{pmatrix}
\enspace \Big{|} \enspace ad-bc\in \go_+^\times, \enspace  d\equiv 1 \pmod{\gn} \Big{\}}.$$

Put $M^{\an}=M_1(\gc,\gn)^{\an}=
\Gamma_1(\gc,\gn)\bs  \mathfrak{H}_F$.
Then we have $Y_1(\gn)^{\an}\simeq\coprod_{i=1}^{h^+} 
M_1(\gc_i,\gn)^{\an}$,  where the ideals $\gc_i$, $1\leq i\leq h^+$, 
form  a set of representatives of  $\Cl_F^{+}$.

\medskip
Put $\mathfrak{H}_F^*= \mathfrak{H}_F\coprod \p^1(F)$. The minimal 
compactification $M^{*\an}$  of  $M^{\an}$  is defined as $M^{*\an}=\Gamma\bs \mathfrak{H}_F^*$. It is an analytic  normal projective  space whose
boundary $M^{*\an}\bs M^{\an}$ is a finite union of closed points,
called the {\it cusps} of $M^{\an}$.  

\medskip
The same way, by replacing $G$ by $G^*$, we define   
$\Gamma_1^1(\gc,\gn)$, $M^{1,\an}=M_1^1(\gc,\gn)^{\an}$ and $M^{1*,\an}$.


\subsection{Analytic Hilbert modular forms.}
\label{fmh-hecke}

For the definition of the $\C$-vector space of  Hilbert modular forms we 
follow   \cite{hida128}.

\begin{defin}\label{weight} An element 
$k=\sum_{\tau\in J_F}k_\tau\tau\in\Z[J_F]$ is called a  weight. We
always assume that the  $k_\tau$ have the same parity  and 
are all $\geq 2$.  We put 
$k_0=\max\{k_\tau |\tau \in J_F\}$, $n_0=k_0-2$, $t=\sum_{\tau\in J_F}\tau$,
$n=\sum_{\tau\in J_F}n_\tau\tau=k-2t$ and 
$m=\sum_{\tau\in J_F}m_\tau\tau=(k_0t-k)/2$.
\end{defin}

 For $z\in \mathfrak{H}_F$, $\gamma=\begin{pmatrix}a&b\\c&d\end{pmatrix}$
we put $j_J(\gamma,z)=c\cdot z^J+d \in D(\C)$, 
where $z^J_\tau=\begin{cases} z_\tau\textrm{, } \tau\in J,\\  
\overline{z}_\tau\textrm{, } \tau\in J_F\bs J.\end{cases}$ 


\begin{defin}  \label{hmaf} The space $G_{k,J}(K_1(\gn))$ of  ad{\'e}lic Hilbert modular 
 forms of weight $k$, level $K_1(\gn)$ and type $J\!\subset\! J_F$ 
at infinity, is  the $\C$-vector space of the   functions 
$g:G(\A) \rightarrow \C$  satisfying the    followings three
conditions : 

(i) $g(axy)=g(x)$ for all $a\in G(\Q)$, $y\in K_1(\gn)$ and 
$x\in G(\A)$.

(ii) $g(x\gamma)= \nu(\gamma)^{k+m-t}j_J(\gamma,\underline{i})^{-k}
g(x)$, for all  $\gamma\in K_{\infty}^+$ and $x\in G(\A)$.

\medskip
 For all $x\in G(\A_f)$ define  $g_x:\mathfrak{H}_F\rightarrow \C$,   by 
$z\mapsto \nu(\gamma)^{t -k-m} j_J(\gamma,\underline{i})^{k}g(x\gamma),$
\noindent where $\gamma\in G(\R)_+$ is such that  $z=\gamma\cdot \underline{i}$. 
By (ii) $g_x$ does not depend on the  particular choice  of  $\gamma$.

\medskip
(iii)  $g_x$ is holomorphic at  $z_\tau$, for $\tau\in J$, and anti-holomorphic 
at  $z_\tau$, for $\tau\in J_F\bs J$ (when $F=\Q$ an extra condition of 
holomorphy at cusps is needed).

The space $S_{k,J}(K_1(\gn))$ of  ad{\'e}lic Hilbert modular  cuspforms
is the subspace of $G_{k,J}(K_1(\gn))$ consisting of functions satisfying 
the following additional condition : 

(iv) $\int_{U(\Q)\bs U(\A)}g(ux)du=0$, for all $x\in G(\A)$ and
 all additive Haar measures    $du$.
\end{defin}

The conditions  (i) and (ii) of the above  definition imply that for all 
$g\in G_{k,J}(K_1(\gn))$ there exists  a  Hecke character
  $\psi$  of  $F$ of conductor  dividing $\gn$ and  of type 
$-n_0t $ at infinity, such that   for all $x\in G(\A)$ and for all 
$z\in D(\Q)D(\widehat{\Z})D(\R)$,  we have $g(zx)=\psi(z)^{-1}g(x)$. 

\begin{defin}  \label{space-hmf}
Let $\psi$ be a Hecke character  of  $F$ of conductor 
dividing $\gn$ and  of type  $-n_0t $ at infinity.  The space  
$S_{k,J}(\gn,\psi)$ (resp. $G_{k,J}(\gn,\psi)$) is defined as the subspace 
of  $S_{k,J}(K_1(\gn))$ (resp. $G_{k,J}(K_1(\gn))$) of elements  $g$ satisfying
$g(zx)=\psi(z)^{-1}g(x)$, for all $x\in G(\A)$ and for all $z\in D(\A)$.
When $J=J_F$ this  space is denoted by $S_k(\gn,\psi)$ (resp. by $G_k(\gn,\psi)$).
\end{defin}

As the characters of the ideal class  group
$\Cl_F= D(\A)/D(\Q)D(\widehat{\Z})D(\R)$  of $F$ form a basis
of the   complex valued functions on this set, we have :
\begin{equation} \label{decomp-fmh}
G_{k,J}(K_1(\gn))=\bigoplus_{\psi} G_{k,J}(\gn,\psi),\enspace\enspace
S_{k,J}(K_1(\gn))=\bigoplus_{\psi} S_{k,J}(\gn,\psi)
\end{equation}
where $\psi$ runs over  the   Hecke characters of  $F$, 
of conductor dividing $\gn$ and infinity type $-n_0t$.
Let $\Gamma$ be a congruence subgroup of $G(\Q)$. We recall the
classical definition :

\begin{defin} \label{def-hmf} The space $G_{k,J}(\Gamma;\C)$ of Hilbert modular 
  forms of weight $k$, level $\Gamma$ and type $J\!\subset\! J_F$ 
at infinity, is  the $\C$-vector space of the   functions 
$g:\mathfrak{H}_F \rightarrow \C$  which are holomorphic at  
$z_\tau$, for $\tau\in J$, and anti-holomorphic at  $z_\tau$, for $\tau\in
 J_F\bs J$, and such that for every $\gamma\in \Gamma$
we have $g(\gamma(z))=\nu(\gamma)^{t-k-m}j_J(\gamma,z)^{-k}g(z)$. 

The space $S_{k,J}(\Gamma;\C)$ of  Hilbert modular 
 cuspforms is the subspace of $G_{k,J}(\Gamma;\C)$, consisting of 
functions vanishing at all cusps.
\end{defin}

Put $x_i=\begin{pmatrix}\eta_i & 0 \\0 &1\end{pmatrix}$, where  $\eta_i$ is 
the id{\`e}le associated to the ideal $\gc_i$, $1\leq i\leq h^+$.  The map 
$g\mapsto (g_{x_i})_{1\leq i\leq h^+}$ (see Def.\ref{hmaf}) induces isomorphisms :
\begin{equation}\label{ad-cl}
G_{k,J}(K_1(\gn))\simeq\bigoplus_{1\leq i\leq h^+} G_{k,J}(\Gamma_1(\gc_i,\gn);\C), 
\enspace\enspace 
S_{k,J}(K_1(\gn))\simeq\bigoplus_{1\leq i\leq h^+} S_{k,J}(\Gamma_1(\gc_i,\gn);\C).
\end{equation}

 Let  $d\mu(z)=\underset{\tau\in J_F}{\prod}y_\tau^{-2}dx_\tau dy_\tau$
be the standard  Haar  measure on  $\mathfrak{H}_F$.

\begin{defin}\label{peterson}
The Petersson inner product of two cuspforms $g,h\in S_{k,J}(K_1(\gn))$
is given by the formula 
$$(g,h)_{\gn}=\sum_{i=1}^{h^+} \int_{\Gamma_1(\gc_i,\gn)\bs\mathfrak{H}_F}
\overline{g_i(z)} h_i(z) y^k d\mu(z),$$
where $(g_i)_{1\leq i\leq h^+}$ (resp. $(h_i)_{1\leq i\leq h^+}$)
is the image of $g$ (resp. $h$) under the isomorphism (\ref{ad-cl}).
\end{defin}

\subsection{Hilbert-Blumenthal abelian varieties.}\label{hbav}
A sheaf over a scheme $S$ which is locally free of rank one over 
 $\go\otimes \cO_S$, is called an {\it invertible $\go$-bundle} on $S$.

\begin{defin} A Hilbert-Blumenthal abelian variety (HBAV) over a 
$\Z[\frac{1}{\Nm(\gd)}]$-scheme $S$
is an abelian scheme $\pi:A\rightarrow S$ of relative dimension $d$ 
together with an injection  $\iota:\go\hookrightarrow \End(A/S)$, 
such that  $\underline{\omega}_{A/S}:=\pi_*
\Omega^1_{A/S}$ is  an invertible $\go$-bundle on $S$.
\end{defin}

Let $\gc$ be a fractional ideal of $F$ and  $\gc_+$ be the cone of 
totally positive elements in  $\gc$.   Given a HBAV $A/S$, the 
functor assigning  to a $S$-scheme $X$ the set  $A(X)\otimes_{\go}\gc$
is representable by  another HBAV, denoted by $A\otimes_{\go}\gc$.
Then $\iota$ yields $\gc\hookrightarrow  \Hom_{\go}(A,A\otimes_{\go}\gc)$. 
The dual of a HBAV $A$ is denoted by $A^t$.

\begin{defin}
{\rm (i)} A $\gc$-polarization on a HBAV  $A/S$ is an $\go$-linear 
isomorphism $\lambda:A\otimes_{\go}\gc \overset{\sim}{\longrightarrow} A^t$, 
such that under the  induced isomorphism $\Hom_{\go}(A,A\otimes_{\go}\gc)
\cong \Hom_{\go}(A,A^t)$ elements of $\gc$ (resp. $\gc_+$)
correspond exactly to symmetric elements (resp. polarizations).

{\rm (ii)} A $\gc$-polarization class $\overline{\lambda}$
is an orbit of $\gc$-polarizations under  $\go_+^\times$. 
\end{defin}

 Let  $(\Gm\otimes \mathfrak{d}^{-1})[\gn]$
be the reduced  subscheme of $\Gm\otimes \mathfrak{d}^{-1}$,  defined
as the intersection 
of the kernels of multiplications  by elements of $\gn$. Its  Cartier dual 
is isomorphic to the   finite    group scheme $\go\!/\!\gn$.

\begin{defin}
 A $\mu_{\gn}$-level structure on a HBAV $A/S$
is an $\go$-linear closed immersion $\alpha:(\Gm\otimes \mathfrak{d}^{-1})[\gn]
\hookrightarrow A$ of  group schemes over $S$.
\end{defin}

\subsection{Hilbert modular varieties.}\label{hmv}
We consider the contravariant functor $\underline{\cM}^1$
(resp. $\underline{\cM}$) from the category of 
$\Z[\frac{1}{\Delta}]$-schemes to the category of sets, 
assigning  to a scheme $S$  the set of isomorphism classes
of triples $(A,\lambda,\alpha)$ (resp. 
 $(A,\overline{\lambda},\alpha)$) where 
$A$ is a HBAV over $S$, endowed with a $\gc$-polarization  $\lambda$
(resp. a  $\gc$-polarization  class $\overline{\lambda}$)  and 
a  $\mu_{\gn}$-level structure $\alpha$. 
Assume the following condition :

\medskip
{\bf (NT)  } $\gn$  does not divide $2$, nor $3$, nor $\Nm(\mathfrak{d})$.

\smallskip
Then $\Gamma_1(\gc,\gn)$ is torsion free, and 
 the functor $\underline{\cM}^1$ is representable 
by a quasi-projective, smooth,  geometrically connected 
$\Z[\frac{1}{\Delta}]$-scheme $M^1=M^1_1(\gc,\gn)$, 
endowed with  an universal HBAV $\pi:\cA\rightarrow M^1$. 
By definition,  the sheaf $\underline{\omega}_{\cA/M^1}=\pi_*
\Omega^1_{\cA/M^1}$ is an invertible $\go$-bundle on $M^1$. 
Consider the first de Rham cohomology sheaf 
$\cH^1_{\dR}(\cA/M^1)=R^1\pi_*\Omega^\bullet_{\cA/M^1}$ on $M^1$.  
The Hodge filtration yields an exact sequence :
$$0 \rightarrow \underline{\omega}_{\cA/M^1} \rightarrow 
\cH^1_{\dR}(\cA/M^1)  \rightarrow 
\underline{\omega}_{\cA/M^1}^{\vee}\otimes\mathfrak{cd}^{-1}  \rightarrow 0.$$
Therefore $\cH^1_{\dR}(\cA/M^1)$ 
is locally free of rank two over $\go\otimes \cO_{M^1}$.

The  functor $\underline{\cM}$ admits a coarse moduli space 
$M=M_1(\gc,\gn)$, which is a  quasi-projective, smooth, geometrically 
connected  $\Z[\frac{1}{\Delta}]$-scheme. The finite group 
$\go_+^\times/\go^{\times 2}_{\gn,1}$ acts properly and discontinuously on
  $M^1$ by  $[\epsilon]:(A,\iota,\lambda,\alpha)/S\mapsto 
(A,\iota,\epsilon\lambda,\alpha)/S$ and the quotient is given by $M$. 
This group acts also on $\underline{\omega}_{\cA/M^1}$ 
 and  on $\cH^1_{\dR}(\cA/M^1)$  by acting on the de Rham 
complex  $\Omega^\bullet_{\cA/M^1}$ ($[\epsilon]$ acts on  
$\underline{\omega}_{\cA/M^1}$ by  $\epsilon^{-1/2}[\epsilon]^*$).

These actions are defined over the integer ring of the number field
$F(\epsilon^{1/2}, \epsilon\in \go_+^\times)$. 

Let   $\go'$ be the integer ring of 
$\widetilde{F}(\epsilon^{1/2}, \epsilon\in \go_+^\times)$. For 
every $\Z[\frac{1}{\Delta}]$-scheme $X$ we put 
$$\textstyle X'=X\times \Spec(\go'[\frac{1}{\Delta}]).$$  

The sheaf of $\go_+^\times/\go^{\times 2}_{\gn,1}$- invariants  of 
 $\underline{\omega}_{\cA/M^1}$ (resp. of $\cH^1_{\dR}(\cA/M^1)$) 
 is locally free of rank  one (resp. two)  over $\go\otimes \cO_{M'}$, and 
is denoted  by $\underline{\omega}$ (resp. $\cH^1_{\dR}$). 

We put $Y=Y_1(\gn)=\coprod_{i=1}^{h^+} M_1(\gc_i,\gn)$ and
$Y^1=Y_1^1(\gn)=\coprod_{i=1}^{h^+} M_1^1(\gc_i,\gn)$,
  where the ideals $\gc_i$, $1\leq i\leq h^+$, form  
 a set of representatives of  $\Cl_F^{+}$.
\subsection{Geometric Hilbert modular forms.}

Under the action of $\go$, the invertible $\go$-bundle  $\underline{\omega}$
on $M'$ decomposes as a direct sum of line bundles 
$\underline{\omega}_{\tau}$, $\tau\in J_F$.  For every 
$k=\sum_\tau k_\tau\tau \in \Z[J_F]$ we define the line bundle 
$\underline{\omega}^k=\underset{\tau}{\otimes}
\underline{\omega}_{\tau}^{\otimes k_\tau}$ on $M'$. 

One should be careful to observe, that the global section of 
$\underline{\omega}^k$ on $M^{\an}$ are given by the cocycle
$\gamma\mapsto \nu(\gamma)^{-k/2}j(\gamma,z)^k$, meanwhile we are interested
in finding a geometric interpretation of the cocycle 
$\gamma\mapsto \nu(\gamma)^{t-k-m}j(\gamma,z)^k$, used in Def.\ref{def-hmf}.

\medskip
The universal polarization class $\overline{\lambda}$ 
endows $\cH^1_{\dR}$ with  a perfect symplectic $\go$-linear pairing. 
Consider the invertible $\go$-bundle  $\underline{\nu}:=
\wedge^2_{\go\otimes \cO_{M'}}\cH^1_{\dR}$ on $M'$. 
Note that $(k+m-t)-\frac{k}{2}=\frac{n_0}{2}t$. 

\begin{defin} Let $R$ be an $\go'[\frac{1}{\Delta}]$-algebra.
A Hilbert modular forms of weight $k$, level $\Gamma$ and
coefficients in $R$, is a global section of $\underline{\omega}^k\otimes 
\underline{\nu}^{-n_0t/2}$ over  $M\times_{\Spec(\Z[\frac{1}{\Delta}])}\Spec(R)$.
 We denote by $G_k(\Gamma;R)= \rH^0(M\times_{\Spec(\Z[\frac{1}{\Delta}])} 
\Spec(R),\underline{\omega}^k\otimes \underline{\nu}^{-n_0t/2})$ the 
$R$-module of these Hilbert modular forms.
\end{defin}

\subsection{Toroidal  compactifications.}\label{toroidal}

The toroidal compactifications of the moduli space of $\gc$-polarized 
HBAV with {\it principal} level structure have been constructed
by Rapoport \cite{rapoport}. Several modifications  need to 
be made in  order  to treat the case of $\mu_{\gn}$-level structure. These 
are described  in \cite{dimdg}Thm.7.2.

Let  $\Sigma$ be a smooth $\Gamma_1^1(\gc,\gn)$-admissible collection of fans 
(see \cite{dimdg}Def.7.1).
 Then,  there exists an open immersion of $M^1$ into a proper and  smooth
$\Z[\frac{1}{\Delta}]$-scheme $\overline{M^1}=M^1_{\Sigma}$, called the 
 toroidal  compactification of $M^1$ with respect to   $\Sigma$. 
 The universal HBAV   $\pi:\cA\rightarrow M^1$ extends uniquely to
a semi-abelian scheme $\overline{\pi}:\mathfrak{G} \rightarrow \overline{M^1}$.
The group scheme $\mathfrak{G} $ is endowed with an action of 
$\go$ and its restriction to $\overline{M^1}\bs M^1$ is a torus.
Moreover, the sheaf $\underline{\omega}_{\mathfrak{G}/\overline{M^1}}$ of
$\mathfrak{G}$-invariants sections of 
$\overline{\pi}_*\Omega^1_{\mathfrak{G}/\overline{M^1}}$ is an 
invertible $\go$-bundle on $\overline{M^1}$,  extending
 $\underline{\omega}_{\cA/M^1}$.

The scheme  $\overline{M^1}\bs M^1$ is a divisor with normal crossings, 
and the formal completion of $\overline{M^1}$ along this divisor
can be  completely determined in terms of $\Sigma$ (see \cite{dimdg}Thm.7.2).
For the sake of simplicity, we will only describe the 
 completion of $\overline{M^1}$ along the connected component 
of $\overline{M^1}\bs M^1$ corresponding to the 
standard cusp at $\infty$. Let $\Sigma^\infty\in \Sigma$ be the 
 fan corresponding to the   cusp at $\infty$. It is a  complete, smooth 
fan of $\gc_+^*\cup \{0\}$, stable by the action of $\go^{\times 2}_{\gn,1}$, 
and containing a finite number of cones modulo this action. 
Put $R_\infty=\Z[q^\xi,\xi\in\gc]$ and 
$S_\infty=\Spec(R_\infty)=\Gm\otimes \gc^*$. Associated to the 
fan $\Sigma^\infty$, there is a  toroidal embedding  
$S_\infty\hookrightarrow S_{\Sigma^\infty}$ (it is obtained 
by gluing the affine toric  embeddings $S_\infty\hookrightarrow 
S_{\infty,\sigma}=\Spec(\Z[q^\xi,\xi\in\gc\cap\check{\sigma}])$ for 
$\sigma\in \Sigma^\infty$). 
Let $S_{\Sigma^\infty}^{\wedge}$ be the formal completion of
 $S_{\Sigma^\infty}$ along $S_{\Sigma^\infty}\bs S_\infty$. 
By construction, the formal completion of $\overline{M^1}$ 
along the connected component  of $\overline{M^1}\bs M^1$ corresponding to the 
standard cusp at $\infty$, is isomorphic to $S_{\Sigma^\infty}^{\wedge}/
\go^{\times 2}_{\gn,1}$. 

Assume that  $\Sigma$ is  $\Gamma_1(\gc,\gn)$-admissible (for the cusp at $\infty$, 
it means that $\Sigma^\infty$ is stable under the action of 
$\go_+^{\times}$). Then  the finite group $\go_+^\times/\go^{\times 2}_{\gn,1}$ 
acts properly and discontinuously  on $\overline{M^1}$, and  the quotient 
$\overline{M}=M_{\Sigma}$ is a proper and smooth $\Z[\frac{1}{\Delta}]$-scheme,
containing $M$ is an open subscheme. Again by construction, the formal
 completion of $\overline{M}$  along the connected component  of 
$\overline{M}\bs M$ corresponding to the  standard cusp at $\infty$, 
is isomorphic to $S_{\Sigma^\infty}^{\wedge}/ \go_{+}^{\times}$.

The invertible $\go$-bundle  $\underline{\omega}_{\mathfrak{G}/\overline{M^1}}$
on $\overline{M^1}$ descends to an invertible $\go$-bundle on $\overline{M}{}'$, 
extending  $\underline{\omega}$. We still denote this extension by 
$\underline{\omega}$. For each $k\in \Z[J_F]$ this gives us an extension of
 $\underline{\omega}^k$ to a line bundle on $\overline{M}{}'$,  still 
denoted by  $\underline{\omega}^k$.

\subsection{$q$-expansion and Koecher Principles.}\label{qdev}

The Koecher Principle states (see \cite{dimdg}Thm.8.3)
\begin{equation}\label{koecher}
\rH^0(M\times\Spec(R), \underline{\omega}^k\otimes \underline{\nu}^{-n_0t/2})=
\rH^0(\overline{M}\times\Spec(R),
\underline{\omega}^k\otimes \underline{\nu}^{-n_0t/2})
\end{equation}

For  simplicity, we will only describe the $q$-expansion at the 
standard (unramified) cusp at $\infty$. For every $\sigma\in \Sigma^\infty$, 
and every $\go'[\frac{1}{\Delta}]$-algebra $R$, the pull-back of 
$\underline{\omega}$ to $S_{\Sigma^\infty}^{\wedge}\times\Spec(R)$ is 
canonically isomorphic to 
$\go\otimes \cO_{S_{\Sigma^\infty}^{\wedge}}\otimes R$.  Thus 
$$\hspace{-5mm}\rH^0(S_{\Sigma^\infty}^{\wedge}\times\Spec(R)/\go_+^\times,
\underline{\omega}^k\otimes \underline{\nu}^{-n_0t/2})\!=\!
\left\{\! \sum_{\xi\in\gc_+\cup\{0\}}\!\!\!\! a_\xi q^\xi\Big{|}
 a_\xi\!\in\! R, a_{u^2\epsilon\xi}=u^k\epsilon^{k+m-t} a_\xi, 
\forall (u,\epsilon)\! \in \! \go^{\times}_{\gn,1} \!\times\!\go_+^\times \right\}$$

By the above construction, to each $g\in G_k(\Gamma;R)$, 
we can associate an element $g_\infty= \sum_{\xi\in\gc_+\cup\{0\}}a_\xi(g)q^\xi$, 
 called the $q$-expansion of  $g$ at the cusp at $\infty$. 
The element $a_0(g)\in R$ is the value of  $g$ at the cusp at $\infty$.  

\begin{prop} Let $R$ be a $\go'[\frac{1}{\Delta}]$-algebra. 

{\rm (i)($q$-expansion Principle)}  $ G_k(\Gamma;R)\rightarrow 
 R[[q^\xi,\xi\in\gc_+\cup\{0\}]],  \enspace g\mapsto g_\infty$ is injective.

(ii) If there exists $g\in G_k(\Gamma;R)$, such that 
$a_0(g)\neq 0$, then  $\epsilon^{k+m-t}-1$ is a zero-divisor in $R$, 
for  all $\epsilon\in \go_+^\times $. 
\end{prop}

\subsection{The minimal compactification.}
There exist a projective, normal $\Z[\frac{1}{\Delta}]$-scheme
$M^{1*}$, containing $M^1$ as an open dense subscheme and 
such that the scheme $M^{1*}\bs M^1$ is  finite and {\'e}tale 
over  $\Z[\frac{1}{\Delta}]$.   Moreover, for each toroidal 
compactification $\overline{M^1}$  of $M^1$ there is a 
natural surjection $\overline{M^1}\rightarrow M^{1*}$, inducing 
the identity map on  $M^1$. 
The scheme $M^{1*}$ is called the  minimal compactification of $M^1$. 
The action of  $\go_+^\times/\go^{\times 2}_{\gn,1}$
on $M^1$ extends to an action on $M^{1*}$, and the 
minimal compactification $M^*$ of $M$ is defined as
the quotient for this action. In general 
 $M^{1*}\rightarrow M^*$ is not {\'e}tale.

We summarize the above discussion in the following commutative diagram :
$$\xymatrix@C=20pt@R=8pt{
\mathfrak{G}\ar[rr]^{\overline{\pi}} && \overline{M^1}\ar@{->>}[rr]
\ar@{->>}[rd] && \overline{M}\ar@{->>}[rd] & \\
              & &                &M^{1*} \ar@{-->>}[rr] &      & M^*\\
\cA \ar[rr]^{\pi}\ar@{^{(}->}[uu]   &&      M^1 \ar@{->>}[rr]\ar@{^{(}->}[uu]\ar@{^{(}->}[ru]    &&   M \ar@{^{(}->}[uu] \ar@{^{(}->}[ru]     & }$$

\subsection{Toroidal compactifications of the Kuga-Sato varieties.}\label{kugasato}
Let $s$ be a positive integer. Let  $\pi_s:\cA^s\rightarrow M^1$ be 
the $s$-fold fiber product of $\pi:\cA\rightarrow M^1$, and 
$(\overline{\pi})_s: \mathfrak{G}^s\rightarrow  \overline{M^1}$
be the $s$-fold fiber product of $\overline{\pi} :\mathfrak{G}
\rightarrow  \overline{M^1}$. 

Let $\widetilde{\Sigma}$ be a  $(\go\oplus\gc)\rtimes 
\Gamma_1^1(\gc,\gn)$-admissible, polarized, equidimensional,  smooth collection of fans, 
above the $\Gamma_1^1(\gc,\gn)$-admissible collection of fans $\Sigma$ of \S\ref{toroidal}. 
Using Faltings-Chai's method
\cite{FaCh}, the main result of \cite{dimtildg} Sect.6 is the following :
there exists an open immersion of a $\cA^s$ into a 
projective smooth $\Z[\frac{1}{\Delta}]$-scheme 
$\overline{\cA^s}=\cA^s_{\widetilde{\Sigma}}$, 
and a  proper, semi-stable homomorphism  
$\overline{\pi_s}:\overline{\cA^s}\rightarrow  \overline{M^1}$
extending $\pi_s:\cA^s\rightarrow M^1$, and such that 
$\overline{\cA^s}\bs\cA^s$ is a relative normal crossing divisor 
above $ \overline{M^1}\bs M^1$. Moreover, $\overline{\cA^s}$
contains $\mathfrak{G}^s$ as an open dense subscheme and 
$\mathfrak{G}^s$ acts on $\overline{\cA^s}$
extending the translation action of $\cA^s$ on itself.

The  sheaf  $\cH^1_{\ldR}(\overline{\cA}/\overline{M^1})=
R^1\overline{\pi_1}_*\Omega^\bullet_{\overline{\cA}/\overline{M^1}}
(\dlog\infty)$ is independent of the particular choice of $\widetilde{\Sigma}$
 above   $\Sigma$ and is endowed with a filtration :
$$0 \rightarrow \underline{\omega}_{\mathfrak{G}/\overline{M^1}} \rightarrow 
\cH^1_{\ldR}(\overline{\cA}/\overline{M^1})  \rightarrow 
\underline{\omega}_{\mathfrak{G}/\overline{M^1}}^{\vee}
\otimes\mathfrak{cd}^{-1}  \rightarrow 0.$$
It descends to a sheaf $\cH^1_{\ldR}$ on $\overline{M}$  which fits
in the following exact sequence : 
$$ 0 \rightarrow \underline{\omega}\rightarrow \cH^1_{\ldR}
 \rightarrow \underline{\omega}^{\vee}
\otimes\mathfrak{cd}^{-1}  \rightarrow 0.$$

\subsection{Hecke operators on modular forms.}\label{hecke} Let 
 $\Z[K_1(\gn)\bs G(\A_f)/K_1(\gn)]$ be the free abelian group
with basis the double cosets of $K_1(\gn)$ in $G(\A_f)$. It is endowed with 
algebra structure, where the product of two basis elements is given by : 
\begin{equation}\label{convolution}
[K_1(\gn)xK_1(\gn)]\cdot [K_1(\gn)yK_1(\gn)]=\sum_i [K_1(\gn)x_iyK_1(\gn)],
\end{equation}
where   $[K_1(\gn)xK_1(\gn)]=\coprod_i K_1(\gn) x_i$. 
For $g\in S_{k,J}(K_1(\gn))$ we put :
$$g|_{[K_1(\gn)xK_1(\gn)]}(\cdot)=\sum_i g(\cdot x_i^{-1}).$$

This defines an action of the algebra $\Z[K_1(\gn)\bs G(\A_f)/K_1(\gn)]$ on 
$S_{k,J}(K_1(\gn))$ (resp. on $G_{k,J}(K_1(\gn))$). Unfortunately, this algebra 
is not commutative  when $\gn\neq\go$.  We will now define a commutative 
subalgebra.  Consider the semi-group :
$$\Delta(\gn)=\left\{\begin{pmatrix} a & b \\ c &
    d\end{pmatrix}\in G(\A_f)\cap 
    \mathrm{M}_2(\widehat{\go})\enspace  |\enspace  d_v\in \go_v^\times,  \enspace 
c_v\in \gn_v\text{, for all } v\text{ dividing } \gn\right\}.$$

The abstract Hecke algebra of level $K_1(\gn)$, is defined as  
$\Z[K_1(\gn)\bs \Delta(\gn)/K_1(\gn)]$ endowed with the convolution product 
$(\ref{convolution})$. This algebra has the following explicit description.

For each ideal $\ga\!\subset\! \go$ we define the  Hecke 
operator $T_{\ga}$ as the finite sum of double  cosets
 $[K_1(\gn)xK_1(\gn)]$ contained  in  the set  $\{x\in
\Delta(\gn)|\nu(x)\go=\ga\}$.
In the  same way, for a prime to $\gn$ ideal $\ga\!\subset\! \go$, 
we define the Hecke operator
$S_{\ga}$ by the double  coset for $K_1(\gn)$
containing the scalar matrix of the id{\`e}le attached  to the ideal 
$\ga$.

 For each finite   place  $v$ of $F$, we have 
$T_v=K_1(\gn)\begin{pmatrix} \varpi_v & 0 \\ 0 & 1 \end{pmatrix} 
K_1(\gn)$, and for each $v$ not  dividing $\gn$ 
we have $S_v=K_1(\gn)\begin{pmatrix} \varpi_v & 0 \\ 0 & \varpi_v 
\end{pmatrix} K_1(\gn)$, where $\varpi_v$ is an uniformizer of $F_v$.

Then, the abstract Hecke algebra of level $K_1(\gn)$ is isomorphic to the 
polynomial algebra in  the  variables   $T_v$, where $v$ 
runs over the prime ideals of $F$, and the variables $S_v$,
where $v$  runs over the prime ideals of $F$ not dividing  $\gn$.
The  action of Hecke algebra obviously  preserves  the  decomposition 
 (\ref{decomp-fmh}) and moreover, $S_v$ acts on $S_{k,J}(\gn,\psi)$
as the scalar $\psi(v)$. 

\medskip
Let $\T(\C)=\T_k(\gn,\psi;\C)$ be the subalgebra of $\End_{\C}(S_{k,J}(\gn,\psi))$
 generate by   the operators $S_v$ for $v \nmid \gn$  and $T_v$ for all $v$ 
(we will see in   \S\ref{conjugates} that  $\T(\C)$ does not  depend on $J$). 

The algebra $\T(\C)$  is commutative, but  not semi-simple in  
general. Nevertheless, for   $v \nmid \gn$  the operators $S_v$ and $T_v$
 are normal with respect to the Petersson inner product (see Def.\ref{peterson}). 
Denote by  $\T'(\C)$ be the subalgebra of $\T(\C)$   generated by the Hecke operators
outside a finite set of places containing those dividing $\gn$. 
The algebra $\T'(\C)$  is  semi-simple, that is to say  $S_{k,J}(\gn,\psi)$ has a basis
 made of eigenvectors for   $\T'(\C)$.

\medskip

We will now describe the relation between Fourier coefficients and eigenvalues 
for the Hecke operators. By (\ref{ad-cl}) we can associate to $g\in S_k(K_1(\gn))$ 
a family of classical cusp forms  $g_i\in S_k(\Gamma_1(\gc_i,\gn);\C)$, where 
$\gc_i$ are  representatives  of the narrow ideal class group   $\Cl_F^+$.

Each form $g_i$ is determined by its $q$-expansion 
at the cusp $\infty$ of $M_1(\gc_i,\gn)^{\an}$.
For each fractional ideal  $\ga=\gc_i \xi$, 
with $\xi\in F_+^\times$,  we put  $c(g,\ga)=\xi^m a_\xi(g_i)$.
By \S\ref{qdev} for each $\epsilon \in \go_+^\times$,  we have 
$a_{\epsilon\xi}=\epsilon^{k+m-t}a_{\xi}$
and  therefore the definition of $c(g,\ga)$ does not depend 
on the  choice of $\xi$ (nor on  the particular choice  of the 
ideals   $\gc_i$; see \cite{hida-padic} IV.4.2.9.). 

\begin{defin}
 We say that $g\in S_k(\gn,\psi)$ is an {\it eigenform}, if it is an  
eigenvector for  $\T(\C)$. We say that an eigenform  $g$ is  {\it normalized}
if  $c(g,\go)=1$.
\end{defin} 

\begin{lemma}{\rm (\cite{hida128} Prop.4.1, \cite{hida-padic} (4.64))}
If $g\in S_k(\gn,\psi)$ is a normalized eigenform, then the  eigenvalue of 
$T_{\ga}$ on $g$ is  equal  to the   Fourier coefficient $c(g,\ga)$.
\end{lemma}

A consequence of this lemma and the $q$-expansion Principle 
(see \S\ref{qdev}), is the Weak Multiplicity One Theorem stating
 that two normalized  eigenforms having  the same   eigenvalues are equal.

\medskip
\noindent{\bf Ordinary modular forms.}
When the weight $k$ is non-parallel, the definition of the 
Hecke operators should be slightly modified. We 
 put $T_{0,v}=\varpi_v^{-m}T_v$ and 
$S_{0,v}=\varpi_v^{-2m}S_v$ (see \cite{hida128} Sect.3 ;
in the applications our base ring will be the  $p$-adic ring
$\cO$ which satisfies  the assumptions of this reference).

The advantage of the Hecke operators $T_{0,v}$ and $S_{0,v}$ is that 
they preserve in an optimal way the  $\cO$-integral 
structures on the space of Hilbert modular forms and on the 
cohomology of the  Hilbert modular variety.

\begin{defin} \label{ord-autom}
A Hilbert modular eigenform is {\it ordinary } at $p$
if, for all primes $\gp$ of $F$ dividing $p$, the image by $\iota_p$ of
its  $T_{0,\gp}$-eigenvalue   is a  $p$-adic unit.
\end{defin}

\medskip
\noindent{\bf Primitive modular forms.}
For each $\gn_1$ dividing $\gn$ and divisible 
by the conductor of $\psi$,  and for all $\gn_2$
dividing $\mathfrak{nn}_1^{-1}$ we consider  the linear map 
$$S_k(\gn_1,\psi)\rightarrow S_k(\gn,\psi),
\enspace g\mapsto g|\gn_2,$$

\noindent where $g|\gn_2$ is determined by the relation 
$c(\ga,g|\gn_2)=c(\ga\gn_2^{-1},g)$.

We define  the subspace $S_k^{\mathrm{old}}(\gn,\psi)$
of  $S_k(\gn,\psi)$ as the subspace  
generated by the images of all these  linear maps. This 
space is preserved  by the Hecke operators outside $\gn$.
We define  the space  $S_k^{\mathrm{new}}(\gn,\psi)$ 
of the primitive modular forms as the orthogonal of 
$S_k^{\mathrm{old}}(\gn,\psi)$
in $S_k(\gn,\psi)$ with respect to the Petersson
inner product (see Def.\ref{peterson}). Because the 
Hecke operators  outside  $\gn$ are normal 
for the Petersson inner  product, the direct sum decomposition 
$S_k(\gn,\psi)=S_k^{\mathrm{new}}(\gn,\psi)
\oplus S_k^{\mathrm{old}}(\gn,\psi)$ is preserved
 by  $\T'(\C)$.
The  Strong Multiplicity One  Theorem, due to
 Miyake in the Hilbert modular case, asserts  that if 
$f\in S_k^{\mathrm{new}}(\gn,\psi)$ 
 is an eigenform  $\T'(\C)$, then it is an eigenform  for $\T(\C)$. 

A normalized primitive eigenform is called a {\it newform}.

\medskip

The pairing $\T(\C)\times S_k(\gn,\psi)\rightarrow \C$, 
$(T,g)\mapsto c(g|_{T},\go)$ is   a perfect  duality (see \cite{hida128} Thm.5.2).

\subsection{External  and  Weyl group  conjugates.}\label{conjugates}

For an element $\sigma\in \Aut(\C)$ we define the
 {\it external conjugate} of $g\in S_k(K_1(\gn))$, 
as the unique element $g^\sigma \in S_k(K_1(\gn))$
satisfying $c(g^\sigma,\ga)\!=\!c(g,\ga)^\sigma$,
for each ideal $\ga$ of $\go$.

\medskip
We  identify  $\{\pm 1\}^{J_F}$ with the Weyl group  $K_{\infty}/K_{\infty}^+$ 
of $G$,  by sending $\epsilon_J\!=\!(-1_J, 1_{J_F\!\bs\!J})$
to $c_JK_{\infty}^+$, where for all $\tau\in J_F$,  $\det(c_{J,\tau})< 0$ 
if and only if $\tau\in J$. The length of $\epsilon_J$ is $|J|$.

We have an action of the Weyl group on the space of Hilbert modular forms.
More precisely, $\epsilon_J$ acts as the double class 
 $[K_1(\gn)c_JK_1(\gn)]$, and  maps   bijectively 
$S_k(K_1(\gn))$ onto $S_{k,J_F\!\bs\!J}(K_1(\gn))$.
The action of $\epsilon_J$ commutes with the action of the  Hecke operators.
For an element $g\in S_k(K_1(\gn))$  we put $g_J=\epsilon_{J_F\!\bs\!J}\cdot g$.

\subsection{Eichler-Shimura-Harder isomorphism.}\label{ESHI}
Let $R$ be  an $\cO$-algebra  and  $V_n(R)$ be the polynomial ring 
over  $R$ in the  variables $(X_\tau,Y_\tau)_{\tau\in J_F}$ which are  
homogeneous  of degree ${n}_\tau$ in $(X_\tau,Y_\tau)$. 
We have a pairing (perfect if $n_0!$ is invertible in $R$)
\begin{equation} \label{theta-n}
\langle\quad,\quad\rangle:V_n(R)\times V_n(R)\rightarrow R
\text{, given by}
\end{equation}
\begin{displaymath}\left\langle\sum_{0\leq j\leq  {n}}
a_jX^{{n}-j}Y^j,\sum_{0\leq j\leq {n}}
b_jX^{{n}-j}Y^j\right\rangle=
\sum_{0\leq j\leq {n}}(-1)^j a_j b_{{n}-j}\binom{{n}}{j}
\text{,where } \binom{{n}}{j}=\prod_{\tau\in
  J_F}\binom{{n}_\tau}{j_\tau},
\end{displaymath} 

The   $R$-module $V_n(R)$  realizes   the  algebraic 
representation  $V_n=\bigotimes_{\tau}(\Sym^{n_\tau}\otimes \det^{m_\tau}) $ of 
$G(R)$. We endow  $V_n(R)$  with an   action of $(M_2(\cO)\cap \GL_2(E))^{J_F}$
given by
$$\gamma.P((X_\tau,Y_\tau)_{\tau\in
  J_F})=\nu(\gamma)^m P((\det(\gamma)\gamma^{-1})^t
(X_\tau,Y_\tau)_{\tau\in J_F}).$$




Let $\bV_n(R)$ be the sheaf of continuous (thus locally constant) sections of  
$$G(\Q)\bs G(\A)\times V_n(R)/K_1(\gn)K_{\infty}^+\rightarrow
G(\Q)\bs G(\A)/K_1(\gn)K_{\infty}^+=Y^{\an},$$
 where $y\in K_1(\gn)K_{\infty}^+$ acts on $V_n(R)$ via its $p$-part $y_p$.

For each $y\in \Delta(\gn)$ the map $[y]:G(\A)\times V_n(R)\rightarrow 
G(\A)\times V_n(R)$,  $(x,v)\mapsto (xy,y_p.v)$ is a homomorphism of sheaves. 
This induces   an action of the   Hecke operator $[K_1(\gn)yK_1(\gn)]$ on   
$\rH^d(Y^{\an},\bV_n(R))$  preserving the cuspidal cohomology 
$\rH_{\mathrm{cusp}}^d(Y^{\an},\bV_n(R))$.

The action  of $\epsilon_J$ on $(M^{\an },\bV^{\an})$ given by 
$\epsilon_J \cdot ((z_J,z_{J_F\!\bs\!J}),v)=((-\overline{z}_J,z_{J_F\!\bs\!J}),v)$,
induces an action of the Weyl group  on $\rH^d (Y^{\an} , \bV^{\an})$
commuting  with the Hecke action.

By  Harder \cite{harder} we know that, if $n\neq 0$, then
$\rH^d_{!}(Y^{\an},\bV_n(\C))=\rH_{\mathrm{cusp}}^d(Y^{\an},\bV_n(\C))$.

By  (\ref{theta-n}) we have  a  Poincar{\'e} pairing
$\langle\quad,\quad \rangle: \rH^d(Y^{\an},\bV_n(R))
 \times \rH^d_c(Y^{\an},\bV_n(R))  \rightarrow R$.

Let $\eta$ be the id{\`e}le corresponding to the ideal $\gn$ and  let  
$\iota=\begin{pmatrix}  0 & 1 \\ -\eta & 0 \end{pmatrix}$ be the Atkin-Lehner 
involution. By putting $[x,y]=\langle x,\iota y \rangle $ we obtain a new pairing
\begin{equation}\label{accoupl-twisted}
[\quad,\quad]: \rH_{!}^d(Y^{\an},\bV_n(R))\times\rH_{!}^d(Y^{\an},\bV_n(R))
\rightarrow R,
\end{equation}
which is Hecke-equivariant. We call it the {\it twisted Poincar{\'e} pairing}.

\medskip

Now we state the  Eichler-Shimura-Harder  isomorphism :

\begin{theo} 
{\rm (Hida, \cite{hida74})}
If  $n\neq 0$, then there exists an isomorphism :
\begin{equation}\label{ESH}
\delta:\bigoplus_{\psi}\bigoplus_{J\subset J_F}S_{k,J}(\gn,\psi)\cong
\rH_{!}^d(Y^{\an},\bV_n(\C)),
\end{equation}
\noindent where $\psi$ runs over the Hecke  characters of conductor dividing 
$\gn$ and  type $-n_0t$ at infinity. This  isomorphism is equivariant for 
the actions of the  Hecke algebra and  the Weyl  group. 
\end{theo}

\medskip
For  each $J\!\subset\! J_F$ let $\widehat{\epsilon}_J:
\{\pm 1\}^{J_F}\rightarrow \{\pm 1\}$  be the  unique character of the Weyl  
group  sending  $\epsilon_\tau=(-1_\tau,1^\tau)$ 
to  $1$, if $\tau\in J$, and to  $-1$  if   $\tau\in J_F\!\bs\!J$. 
The restriction of the Eichler-Shimura-Harder isomorphism (\ref{ESH})
to $S_{k,J}(\gn,\psi)$, followed by the projection on the 
$(\psi,\widehat{\epsilon}_J)$-part yields a Hecke equivariant isomorphism
\begin{equation}\label{ESH-epsilon}
\delta_J:S_{k,J}(\gn,\psi)\cong \rH_{!}^d
(Y,\bV_n(\C))[\psi,\widehat{\epsilon}_J].
\end{equation}

\medskip
Moreover, after twisting by the complex conjugation $c$ on the coefficients, 
we still have   a direct sum decomposition :
\begin{equation}\label{hodge-decomp}
\rH^d(M^{\an }, \bV_n(\C))=\bigoplus_{J\subset J_F} 
\rH^d(M^{\an }, \bV_n(\C))[\widehat{{\epsilon}_J\!\otimes\! c}].
\end{equation}

This  decomposition is  finer than the usual Hodge decomposition,
whose graded are given by ($0\leq a\leq d$) :
$$ \gr^a\rH^d(M^{\an }, \bV_n(\C))=
\bigoplus_{J\!\subset\! J_F, |J|=a} \rH^d(M^{\an },
 \bV_n(\C))[\widehat{{\epsilon}_J\!\otimes\! c}].$$

The transcendental decomposition (\ref{hodge-decomp}) has an algebraic 
interpretation, via the so-called BGG complex, that we will
describe in the next section.

\section{Hodge-Tate weights of the Hilbert modular varieties.}

 The aim of this  section  is to determine   the   Hodge-Tate weights 
of the $p$-adic {\'e}tale  cohomology of the Hilbert modular variety 
$\rH^\bullet(M_{\overline{\Q}_p},\bV_n(\overline{\Q}_p))$, as well as
those of the $p$-adic Galois representation associated to a
Hilbert modular form. In all this section we assume 

\medskip
{\bf (I)} $p$ does not divide $\Delta= \Nm(\gn\gd)$.

\medskip
The proof relies on Faltings'  Comparison Theorem \cite{Fa-jami}
relating  the {\'e}tale cohomology of $M$ with coefficients in the 
local system $\bV_n(\Q_p)$ to the de Rham logarithmic cohomology
of the corresponding vector bundle $\overline{\cV}_n$ over 
a smooth toroidal compactification $\overline{M}$ of $M$.
  The Hodge-Tate weights are given  by the jumps of 
the Hodge filtration of the associated de Rham complex. 
These are computed, following 
\cite{FaCh},  using the so-called  Bernstein-Gelfand-Gelfand complex (BGG complex).

Instead of using Faltings'  Comparison Theorem, one can 
use Tsuji's result for  the  {\'e}tale cohomology with 
constant coefficients of the Kuga-Sato variety  $\cA^s$ 
($s$-fold fiber product of the universal abelian variety $\cA$ 
above the fine moduli space $M^1$ associated to $M$; see
\cite{dimtildg} Sect.6 for the construction of  toroidal compactifications
of  $\cA^s$).

\bigskip
For each  subset $J$ of $J_F$ we put 
$p(J)=\sum_{\tau\in J}(k_0\!-\!m_\tau\!-\!1)\tau+
\sum_{\tau\in J_F\bs J}m_\tau\tau\in\Z[J_F]$ and for each
$a=\sum_{\tau\in J_F}a_\tau\tau\in\Z[J_F]$ we put
$|a|=\sum_{\tau\in J_F}a_\tau\in \Z$.

\subsection{Motivic weight of the  cohomology.} \label{MWC}
Consider the smooth sheaf $R^1\pi_* \overline{\Q}_p$
on $M^1$, where $\pi:\cA\rightarrow M^1$ is the universal HBAV. It corresponds
to a representation of the fundamental group of $M^1$ in $G(\overline{\Q}_p)$.
By composing this representation with the  algebraic representation $V_n$  of $G$
of highest weight $n$ (see \S\ref{ESHI}),
we obtain a smooth sheaf on $M^1$ (thus on $Y^1$). It descends to  
a smooth sheaf on $Y$ denoted  $\bV_n(\overline{\Q}_p)$. 

Let  $W_f=\underset{\ga\subset \go}{\bigcap}
\ker(T_{\ga}-c(f,\ga))$ be the subspace 
of  $\rH^d(Y_{\overline{\Q}},\bV_n(\overline{\Q}_p))$ 
corresponding  to the   Hilbert modular newform  $f\in S_k(\gn,\psi)$.
Put $s=\sum_\tau(n_\tau+2m_\tau)=d n_0$.

\begin{prop}  $W_f$ is pure of weight $d+s$, that is to say  for all prime 
$l\nmid p\Delta$ the eigenvalues  of the geometric Frobenius 
 $\Frob_l$ at $l$ are Weil numbers of  absolute  value  $l^{\frac{d+s}{2}}$. 
\end{prop}

{\bf Proof : } As $f$ is cuspidal $W_f\!\subset\! 
\rH^d_!(Y_{\overline{\Q}},\bV_n(\overline{\Q}_p))$. We recall that 
 $Y_{\overline{\Q}}$ is a disjoint union  of its connected components 
 $M_{\overline{\Q}}=M_1(\gc_i,\gn)_{\overline{\Q}}$,
where the $\gc_i$'s form  a set of representatives of $\Cl_F^+$. 
Let $\gc$ be one of the $\gc_i$'s and  $M^1=M_1^1(\gc,\gn)$. 
For $*=\varnothing,c$ we have 
$$\rH^0(\go_+^\times/\go^{\times 2}_{\gn,1}, 
\rH^d_*(M_{\overline{\Q}}^1,\bV_n(\overline{\Q}_p)))=
\rH^d_*(M_{\overline{\Q}},\bV_n(\overline{\Q}_p)),$$
 and therefore, it is enough to 
prove that $\rH^d_!(M_{\overline{\Q}}^1,\bV_n(\overline{\Q}_p))$
is pure of weight $d+s$. We use  Deligne's method
\cite{De2}. Let $\pi: \cA\rightarrow  M^1$ be the universal
abelian variety (see \S\ref{hmv}). The  sheaf 
$\bV_n(\overline{\Q}_p)$ corresponds to the representation 
$\bigotimes_{\tau\in J_F}\Sym^{n_\tau}\otimes \det^{m_\tau}$
of the  group $G^*$ and can therefore be cut out by algebraic 
correspondences in   $(R^1\pi_*\overline{\Q}_p)^{\otimes s}$. Let  
 $\pi_s: \cA^s\rightarrow  M^1$ be the Kuga-Sato variety. By the 
Kunneth's formula we have
$$\rH^d_!(M_{\overline{\Q}}^1,(R^1\pi_*\overline{\Q}_p)^{\otimes s})\!\subset\!
\rH^d_!(M_{\overline{\Q}}^1,R^s\pi_{s*}\overline{\Q}_p)\!\subset\!
\rH^{d+s}_!(\cA_{\overline{\Q}}^s,\overline{\Q}_p)\!\subset\!
\rH^{d+s}(\overline{\cA_{\overline{\Q}}^s},\overline{\Q}_p), $$

\noindent where the middle inclusion comes from the degeneration 
of the of the Leray spectral sequence  $\rE_{2*}^{i,j}=
\rH^i_*(M_{\overline{\Q}}^1,R^j\pi_{s*}\overline{\Q}_p)\Rightarrow 
\rH^{i+j}_*(\cA_{\overline{\Q}}^s,\overline{\Q}_p)$ 
for $*=\varnothing,c$ (see \cite{De2}).
The proposition is then a consequence of the Weil conjectures
for the eigenvalues of the Frobenius, proved by  Deligne  \cite{De1}.
\hfill$\square$

\subsection{The Bernstein-Gelfand-Gelfand complex over $\overline{\Q}$.}
In  this and the next sections we  give, 
following Faltings \cite{Fa},   an algebraic construction 
of the  transcendental  decomposition of the 
Betti cohomology described in (\ref{hodge-decomp}).

In this section all the objects are defined over a
characteristic zero field  splitting  $G$.

\medskip
Let $\mathfrak{g}$, $\gb$, $\gt$ and $\gu$ 
denote the Lie algebras of  $G$, $B$, $T$ and $U$,  respectively.
Consider the canonical splitting 
$\mathfrak{g}=\gb\oplus \gu^-$.
Let  $U(\mathfrak{g})$, $U(\gb)$ be  the enveloping  algebras of 
$\mathfrak{g}$ and  $\gb$, respectively.

\medskip
 The aim of this section is to write down a  resolution of $V_n$  of the type :
$$ 0 \leftarrow  V_n\leftarrow U(\mathfrak{g})
\otimes_{ U(\gb)}K_n^{\bullet},$$
where the  $K_n^j$  are finite dimensional semi-simple 
$\gb$-modules, with explicit simple components.

We start  by the case $n=0$.
If we put $K_0^j=\wedge^j(\mathfrak{g}/\gb)$ we 
obtain the so-called   {\it bar-resolution} of $V_0$.
Note that $\wedge^i(\mathfrak{g}/\gb)$ is a $\gb$-module
with trivial action of  $\gu$, therefore 
$K_0^j=\oplus W_{\mu}$, where $\mu$  runs over the weights  of $B$
that are  sum of $j$ distinct negative roots.

 By tensoring  this resolution with $V_n$ we  obtain the Koszul's complex :
\begin{equation}\label{bgg0}
 0 \leftarrow  V_n\leftarrow U(\mathfrak{g})
\otimes_{ U(\gb)} \left(\wedge^j(\mathfrak{g}/
\gb)\otimes V_n|_{\gb}\right),
\end{equation}
which is a resolution of  $V_n$ by   $\gb$-modules 
$\wedge^i(\mathfrak{g}/ \gb)\otimes V_n|_{\gb}$,  
 not   semi-simple in  general. 

\medskip

The BGG complex that we are going to define is  a direct factor 
of the Koszul's complex, cut by the action of the center 
$U(\mathfrak{g})^G$ of $U(\mathfrak{g})$.

Denote by   $\chi_n$ the character of  $U(\mathfrak{g})^G$
 corresponding to the  weight $n$. It is a classical  result  that 

\begin{lemma} 
$\chi_n=\chi_{\mu}$, if and only if there 
exists $J\!\subset\! J_F$ such that $\mu=\epsilon_J(n+t)-t $.
\end{lemma}

By taking the $\chi_n$-part of the bar resolution (\ref{bgg0}) 
of  $V_n$ we obtain the complex  :
\begin{equation}\label{bgg}
 0 \leftarrow V_n \leftarrow U(\mathfrak{g})
\otimes_{ U(\gb)}K_n^{\bullet}, \enspace
\mathrm{with}  \enspace K_n^i=
\bigoplus_{J\!\subset\! J_F, |J|=i} W_{\epsilon_J(n+t )-t },
\end{equation}
which  is  still  a resolution of $V_n$, as it is a direct factor 
of a   resolution. We call this resolution  the BGG complex.

\subsection{Hodge-Tate decomposition of  $\rH^\bullet(M\otimes 
\overline{\Q}_p,\bV_n(\overline{\Q}_p))$.} \label{hodge-tate} ${}$
In this paragraph we summarize the results  of \cite{dimtildg} Sect.7. 
The algebraic groups $G$, $B$, $T$ and $D$ of \S\ref{hmfv}  
have models over $\Z$, denoted by the same letters. 
For every scheme $X$, we put $X'=X\times \Spec(\go'[\frac{1}{\Delta}])$.

\medskip
By \S\ref{kugasato},  we can extend the  vector bundles $\underline{\omega}$ and 
$\cH^1_{\dR}$ to $\overline{M}{}'$. Only the 
construction depends on a choice of a toroidal compactification
 $\overline{\pi}:\overline{\cA}\rightarrow \overline{M^1}$  of 
$\pi:\cA\rightarrow M^1$. 

The sheaf  $\underline{\gM_D}=\underline{\Isom}_{\go\otimes 
\cO_{\overline{M}{}'}} (\underline{\omega},\go\otimes \cO_{\overline{M}{}'})$ 
is a $D'$-torsor over $\overline{M}{}'$ (for the Zariski topology).
We have a functor $\cF_D$ from the category of algebraic representations
of $D'$  to the  category of  vector bundles on $\overline{M}{}'$
which are direct sum of invertible bundles. To an  algebraic representation
$W$ of $D'$, $\cF_D$ associates the fiber product 
$\overline{\cW}:=\gM_D\overset{D'}{\times} W$. 

\medskip
The sheaf  $\underline{\gM_B}=\underline{\Isom}_{\go\otimes 
\cO_{\overline{M}{}'}}^{\fil}(\cH^1_{\ldR},
(\go\otimes \cO_{\overline{M}{}'})^2)$ is a $B'$-torsor over $\overline{M}{}'$.
We have a functor $\cF_B$ from the category of algebraic representations
of $B'$  to the  category of filtered vector bundles on $\overline{M}{}'$
whose graded are sums of invertible bundles. To an  algebraic representation
$V$ of $B'$, $\cF_B$ associates the fiber product 
$\overline{\cV}:=\gM_B\overset{B'}{\times} V$. 

A representation of $G$ (resp. $T$) can be considered as a representation 
of $B$ by restriction (resp. by making $U$ act trivially). Thus, we may define 
the  filtered vector  bundle $\overline{\cV}_n$ on $\overline{M}{}'$
associated to the algebraic representation $V_n$ of $G$, and the 
invertible bundle  $\overline{\cW}_{n,n_0}$ on  $\overline{M}{}'$ associated 
to the algebraic representation of $T=D\times D$, given by $(u,\epsilon)\mapsto 
u^n\epsilon^m$.

\medskip
The sheaf  $\underline{\gM_G}=\underline{\Isom}_{\go\otimes 
\cO_{\overline{M}{}'}}(\cH^1_{\ldR},
(\go\otimes \cO_{\overline{M}{}'})^2)$ is a $G'$-torsor over $\overline{M}{}'$.
We have a functor $\cF_G$ from the category of algebraic representations
of $G'$  to the  category of flat vector bundles on $\overline{M}{}'$
(that is  vector bundles  endowed with an 
integrable quasi-nilpotent logarithmic connection). 
To any algebraic representation $V$ of $B'$, $\cF_G$ associates the 
fiber product  $\overline{\cV}{}^{\nabla}:=\gM_G\overset{G'}{\times} V$.
For $j\!\in\!\n$, we  put
$\cH_{\ldR}^{j}(\overline{M}',\overline{\cV})=
R^j\phi_*(\overline{\cV}{}^{\nabla}\otimes \Omega^\bullet_{\overline{M}'}(\dlog\infty))$,
where $\phi: \overline{M}'\rightarrow \Spec(\go'[\frac{1}{\Delta}])$
denotes the structural homomorphism.

\medskip
By the Faltings' Comparison Theorem \cite{Fa-jami}, 
the $\G_{\Q_p}$-representation $\rH^\bullet 
(M^1_{\overline{\Q}_p},\bV_n(\overline{\Q}_p))$ is crystalline, hence de Rham, 
and we have a canonical isomorphism
$$\rH^\bullet(M^1_{\overline{\Q}_p},\bV_n(\overline{\Q}_p))\otimes B_{\dR}\cong
\cH_{\ldR}^\bullet(\overline{M^1}_{/\overline{\Q}_p},\overline{\cV}_n)\otimes B_{\dR}.$$

By \cite{dimtildg} Sect.7, the  Hodge to de Rham spectral sequence
$$\rE_1^{i,j}=\rH^{i+j}(\overline{M^1}_{/\overline{\Q}_p},\gr^i 
(\overline{\cV}_n\otimes 
\Omega^\bullet_{\overline{M^1}}(\dlog\,\infty)))
\Rightarrow \cH_{\ldR}^{i+j}(\overline{M^1}_{/\overline{\Q}_p},\overline{\cV}_n), $$
degenerates at $\rE_1$ (the filtration is 
the tensor product of the two Hodge filtrations). 
In order to compute the  jumps of the resulting filtration we introduce 
the  BGG complex :
$$\overline{\mathcal{K}}_n^i=\bigoplus_{J \!\subset\! J_F, |J|=i} 
\overline{\cW}_{\epsilon_J(n+t)-t,n_0}.$$
The fact that $\overline{\mathcal{K}}_n^\bullet$ is a complex
follows from  (\ref{bgg}) and from the following isomorphism
(see  \cite{FaCh} Prop.VI.5.1)
\begin{equation}\label{op-diff}
\Hom_{U(\mathfrak{g})}(U(\mathfrak{g})
\otimes_{U(\gb)}W_{1}),U(\mathfrak{g}) \otimes_{U(\gb)}W_{2})\rightarrow 
\mathrm{Diff.Op.}(\overline{\cW}_{2},\overline{\cW}_{1}).
\end{equation}
Define a  filtration on  $\overline{\mathcal{K}}_n^\bullet$  by 
$\Fil^i\overline{\mathcal{K}}_n^\bullet=\bigoplus_{J \!\subset\! J_F,
 |p(J)|\geq i}  \overline{\cW}_{\epsilon_J(n+t)-t,n_0}.$

As  the image  of the  Koszul's complex (\ref{bgg0}) by the 
contravariant functor  $W\mapsto \overline{\cW}$ 
is  equal to the de Rham  complex, and as the BGG complex
is  a direct (filtered) factor   of the Koszul's complex, we
have :

\begin{theo}\label{hodge-pad} {\rm (\cite{dimtildg} Thm.7.8)}
{\rm (i)} There is a  quasi-isomorphism of filtered complexes 
$$\overline{\mathcal{K}}_n^\bullet\hookrightarrow \overline{\cV}_n\otimes 
\Omega^\bullet_{\overline{M^1}}(\dlog\,\infty).$$
 
{\rm (ii)} The spectral sequence  given by the  Hodge filtration 

$$\rE_1^{i,j}=\bigoplus_{J \!\subset\! J_F, |p(J)|=i} 
\rH^{i+j-|J|} (\overline{M^1}_{/\overline{\Q}_p}, 
\overline{\cW}_{\epsilon_J(n+t)-t,n_0})
\Rightarrow
\cH_{\ldR}^{i+j}(\overline{M^1}_{/\overline{\Q}_p},\overline{\cV}_n)$$
degenerates at $\rE_1$.

{\rm (iii)} For all  $j\leq d$, the Hodge-Tate  weights 
 of the  $p$-adic  representation
$\rH^j(M^1_{\overline{\Q}_p},\bV_n(\overline{\Q}_p))$ 
belong to the set  $\{|p(J)|\enspace , \enspace |J|\leq j\}$.
\end{theo}

\subsection{Hecke operators on the cohomology.}\label{hecke2}
We describe  the standard Hecke  operator $T_{\ga}$ as a correspondence 
on $Y^1$. We are indebted to M. Kisin for pointing us out that the usual 
definition of Hecke operators  on $Y$ extends to $Y^1$ (see \cite{KiLa}\S1.9-1.11). 
Note that the  corresponding Hecke action on  analytic modular forms for $G^*$ 
(see \S\ref{hecke}) is not easy to write down, because the double 
class for the Hecke operator $T_v$ does not belong to $G^*(\A_f)$, 
unless $v$ is inert in $F$.

Recall that $Y_1^1(\gn)=\coprod_{i=1}^{h^+} M_1^1(\gc_i,\gn)$,
where $\gc_1$,..,$\gc_{h^+}$ are a set of representatives of $\Cl_F^+$.

Assume that $\gc_i\ga$ and $\gc_j$ have the same class in $\Cl_F^+$. 
Then, consider the contravariant functor $\underline{\cM}^1_{\ga}$
 from the category of  $\Z[\frac{1}{\Delta}]$-schemes to the category of sets, 
assigning to a scheme $S$  the set of isomorphism classes of quintuples 
$(A,\lambda,\alpha,C,\beta)$ where  $(A,\lambda,\alpha)/S$ is a $\gc_i$-polarized 
HBAV  with $\mu_{\gn}$-level structure, $C$ is a closed subscheme of $A[\ga]$ 
which is $\go$-stable, disjoint from  $\alpha(\Gm\otimes \mathfrak{d}^{-1})$ and 
locally isomorphic to the constant group scheme $\go/\ga$ over $S$, and $\beta$ is an 
$\go^{\times 2}_{\gn,1}$-orbit of isomorphisms $(\gc_i\ga,(\gc_i\ga)_+) 
\overset{\sim}{\longrightarrow} (\gc_j,\gc_{j+})$.

We have a projection  $\underline{\cM}^1_{\ga}\rightarrow
\underline{\cM}^1$, $(A,\lambda,\alpha,C,\beta)\mapsto (A,\lambda,\alpha)$
which is relatively representable by $\pi_1:M^1_{\ga}(\gc_i,\gn)\rightarrow
M^1_1(\gc_i,\gn)$. We have also a projection $\pi_2:M^1_{\ga}(\gc_i,\gn)
\rightarrow M^1_1(\gc_j,\gn)$, coming from 
$(A,\lambda,\alpha,C,\beta)\mapsto(A/C,\lambda',\alpha')$, where 
$\alpha'$  is the composed map of $\alpha$ and $A\rightarrow A/C$, and 
$\lambda'$ is a $\gc_j$-polarization  of $A/C$ (defined via $\lambda$ and $\beta$).

Put $Y^1_{\ga}=\coprod_{i=1}^{h^+} M^1_{\ga}(\gc_i,\gn)$. 
As $\gc_i\mapsto \gc_j\simeq \gc_i\ga$ is a permutation of $\Cl_F^+$, we get 
two finite projections $\pi_1,\pi_2:Y^1_{\ga}\rightarrow Y^1$. 
$$\xymatrix@C=30pt@R=20pt{
\cA\ar[d]_{\pi}  &\cA_{\ga}\ar[d]^{\pi_{\ga}} \ar[r]\ar[l] & \cA \ar[d]^{\pi}\\
Y^1 & Y^1_{\ga}\ar[r]^{\pi_2}\ar[l]_{\pi_1} & Y^1 }$$

From this diagram we obtain 
$\pi_2^*\cH^1_{\dR}\rightarrow \pi_1^*\cH^1_{\dR}$.
Therefore, for every  algebraic representation $V$ of $G$, we have 
$\pi_2^*\cV^\nabla\rightarrow \pi_1^*\cV^\nabla$. 
By composing this morphisms by  $\pi_{1*}$ and taking  the trace, we obtain 
$\cV^\nabla \rightarrow  \pi_{1*}\pi_2^*\cV^\nabla\rightarrow
 \pi_{1*}\pi_1^*\cV^\nabla \rightarrow \cV^\nabla$, and thus 
we obtain an action of $T_{\ga}$ on $\rH^\bullet(Y^1,\cV^\nabla)$.

The same way from the above diagram we obtain 
$\pi_2^*\underline{\omega}\rightarrow \pi_1^*\underline{\omega}$ and 
 $\pi_2^*\underline{\nu}\rightarrow \pi_1^*\underline{\nu}$.
Therefore, for each  algebraic  representation  $W$ of $T$, we get 
$\pi_2^*\cW\rightarrow \pi_1^*\cW$.
In order to define the good action of  $T_{\ga}$ on on Hilbert modular forms, we 
should modify slightly the last arrow :
 we decompose  $\cW$ as $(\cW\underline{\omega}^{-2t})\underline{\omega}^{2t}$ 
and  we define  $\pi_{2*}(\cW\underline{\omega}^{-2t})
\rightarrow\pi_{1*}(\cW\underline{\omega}^{-2t})$ as above and 
$\pi_{2*}\underline{\omega}^{2t}\rightarrow\pi_{1*}\underline{\omega}^{2t}$
as in \cite{KiLa}\S1.11 (via the Kodaira-Spencer isomorphism 
$\Omega^1_{Y^1}\simeq \underline{\omega}^2\otimes_{\go}\gd\gc^{-1}$). 
Thus we obtain 
$\cW \rightarrow  \pi_{1*}\pi_2^*\cW\rightarrow \pi_{1*}\pi_1^*\cW\rightarrow\cW,$
and an  action of  $T_{\ga}$  on  $\rH^\bullet(Y^1,\cW)$. 

In particular, we obtain an action of  $T_{\ga}$  on the space 
of geometric Hilbert modular forms for $G^*$, 
 $\rH^0(Y^1,\underline{\omega}^k\otimes \underline{\nu}^{-n_0t/2})$.  As it 
has been observed in \cite{KiLa}1.11.8, this action is given by the projection 
$$\textstyle\frac{1}{[\go_+^\times:\go_{\gn,1}^\times]}\sum_{[\epsilon]\in
\go_+^\times/\go_{\gn,1}^\times} [\epsilon]\cdot:\rH^0(Y^1,\underline{\omega}^k\otimes \underline{\nu}^{-n_0t/2}) \rightarrow \rH^0(Y,\underline{\omega}^k\otimes \underline{\nu}^{-n_0t/2}), $$
followed by the usual Hecke operator on the space of  Hilbert modular forms 
(see \S\ref{hecke}).

\subsection{Hodge-Tate weights of $\otimes\Ind_F^{\Q}\rho$
 in the crystalline case.} \label{HT-induced}

We first recall the notion of induced representation.
Let $V_0$ be vector space over a field $L$, and let
 $\rho_0:\G_F\rightarrow \GL(V_0)$ be a linear representation.
The  induced representation $\Ind_F^{\Q}\rho_0$ of $\rho_0$ from $F$ to  $\Q$ 
is by definition  the $L$-vector space  
 $$\Hom_{\G_F}(\G_{\Q},V_0):= \{\phi_0:\G_{\Q}\rightarrow V_0\enspace 
|\enspace  \forall h\in\G_F,\enspace g\in \G_{\Q},\enspace 
\phi_0(gh)=\rho_0(h^{-1})(\phi_0(g))\},$$ 
where  $g\in \G_{\Q}$ acts on  $\phi_0\in \Hom_{\G_F}(\G_{\Q},V_0)$ by 
 $g\cdot \phi_0 (\cdot)=\phi_0(g^{-1}\cdot)$. 

\smallskip
For any fixed decomposition  $\G_{\Q}= \underset{\tau\in J_F}{\coprod}\widetilde{\tau}\G_F$,
the map  $\phi\mapsto (\phi(\widetilde{\tau}))_{\tau}$
gives an isomorphism between  $\Hom_{\G_F}(\G_{\Q},V_0)$ and 
 the direct sum  $\underset{\tau}{\oplus} V_\tau$ 
(where each $V_\tau$ is isomorphic to  $V_0$).
Via this identification,  the action of $\G_{\Q}$ on
 $\underset{\tau}{\oplus} V_\tau$ is given by :
$$ (\Ind_F^{\Q}\rho_0)(g)((v_\tau)_{\tau})=
(\rho_0(\widetilde{\tau}^{-1}g \widetilde{\tau_g})(v_{\tau_g}))_{\tau},$$
where $g^{-1}\widetilde{\tau}\in \widetilde{\tau_g}\G_F$. In fact
$(\phi_0(\widetilde{\tau}))_{\tau}\overset{g}{\mapsto}
(\phi_0(g^{-1}\widetilde{\tau}))_{\tau}=
(\rho_0(\widetilde{\tau}^{-1}g\widetilde{\tau_g})
(\phi(\widetilde{\tau_g})))_{\tau}$.

Keeping  the same  notations we define, following  Yoshida \cite{yoshida}, 
the tensor induced  representation $\otimes\Ind_F^{\Q}\rho_0:
\G_{\Q}\rightarrow \GL(\underset{\tau}{\otimes} V_\tau)$ as :
$$(\otimes\Ind_F^{\Q}\rho_0)(g)(\underset{\tau}{\otimes}  v_\tau)=
\underset{\tau}{\otimes} \rho_0(\widetilde{\tau}^{-1}g\widetilde{\tau_g})(v_{\tau_g}).$$

\begin{rque} For each $g\in \G_{\Q}$ the map  $\tau\mapsto \tau_g$
is a  permutation of $J_F$, and it  is trivial if and only if 
$g\in \G_{\widetilde{F}}$. Therefore, for each  $g\in \G_{\widetilde{F}}$,
we have  $(\otimes\Ind_F^{\Q}\rho_0)(g)=\underset{\tau}{\otimes}
\rho_0(\widetilde{\tau}^{-1}g\widetilde{\tau})$.  
Moreover for all $g,g'\in \G_{\Q}$  we have $(\tau_g)_{g'}=\tau_{gg'}$.
\end{rque}

\begin{defin} The {\it internal conjugate} $g_\tau$ 
of $g$ by $\tau\in J_F$, is defined as the unique element $g_\tau \in 
S_{k^\tau,J^\tau}(\tau(\gn),\psi_\tau)/\tau(F)$
satisfying $c(g_\tau,\ga)\!=\!c(g,\tau(\ga))$,
for each ideal $\ga$ of $\go$, where $k^\tau=\sum_{\tau'} k_{\tau\tau'}\tau'$ 
and  $\psi_\tau(\ga)\!=\!\psi(\tau(\ga))$. 
\end{defin} 

If  $\rho=\rho_{f,p}$ by the  previous remark we have 
$(\otimes\Ind_F^{\Q}\rho)(g)=
\underset{\tau}{\otimes} \rho_{f_{\tau}}(g)$, for each  $g\in \G_{\widetilde{F}}$.

 Brylinski and Labesse \cite{BL} have shown (see \cite{Ta} for 
this formulation)

\begin{theo} {\rm (Brylinski-Labesse)}\label{BL-thm}
The restrictions to $\G_{\widetilde{F}}$ of the two $\G_{\Q}$-representations  
 $W_f$ and $\otimes\Ind_F^{\Q}\rho$ have the same characteristic polynomial.
\end{theo}

\begin{cor}\label{hodge-complex}  {\rm(\cite{dimtildg} Cor.7.9)}
{\rm (i)} The spectral sequence  given by the  Hodge filtration 

$$\rE_1^{i,j}=\bigoplus_{J \!\subset\! J_F, |p(J)|=i} 
\rH^{i+j-|J|} (\overline{M}_{/\overline{\Q}_p}, \overline{\cW}_{\epsilon_J(n+t)-t,n_0}) \Rightarrow
\cH_{\ldR}^{i+j}(\overline{M}_{/\overline{\Q}_p},\overline{\cV}_n)$$
degenerates at $\rE_1$ and is Hecke equivariant.

{\rm (ii)}
 The   Hodge-Tate weights of $W_f$ are the integers
$\vert p(J)\vert$, $J\!\subset\! J_F$, counted with multiplicity.
\end{cor}

{\bf Proof :}  
(i) By taking the invariants of the Hodge filtration of  $\overline{\cV}_n \otimes 
\Omega^\bullet_{\overline{M^1}}(\dlog\,\infty)$ by the Galois group
 of  the {\'e}tale covering  $\overline{M^1}\rightarrow \overline{M}$,
we  obtain a filtration of the complex $\overline{\cV}_n \otimes 
\Omega^\bullet_{\overline{M}}(\dlog\,\infty)$ on $\overline{M}{}'$, 
 still called  the  Hodge filtration. The  same way, we define 
the BGG complex over $\overline{M}{}'$ by taking the  invariants of the BGG complex
over  $ \overline{M^1}$. The   associated spectral sequence is given by 
the invariants of the spectral sequence of  Thm.\ref{hodge-pad}(ii). 
We have now to  see that it is Hecke equivariant.

The Hecke operator $T_{\ga}$ extends to $\overline{Y^1}$. One way to define it 
is to take the schematic closure of $T_{\ga}\subset Y^1\times Y^1$ in 
$\overline{Y^1}\times \overline{Y^1}$. Another way is to take a 
toroidal compactification $\overline{Y^1_{\ga}}$ of $Y^1_{\ga}$ over  the
toroidal compactification  $\overline{Y^1}$ of $Y^1$. In both cases we obtain 
an action of  $T_{\ga}$ on   $\rH^\bullet(\overline{Y^1},\overline{\cW})$ and 
on $\rH^\bullet(\overline{Y^1},\overline{\cV}{}^\nabla)$.  Although it is not clear
in general that these extended Hecke operators commute. Nevertheless, they
commute on the right  hand side  of  Thm.\ref{hodge-pad}(ii), because
by Faltings' Comparison Theorem this side  is independent
of the  toroidal compactification.
Since the spectral sequence  of Thm.\ref{hodge-pad}(ii) degenerates at $\rE_1$, 
 they should also commute  on the left hand side. 

(ii) We have   $\overline{\cW}_{\epsilon_J(n+t)-t,n_0}=
\underline{\omega}^{-\epsilon_J(n+t)+t}\otimes \underline{\nu}^{p(J)}$. 
 It follows from Thm.\ref{hodge-pad} (as in \cite{FaCh} Thm.5.5
 and \cite{MoTi} Sect.2.3) that the jumps of the Hodge filtration are 
among $|p(J)|$, $J\!\subset\! J_F$. 

Moreover $ \gr^{|p(J)|}\rH^{d}(\overline{M}_{/\overline{\Q}_p},\overline{\cV}_n\otimes 
\Omega^\bullet_{\overline{M}}(\dlog\,\infty))= \rH^{d-|J|}(\overline{M}_{/\overline{\Q}_p},
\underline{\omega}^{-\epsilon_J(n+t)+t}\otimes \underline{\nu}^{p(J)})$.

It is enough to see that the  $\overline{\Q}_p$-vector space   
$\rH^{d-|J|}(\overline{Y}_{\overline{\Q}_p},
\underline{\omega}^{-\epsilon_J(n+t)+t}\otimes \underline{\nu}^{p(J)})
[f]$ is of dimension $1$, for all $J\!\subset\! J_F$.

Because of the existence of a  BGG complex over $\overline{\Q}$ 
 giving by base change the BGG complexes  over $\overline{\Q}_p$ and 
$\C$, we have an Hecke-equivariant isomorphism
$$\rH^{d-|J|}(\overline{Y}_{\overline{\Q}_p},
\underline{\omega}^{-\epsilon_J(n+t)+t}\otimes 
\underline{\nu}^{p(J)})  \otimes_{\overline{\Q}_p} {\C}=
\rH^d(Y^{\an }, \bV_n(\C))[\widehat{{\epsilon}_J\!\otimes\! c}].$$

Finally,   the 
$f$-part  $\rH^d(Y^{\an }, \bV_n(\C))[\widehat{{\epsilon}_J\!\otimes\! c},f]$ is 
equal to $\rH^d_{!}(Y^{\an }, \bV_n(\C))[\widehat{{\epsilon}_J\!\otimes\! c},f]$
and is therefore one dimensional by (\ref{ESH-epsilon}), for all $J\!\subset\! J_F$. 
\hfill $\square$

\begin{rque}
1) We proved that  $W_f$ is pure of weight $d(k_0-1)$. 
The set  of its   Hodge-Tate weights is stable by the symmetry
$h \mapsto d(k_0-1)-h$, corresponding  $|p(J)|\mapsto |p(J_F\!\bs\!J)|$. 
This symmetry is induced by the  Poincar{\'e} duality 
 $W_f\times W_f\rightarrow \overline{\Q}_p(-d(k_0-1))$.

2) If $F$ is a real quadratic field and  $\tau$ denotes the non-trivial embedding of $F$, then
  the Hodge-Tate weights of   $W_f$ are given by $m_\tau, k_0-m_\tau-1, k_0+m_\tau-1, 2k_0-m_\tau-2$.
\end{rque}

\subsection{Hodge-Tate weights of $\rho$ in the crystalline case.}
\label{HT-crys}
The embedding $\iota_p:\overline{\Q}\hookrightarrow \overline{\Q}_p$
allows us to identify $J_F$ with
$\Hom_{\Q\mathrm{-alg.}}(F,\overline{\Q}_p)$.
For each prime $\gp$ of $F$ dividing $p$, we put 
$J_{F,\gp}=\Hom_{\Q_p\!\mathrm{-alg.}}
(F_{\gp},\overline{\Q}_p)$. Thus we get a partition 
$J_F=\underset{\gp}{\coprod}J_{F,\gp}$.
Let $D_{\gp}$ (resp. $I_{\gp}$) be a decomposition (resp.  inertia) subgroup of $\G_F$  at $\gp$.

The following result is due to Wiles  if $k$ is parallel, and 
to Hida  in the  general case. 

\begin{theo} {\rm (Wiles \cite{wiles}, Hida \cite{hida-jami})}\label{ord}
Assume that  $f$ is  ordinary  at $p$ (see  Def.\ref{ord-autom}).

Then  $\rho_{|D_{\gp}}$ is reducible and  :
$$\hspace{-7cm}\mathrm{\bf (ORD)}\hspace{3cm} \rho_{f|I_{\gp}}\sim
\begin{pmatrix} \varepsilon_{\gp}  & * \\
0  & \delta_{\gp} \end{pmatrix},$$
where $\delta_{\gp}$ (resp. $\varepsilon_{\gp}$) is obtained  by  composing the  
class field theory map  $I_{\gp}\rightarrow \go_{\gp}^\times$
with the map  $\go_{\gp}^\times \rightarrow  \overline{\Q}_p^\times$, \enspace
 $x\mapsto \underset{\tau\in J_{F,\gp} }{\prod}\tau(x)^{-m_\tau}$
(resp. $x\mapsto \underset{\tau\in J_{F,\gp} }{\prod}\tau(x)^{-(k_0-m_\tau-1)}$). 
\end{theo}

\medskip
 Breuil \cite{breuil} has shown that  if  $p>k_0$ and  
$p$ does not divide $\Delta$, then $\rho$ 
is crystalline at each prime $\gp$ of $F$ dividing $p$, 
with   Hodge-Tate weights  between $0$ and $k_0-1$. 

\begin{cor} Assume  $p>k_0$ and  that  $p$ does not divide $\Delta$. Then 
for each prime $\gp$ of $F$ dividing $p$, $\rho_{|D_{\gp}}$ is crystalline with 
 Hodge-Tate weights the  $2[F_{\gp}:\Q_p]$ integers 
$(m_{\tau},k_0-m_{\tau}-1)_{\tau\in J_{F,\gp}}$.
\end{cor}

{\bf Proof : }
Assume first  that $n\neq 0$. Let $K$ be a CM  quadratic extension   of $F$,
 in which each   prime $\gp$ of $F$ 
(dividing $p$) splits as $\gp=\mathfrak{P}\mathfrak{P}^c$.
Blasius and Rogawski \cite{BR} have then  constructed  a pure motive   
over  $K$ with coefficients in $E$, of Hodge weights 
$(m_{\tau},k_0-m_{\tau}-1)_{\tau\in J_F}$
and  whose $p$-adic  realization  is isomorphic  to the 
restriction of $\rho$ to $\G_K$. This shows that $\rho_{|D_{\gp}}$ 
is de Rham, and even   crystalline for   $p$ big enough. 

By Faltings' Comparison  Theorem the Hodge  weights  of this   
motif correspond naturally  
(via   $\iota_p:\overline{\Q} \hookrightarrow \overline{\Q}_p$) 
to the  Hodge-Tate weights of its  $\mathfrak{P}$-adic realization,
 which are the same as the Hodge-Tate weights of $\rho$. This proves
 the corollary for  $n\neq 0$.

If $n=0$ (or more generally if  $k$ is parallel) we
can complete the proof using the   following
\begin{lemma} 
Let $a$ and $b$ be two positive integers  and let 
$(a_{\tau})_{\tau\in J_F}$ (resp.  
$(b_{\tau})_{\tau\in J_F}$)  be  integers  
satisfying  $0\leq 2a_{\tau}<a$ (resp. $0\leq 2b_{\tau}<b$).
Assume that the following two sets are equal (with multiplicities)

$\left\{\sum_{\tau\in J} a_{\tau}+\sum_{\tau\in J_F\bs J} 
(a-a_{\tau}) \text{ , } J\!\subset\! J_F\right\}=
\left\{ \sum_{\tau\in J} b_{\tau}+\sum_{\tau\in J_F\bs J} 
(b-b_{\tau}) \text{ , } J\!\subset\! J_F\right\}.$

 Then $a=b$ and  we have equality (with multiplicities)
 $\{a_{\tau} \text{ , } \tau\in J_F\}=
\{b_{\tau} \text{ , } \tau\in J_F\}$.
\end{lemma}

Using this lemma together with  Thm.\ref{BL-thm} and  Cor.\ref{hodge-complex}(ii) 
we obtain the Hodge-Tate weights of $\rho$  at the primes   $\gp$
dividing $p$,  up to permutation.
In particular, we know exactly the Hodge-Tate weights of $\rho$ when
$k$ is parallel. \hfill $\square$

\subsection{Fontaine-Laffaille weights of $\overline{\rho}$ 
in the crystalline case.} 

Our aim is to  find the weights of   $\overline{\rho}_{|I_{\gp}}$  for 
$\gp$ dividing  $p$. If $f$ is ordinary at $p$ we know  by 
Thm.\ref{ord}, that   $\rho_{|D_{\gp}}$ is reducible  and  by a simple 
 reduction modulo $\cP$ we obtain the weights of  $\overline{\rho}_{|I_{\gp}}$.

\begin{prop}  Assume  $p>k_0$ and  that  $p$ does not divide $\Delta$.
 Then  $\overline{\rho}$ is  crystalline  at each $\gp$
dividing $p$ with  Fontaine-Laffaille weights  $(m_{\tau},k_0-m_{\tau}-1)_{\tau\in J_{F,\gp}}$.
\end{prop}

{\bf Proof  : } It is a consequence of the theory of Fontaine and Laffaille
\cite{FoLa},  and of the computation of  Hodge-Tate weights of 
 $\rho_{|D_{\gp}}$ in   {\S}\ref{HT-crys}.

Consider a Galois stable lattice $\cO^2$ in the crystalline representation 
$\rho$, as well as the  sub-lattice   $\cP^2$. The representation 
 $\overline{\rho}$  is equal to the  quotient  of these  two lattices. 
It is crystalline, as a sub-quotient of a  crystalline representation.
Its     weights are determined by the associated 
filtered  Fontaine-Laffaille module. Since  the Fontaine-Laffaille's 
functor is exact, it is a  quotient of the  Fontaine-Laffaille's filtered modules 
 associated to the two lattices. By compatibility 
of the filtrations on these two lattices, and by the condition $p>k_0$, 
the graded of the quotient have  the right   dimension. See the 
Appendix for more details. \hfill$\square$

\begin{cor} \label{poids-modere}
Let $\gp$ be a  prime of $F$ above  $p$.
 Then
$$\overline{\rho}_{|I_{\gp}}\sim
\begin{pmatrix} \varepsilon_{\gp}  & * \\
0  & \delta_{\gp} \end{pmatrix},$$
\noindent
where $\varepsilon_{\gp},\delta_{\gp}:
I_{\gp}\rightarrow 
\overline{\F}_p^\times $ are two tame characters of 
level $|J_{F,\gp}|$, whose 
product equals the $(1-k_0)$th power  of the modulo $p$ cyclotomic character 
and whose sum has   Fontaine-Laffaille weights 
 $(m_{\tau},k_0-m_{\tau}-1)_{\tau\in J_{F,\gp}}$.
\end{cor}

\section{Study  of the images  of $\overline{\rho}$
and   $\Ind_F^{\Q}\overline{\rho}$.}

In all this section we assume that $p \!>\! k_0$ and 
that $p$ does not divide $6\Delta$.
 
Let  $\omega:\G_{\Q}\rightarrow \F_p^\times$ be the   modulo $p$
cyclotomic character and let $\pr:\GL_2(\kappa)\rightarrow \PGL_2(\kappa)$
be  the canonical  projection.  

\subsection{Lifting of characters and irreducibility criterion 
for $\overline{\rho}$.}

\begin{prop}\label{Irr}
{\rm(i)} For all but finitely many primes $p$ $\mathrm{\bf (Irr_{\overline{\rho}})}$ holds, that is 
 $\overline{\rho}=\overline{\rho}_{f,p}$  is absolutely irreducible.

{\rm(ii)}  Assume that $k$  is non-parallel.  If  for all  $J\!\subset\! J_F$
there exists  $\epsilon\in\go_+^{\times}$,
$\epsilon-1\in\gn$ such that  $p$   does not divide  the non zero 
integer  $\Nm(\epsilon^{p(J)}\!-\!1)$, then  $\mathrm{\bf (Irr_{\overline{\rho}})}$ holds.
\end{prop}

\begin{rque} Assume that   $k=k_0t$ is parallel and that for all 
$\varnothing\!\subsetneq\!J\!\subsetneq\! J_F$, there exists  $\epsilon\in\go_+^{\times}$,
$\epsilon-1\in\gn$ such that  $p$   does not divide  the non zero 
integer  $\Nm(\epsilon^{p(J)}\!-\!1)$. Then 
we expect  $\overline{\rho}=\overline{\rho}_{f,p}$ to be   absolutely
irreducible, unless $p$
 divides 
the  constant term  of an Eisenstein   series  of weight $k$ and level
dividing $\gn$, that is the numerator of the value  at $1-k_0$ of the  $L$-function of
a  finite order Hecke character  of  $F$ of conductor dividing
$\gn$ (see \cite{FaJo}\S3.2 for the case $F=\Q$). 
 \end{rque}

{\bf Proof : } As  $\overline{\rho}$ is  totally odd, 
if it is  irreducible, then it is  absolutely irreducible.
 Assume   that   $\overline{\rho}$  is reducible :  
$\overline{\rho}^{\mathrm{s.s.}}=\varphi_{\gal}\oplus
 \varphi_{\gal}'$. The characters  $\varphi_{\gal},\varphi_{\gal}':
\G_F\rightarrow \kappa^\times$ are unramified outside $\gn p$ and  
$\varphi_{\gal}\varphi_{\gal}'=\det(\overline{\rho})=\overline{\psi}_{\gal}
\omega^{-1}$  (recall that  $\psi$ is a Hecke character of infinity type $-n_0t$). 
Denote by $\widehat{\go}^{\times}_{\gn,1}$
the subgroup of $\widehat{\go}^{\times}$ of elements 
$\equiv 1\pmod{\gn}$. Then $\widehat{\go}^{\times}_{\gn,1}$ is a product of
its $p$-part $\prod_{\gp\mid p}\go_{\gp}^\times$
and its part  outside $p$, denoted by  $\widehat{\go}^{\times(p)}_{\gn,1}$.

By the global class field theory, the Galois group of  the  maximal 
$\gn$-ramified (resp. $\gn p^\infty$-ramified)
abelian extension  of $F$ is isomorphic to 
$\Cl^+_{F,\gn}=\A_F^\times/ F^\times 
\widehat{\go}^{\times}_{\gn,1}  D(\R)_+$
(resp.  $\Cl^+_{F,\gn p^\infty}:=\limproj 
\Cl^+_{F,\gn p^r}=\A_F^\times/\overline{ F^\times 
\widehat{\go}^{\times(p)}_{\gn,1}D(\R)_+}$). 
We choose the convention in which an uniformizer corresponds to a 
geometric Frobenius. We have the following  exact sequence  
\begin{equation}\label{Irr2}
1 \rightarrow 
(\prod_{\gp\mid p}\go_{\gp}^\times)/
\overline{\{\epsilon\in\go_+^\times | \epsilon-1\in \gn\}}
\rightarrow  \Cl^+_{F,\gn p^\infty}
\rightarrow  \Cl^+_{F,\gn}\rightarrow 1.
\end{equation}
By  Cor.\ref{poids-modere}, for each $\gp\mid
 p$, $\varphi_{\gal}\oplus  \varphi_{\gal}'$ is  crystalline at  $\gp$ 
of  weights $(m_{\tau},k_0\!-\!m_{\tau}\!-\!1)_{\tau\in J_{F,\gp}}$.

By (\ref{Irr2}) for every  $\epsilon\in\go_+^\times $, 
$\epsilon-1\in\gn$ we have the following equality in $\kappa$ :
\begin {equation*}
1=\varphi_{\gal}(\epsilon)=\prod_{\gp\mid p}
\varphi_{\gal,\gp}(\epsilon)=
\underset{\gp\mid p}{\prod}\enspace
\underset{\tau\in J_{F,\gp}}{\prod}\tau(\epsilon)^{m_\tau
\text{ or } (k_0-m_\tau-1)}=\epsilon^{p(J)},
\end{equation*}
for some subset $J\subset J_F$. Because of the assumption $p>k_0$,
if $k$ is  non-parallel, then  $\epsilon^{p(J)}\!\neq\!1$ for all
 $J\subset J_F$. Thus   we obtain (ii)  and (i)  when 
$k$ is  non-parallel.

\bigskip
Assume now that  $k=k_0t$ is parallel and that for all 
$\varnothing\!\subsetneq\!J\!\subsetneq\! J_F$, there exists  $\epsilon\in\go_+^{\times}$,
$\epsilon-1\in\gn$ such that  $p$   does not divide  the non zero 
integer  $\Nm(\epsilon^{p(J)}\!-\!1)$. The same arguments
 as above show that the  restriction  to 
$\prod_{\gp\mid p}\go_{\gp}^\times$ 
 of the character  $\varphi_{\gal}$ (resp. $\varphi_{\gal}'$) 
$ \Cl^+_{F,\gn p^\infty} \rightarrow \kappa^\times$ 
is trivial (resp.  given by the  $(1-k_0)$-th  power of the norm).
By the following lemma (applied to $P=\Cl^+_{F,\gn p^\infty}$ and 
$Q=(\prod_{\gp\mid p}\go_{\gp}^\times)/
\overline{\{\epsilon\in\go_+^\times | \epsilon-1\in \gn\}} $)  there exists an  unique character
 $\widetilde{\varphi}_{\gal}$ (resp. $\widetilde{\varphi}_{\gal}'$)
$:\Cl^+_{F,\gn p^\infty} \rightarrow \mathcal{O}^\times$
lifting  $\varphi_{\gal}$ (resp. $\varphi_{\gal}'$) and  whose  restriction to 
$\prod_{\gp\mid p}\go_{\gp}^\times$ is also
trivial (resp.   given by the $(1-k_0)$-th  power of the norm).

\begin{lemma} \label{lifting}
Let $P$ be an abelian group  and  $Q$ be a subgroup, such that 
the  factor  group $P/Q$ is finite. Let $\varphi_{P}:P \rightarrow 
\kappa^\times$  and $\widetilde{\varphi_{Q}}:Q \rightarrow 
\mathcal{O}^\times$  be two  characters such  that 
$\varphi_{P}|_{Q}= \widetilde{\varphi_{Q}}\mod{p}$. Then, there  exists an unique 
character  $\widetilde{\varphi_{P}}:P \rightarrow \mathcal{O}^\times$, 
whose restriction to $Q$ is $\widetilde{\varphi_{Q}}$ and such that 
$\widetilde{\varphi_{P}}\mod{p}=\varphi_{P}$.
\end{lemma}

For $x\in \A_F^\times$, we put   $\varphi(x):= \widetilde{\varphi}_{\gal}(x)$ and 
 $\varphi'(x):= \widetilde{\varphi}_{\gal}'(x)x_p^{-k}
x_{\infty}^{k}$.  Then $\varphi$ (resp. $\varphi'$)  is a  
 Hecke character of $F$, of conductor dividing $\gn$ 
and infinity type $0$ (resp. $(1-k_0)t$).  It is crucial to observe that
 there are only finitely many such  $\varphi$ and  $\varphi'$.

Assume now that for infinitely many primes $p$,  
$\overline{\rho}$   is reducible. Then there exist 
Hecke characters $\varphi$ and  $\varphi'$ as above, such that 
for infinitely many  primes $p$ we have  $\overline{\rho}^{\mathrm{s.s.}}\equiv
\varphi_{\gal}\oplus\varphi_{\gal}'\pmod{\cP}$.
Hence for each prime $v$ of $F$ not dividing $\gn$ we have  
$c(f,v)\equiv\varphi(\varpi_{v})+\varphi'(\varpi_{v})\pmod{\cP}$ 
for infinitely many  $\cP$'s and  hence
$c(f,v)=\varphi(\varpi_{v})+\varphi'(\varpi_{v})$. By the Cebotarev
Density Theorem  we obtain
$\rho^{\mathrm{s.s.}}=\varphi\oplus\varphi'$.
This  contradicts the absolute irreducibility of $\rho$ (see \cite{Ta2}). \hfill$\square$

\subsection{The  exceptional  case.}\label{exceptional}

The aim of this paragraph is to find a bound for the primes $p$ such
that  $\pr(\overline{\rho}_{f,p}(\G_F))$ is  isomorphic to one  of the groups 
$A_4$, $S_4$ or $A_5$. We will only use the fact that 
 the elements of these groups   are of order  at most  $5$.

Assume that $\pr(\overline{\rho}_{f,p}(\G_F))\cong A_4, S_4$ or $A_5$.
By Cor.\ref{poids-modere} there exist
$\epsilon_\tau\!\in\!\{\pm 1\}$, $\tau\in J_F$, such that   for all 
$\gp\mid p$ and  for any generator $x$ of 
$\F_{p^h}^\times$, where $h=|J_{F,\gp}|$,  the element
$$\prod_{\tau\in \Gal(\F_{p^h}\!/\!\F_{p})}\!\!\!
\tau(x)^{\epsilon_\tau(k_\tau-1)}\in \F_{p^h}^\times$$
\noindent belongs  to  $\pr(\overline{\rho}(I_{\gp}))$ 
and    is therefore of order   at most  $5$  (if {\bf (ORD)} holds
 we may   assume that   $\epsilon_\tau=1$ for all $\tau$). Denote by 
$\tau_1,...,\tau_{h}$ the elements of $J_{F,\gp}$. Then 
$$\epsilon_{\tau_1}(k_{\tau_1}-1)+\epsilon_{\tau_2} p(k_{\tau_2}-1)+...+
\epsilon_{\tau_{h}}  p^{h-1}(k_{\tau_{h}}-1)\in \Z/(p^h-1)$$
 is of order  $\leq 5$, hence 
$5((k_{\tau_1}-1)+p(k_{\tau_2}-1)+...+p^{h-1}
(k_{\tau_{h}}-1))\geq p^h-1$.

If  we replace the generator $x$  by $x^p,x^{p^2},...,x^{p^{h-1}}$
and then sum these  inequalities we find 
$5\underset{\tau\in J_{F,\gp}}{\sum}(k_\tau-1)
\geq |J_{F,\gp}|(p-1)$. Hence 
 $5\underset{\tau\in J_F}{\sum}(k_\tau-1) \geq d(p-1)$. 

We conclude that $\pr(\overline{\rho}(\G_F))$ cannot be 
isomorphic to $A_4$, $S_4$ or $A_5$ if 

$$d(p-1)>5\sum_{\tau\in J_F}(k_\tau-1).$$

Note  that this assumption follows from  {\bf (II)} if  $d\geq 5$.

\subsection{The dihedral case.}\label{dihedral}

In this paragraph we study the case when  $\pr(\overline{\rho}_{f,p}(\G_F))$
 is isomorphic to the dihedral group
 $D_{2n}$, where $n\geq 3$ is an  integer  prime to  $p$. Let
$C_n$ be  the cyclic subgroup  of order  $n$ of $D_{2n}$.
Since  $\pr^{-1}(C_n)$ is a commutative group containing only semi-simple
 elements ($p$ does not divide $n$), it is diagonalizable.
 Since  $\pr^{-1}(D_{2n}\bs C_n)$ is 
contained in the normalizer of $\pr^{-1}(C_n)$, it
is contained in the set of anti-diagonal   matrices. 

Let $\varepsilon:D_{2n}\rightarrow \{\pm 1\}$ be the signature map and 
let $K$ be the fixed field  of $\ker(\varepsilon\circ \pr 
\circ \overline{\rho})$. The extension $K/F$ is  quadratic  and 
 unramified outside $\gn p$.

\medskip
Let $c$ be the  non-trivial  element of the Galois group  
$\Gal(K/F)$. As $\overline{\rho}$ is absolutely irreducible, 
but  $\overline{\rho}_{|\G_K}$ is not, there exists a 
character $\varphi_{\gal}:\G_K \rightarrow \kappa^\times$
distinct from  its Galois conjugate  $\varphi_{\gal}^c$ and such that 
$\overline{\rho}_{|\G_K}=\varphi_{\gal}\oplus \varphi_{\gal}^c$.

\begin{lemma}  Let  $\gp$ be a  prime of  $F$ dividing $p$.
Assume $p\neq 2k_\tau-1$,  for  $\tau\in J_{F,\gp}$. Then

\rm{(i)}  the field $K$ is  unramified at  $\gp$,

\rm{(ii)} the prime  $\gp$ splits in $K$ as $\mathfrak{P}\mathfrak{P}^c$, 
and  $\varphi_{\gal}$  is crystalline at $\mathfrak{P}$ (resp. $\mathfrak{P}^c$) 
of weights $(p_\tau)_{\tau\in J_{F,\gp}}$ (resp. 
$(q_\tau)_{\tau\in J_{F,\gp}}$), where $\{p_\tau,q_\tau\}=
\{m_{\tau},k_0-m_{\tau}-1\}$ for each  $\tau\in J_{F,\gp}$.
\end{lemma}

{\bf Proof : } (i)
Otherwise $\overline{\rho}(I_{\gp})$ would contain at least 
one anti-diagonal matrix and the  basis vectors would not be 
eigen for  $\overline{\rho}(I_{\gp})$. But the group 
$\overline{\rho}(I_{\gp})$ has  a common eigenvector. 
Hence, the elements  of  $\pr( \overline{\rho}(I_{\gp}))$ 
would be  of order $\leq 2$. Using the computations of 
\S\ref{exceptional} and $p>k_0$, we find that for all 
$\tau\in J_{F,\gp}$ we have  $2(k_\tau-1)= p-1$.
Contradiction.

(ii) By  Cor.\ref{poids-modere}, $\varphi_{\gal}\oplus\varphi_{\gal}^c$
is crystalline of  weights
$(m_{\tau},k_0-m_{\tau}-1)_{\tau\in J_{F,\gp}}$.
For each prime $\mathfrak{P}$ of $K$ dividing $\gp$, we have
 $(\varphi_{\gal}\varphi_{\gal}^c)_{|I_{\mathfrak{P}}}=
\omega^{1-k_0}_{|I_{\mathfrak{P}}}$. Since $m_{\tau}<k_0-m_{\tau}-1<p-1$
we deduce that ${\varphi_{\gal}}_{|I_{\mathfrak{P}}}\neq
{\varphi_{\gal}^c}_{|I_{\mathfrak{P}}}\cong
{\varphi_{\gal}}_{|I_{\mathfrak{P}^c}}$ and so $\gp$ splits in $K$. 
\hfill $\square$

\medskip

Let $\gO$ be the integer ring of $K$, and 
 $\widehat{\gO}$ its profinite completion.
Denote by $\widehat{\gO}^{\times}_{\gn,1}$
the subgroup of $\widehat{\gO}^{\times}$ of elements 
$\equiv 1\pmod{\gn}$. Then 
$\widehat{\gO}^{\times}_{\gn,1}$ is a product of
its $p$-part $\prod_{\mathfrak{P}\mid p}\gO_{\mathfrak{P}}^\times$
and its part  outside $p$, denoted by 
$\widehat{\gO}^{\times(p)}_{\gn,1}$.

By the global class field theory, the Galois group of  the  maximal 
$\gn$-ramified (resp. $\gn p^\infty$-ramified)
abelian extension  of $K$ is isomorphic to 
$\Cl_{K,\gn}:=\A_K^\times/ K^\times 
\widehat{\gO}^{\times}_{\gn,1}  K_{\infty}^\times$
(resp.  to $\Cl_{K,\gn p^\infty}:=\A_K^\times/\overline{ K^\times 
\widehat{\gO}^{\times(p)}_{\gn,1}K_{\infty}^\times}$). 
 We have the following  exact sequence  :
\begin{equation}\label{CM}
1 \rightarrow 
(\prod_{\mathfrak{P}\mid p}\gO_{\mathfrak{P}}^\times)/
\overline{\{\epsilon\in\gO^\times | \epsilon-1\in \gn\}}
\rightarrow  \Cl_{K,\gn p^\infty}
\rightarrow  \Cl_{K,\gn}\rightarrow 1
\end{equation}

\begin{prop} \label{not-CM}
{\rm (i)} Assume that  for all $\tau\in J_F$, $p\neq 2k_\tau-1$ and 
that $\pr(\overline{\rho}(\G_F))$ is dihedral. Let $K/F$ be the 
quadratic extension defined above. Then 

 $\relbar$ either $K$ is  CM and  there exists 
a  Hecke character   $\varphi$ of  $K$ of conductor of norm
dividing $\gn\Delta_{K/F}^{-1}$ and  infinity type 
$(m_\tau,k_0\!-\!1\!-\!m_\tau)_{\tau\in J_F}$ such that $\rho\equiv \Ind_{K}^{F}
\varphi \pmod{\cP}$,

 $\relbar$ either $K$ is not  CM and we can choose  
places $\widetilde{\tau}$ of $K$ above each $\tau\in J_F$
such that for all $\epsilon\in\gO^{\times}$, $\epsilon-1\in\gn$ the prime  $p$  divides 
 $\displaystyle \N_{K/\Q}\left(\prod_{\tau\in J_F}
\widetilde{\tau}(\epsilon)^{m_{\tau}}\widetilde{\tau}(c(\epsilon))^{k_0-m_{\tau}-1}\!-\!1\right)$. 

 {\rm (ii)} Assume that   $f$ is not a theta series. Then for 
all but finitely many primes $p$ the group
 $\pr(\overline{\rho}_{f,p}(\G_F))$ is not dihedral.
\end{prop}

\begin{rque}
(i) The  primes $p$ for which the   congruence $\rho\equiv \Ind_{K}^{F} 
\varphi\pmod{\cP}$ may occur should  be  controlled by the special value  
 of the $L$-function  associated to the  CM character   $\varphi/\varphi^c$ 
(in the  elliptic case it is proved by  Hida \cite{hida64} and 
Ribet \cite{ribet4}; see also  Thms A and B). 

(ii) We would like to thank E. Ghate for having pointed us out the
possible existence of dihedral primes for non CM fields $K$. It
would be interesting to explore the converse statement, that is to say 
to try to construct,  for a given prime $p$ dividing the above norms,  a newform $f$ 
such that $\pr(\overline{\rho}_{f,p}(\G_F))$ is dihedral. 
\end{rque}

{\bf Proof : } (i) By (\ref{CM}) and the above lemma, we have 
$\varphi_{\gal}: \Cl_{K,\gn p^\infty} \rightarrow \kappa^\times$
whose restriction to 
$\prod_{\mathfrak{P}\mid p}\gO_{\mathfrak{P}}^\times $
is given by the reduction modulo $p$ of   an algebraic character
$x\mapsto x^{\widetilde{k}}$, where  $\widetilde{k}=\sum_{\tau\in J_F}m_{\tau}
\widetilde{\tau}+(k_0-m_{\tau}-1)\widetilde{\tau}\circ c$, for some
choice of places $\widetilde{\tau}$ of $K$ above $\tau\in J_F$. 

We observe that the  character   $x\mapsto x^{\widetilde{k}}$  is
trivial on $\go_+^{\times}$,  whereas 
it is only trivial modulo $p$ on  $\{\epsilon\in\gO^\times | \epsilon-1\in \gn\}$. 
The case  when $K$ is not  CM follows imediately.

Assume now that $K$ is CM. In this case 
$\{\epsilon\in\go_+^\times | \epsilon-1\in \gn\}$ is 
a finite index subgroup of $\{\epsilon\in\gO^\times | \epsilon-1\in
\gn\}$. Since $\ker(\cO^\times\rightarrow \kappa^\times)$ does not
contain elements of finite order, the above  character is 
trivial on $\{\epsilon\in\gO^\times | \epsilon-1\in \gn\}$.

By the  lemma \ref{lifting} (applied to $P=\Cl_{K,\gn p^\infty}$ and 
$Q=(\prod_{\mathfrak{P}\mid p}\gO_{\mathfrak{P}}^\times)/
\overline{\{\epsilon\in\gO^\times | \epsilon-1\in \gn\}}$) there exists a lift  
$\widetilde{\varphi}_{\gal}:\Cl_{K,\gn p^\infty} \rightarrow
\cO^\times $ whose restriction to 
$\prod_{\mathfrak{P}\mid p}\gO_{\mathfrak{P}}^\times$ is given by 
$x\mapsto x^{\widetilde{k}}$. 

We  put    $\varphi(x):= \widetilde{\varphi}_{\gal}(x)x_p^{-\widetilde{k}}
x_{\infty}^{\widetilde{k}}$.  Then $\varphi$   is a  Hecke character of $K$ as desired.

(ii) There are finitely many fields $K$ as above. For those $K$ that are
not CM it is enough to choose $\epsilon\in\gO^\times$, $\epsilon-1\in
\gn$
of infinite order in $\gO^\times\!/\!\go^\times$. 

 For each of the CM fields  $K$ that are only finitely many
 characters $\varphi$ as above.
Therefore, if  $\pr(\overline{\rho}(\G_F))$ is dihedral
for infinitely many  primes $p$, then there would exist $K$ and 
$\varphi$ as above, such that the congruence 
$\rho\equiv \Ind_{K}^{F} \varphi \pmod{\cP}$ happens for  infinitely many  
$\cP$'s. Hence  $f$ would be  
equal to the theta series associated to $\varphi$. \hfill $\square$

\subsection{The image  of $\overline{\rho}$ is ``large''.}

\begin{theo}\label{dickson} $\mathrm{(Dickson)}$
{\rm(i)} An  irreducible subgroup  of  $\PSL_2(\kappa)$  of order divisible 
by  $p$ is  conjugated   inside  $\PGL_2(\kappa)$  to $\PSL_2(\F_{q})$
 or to $\PGL_2(\F_{q})$,  for some  power $q$  of $p$.

{\rm(ii)} An  irreducible subgroup  of  $\PSL_2(\kappa)$  of order prime to 
$p$  is either dihedral, either isomorphic to one of the 
groups $A_4$, $S_4$ or $A_5$.
\end{theo}

As an application of this theorem, Prop.\ref{Irr}, 
Prop.\ref{not-CM} and {\S}\ref{exceptional} we obtain  the 
following 

\begin{prop}\label{LI}
 Assume that $f\in S_k(\gn,\psi)$ is a newform, which is 
not a theta series. Then for  all but finitely many primes $p$,
the image of the $p$-adic representation 
$\overline{\rho}$ associated to $f$ is  large, in the following sense :

\smallskip
$\mathrm{\bf (LI_{\overline{\rho}})}$ there exists a power $q$ of $p$ 
such that $\SL_2(\F_q)\!\subset\! 
\im(\overline{\rho}) \!\subset\! \kappa^\times \GL_2(\F_q)$.
\end{prop}

\medskip
Let $\widehat{F}$ be the compositum of   $\widetilde{F}$
and  of the subfield of $\overline{\Q}$ fixed  by the Galois group 
$\ker(\overline{\psi}\omega^{n_0})$.
The extension $\widehat{F}/F$ is Galois and unramified at $p$,
because $\widetilde{F}$  is unramified at $p$ and $\psi$ is of conductor 
prime to $p$. Therefore  $\G_{\widehat{F}}$ is a normal 
subgroup  of $\G_F$  containing  the inertia subgroups  
$I_{\gp}$, for all $\gp$ dividing $p$.

We put $\mathcal{D}=\det(\overline{\rho}(\G_{\widehat{F}}))=
(\F_p^\times)^{1-k_0}$.

\begin{prop}\label{oneform} Assume $\mathrm{\bf (LI_{\overline{\rho}})}$.
Then  there exists a power $q$ of $p$ such that,

either   $\overline{\rho}(\G_{\widehat{F}})=\GL_2(\F_q)^\mathcal{D}:=
\{x\in \GL_2(\F_q)\enspace |\enspace  \det(x)\in \mathcal{D}\}$,

either  $\overline{\rho}(\G_{\widehat{F}})=
(\F_{q^2}^\times \GL_2(\F_q))^\mathcal{D}:=
\{x\in\F_{q^2}^\times \GL_2(\F_q) 
\enspace|\enspace \det(x)\in \mathcal{D}\}$ .
\end{prop}

{\bf Proof : }  We first show that  $\pr(\overline{\rho}(\G_{\widehat{F}}))$
is still  irreducible of order  divisible by  $p$.  By 
$\mathrm{\bf (LI_{\overline{\rho}})}$  the group $\pr(\overline{\rho}(\G_F))$ 
 is isomorphic  to   $\PSL_2(\F_q)$ or $\PGL_2(\F_q)$. The group 
$\pr(\overline{\rho}(\G_{\widehat{F}}))$
is a   non-trivial normal subgroup of $\pr({\im}(\overline{\rho}))$ 
(because it  contains  $\pr(\overline{\rho}(I_p))$ and  $p\!>\!k_0$; 
see  Cor.\ref{poids-modere}). As $\PSL_2(\F_q)$ is a simple group of
 index   $2$ in  the group $\PGL_2(\F_q)$, we deduce 
$$\PSL_2(\F_q)\subset \pr(\overline{\rho}(\G_{\widehat{F}}))\subset
 \pr(\overline{\rho}(\G_F)) \subset\PGL_2(\F_q).$$

\begin{lemma}\label{sl2}
Let $H$ be a group of center $Z$ and let $\pr:H\rightarrow H/Z$ the 
canonical projection. Let $P$ and $Q$ be two subgroups of $H$ such that 
 $\pr(P)\supset \pr(Q)$. Assume moreover that  $Q$ does not have 
non-trivial abelian quotients. Then $P\supset Q$.
\end{lemma}

It follows from this lemma that $\overline{\rho}(\G_{\widehat{F}})\supset
\SL_2(\F_q)$, hence 
 $$(\kappa^\times \GL_2(\F_q))^\mathcal{D}\supset 
 \overline{\rho}(\G_{\widehat{F}}) \supset \GL_2(\F_q)^\mathcal{D}.$$

Since  $[(\kappa^\times \GL_2(\F_q))^\mathcal{D}:
\GL_2(\F_q)^\mathcal{D}]\leq 2$ we're done . \hfill$\square$

\bigskip
Let $y\in \F_{q^2}\backslash \F_q$ be  such that   $y^2\in \F_q$. 
Then  $(\F_{q^2}^\times \GL_2(\F_q))^\mathcal{D}=
\GL_2(\F_q)^\mathcal{D}\amalg(y\GL_2(\F_q))^\mathcal{D}$ and  hence
$\tr((\F_{q^2}^\times \GL_2(\F_q))^\mathcal{D})= \F_q\cup y\F_q$.
Therefore,  the $\F_p$-algebra generated by the traces of the 
elements of $(\F_{q^2}^\times \GL_2(\F_q))^\mathcal{D}$ is 
 $\F_{q^2}$, while 
 $\pr((\F_{q^2}^\times \GL_2(\F_q))^\mathcal{D})\subset \PGL_2(\F_q)$.
This reflects the existence of a congruence with a form having  inner twists.

\subsection{ The image  of $\Ind_F^{\Q}\overline{\rho}$ is 
 ``large''.}\label{ind-image}

We assume in this paragraph that $\mathrm{\bf (LI_{\overline{\rho}})}$ holds.

\smallskip
By Prop.\ref{oneform} there exists a power $q$ of $p$ such that 
$\pr(\overline{\rho}(\G_{\widehat{F}}))=\PSL_2(\F_q)$ or $\PGL_2(\F_q)$. 
Consider the representation  $\pr(\Ind_F^{\Q}\overline{\rho}): 
\G_{\widehat{F}}\rightarrow\PGL_2(\F_q)^{J_F}$.

Any automorphism of the simple  group $\PSL_2(\F_q)$ is the
composition of  a conjugation  by an element of $\PGL_2(\F_q)$ and 
of a Galois automorphism of $\F_q$. By a lemma of Serre 
(see \cite{ribet3}), there exist a partition $J_F=\coprod_{i\in I}J_F^i$
 and for each  $i\in I$, and for each  $\tau\in J_F^i$   an element
  $\sigma_{i,\tau}\in \Gal(\F_q\!/\!\F_p)$ such that 
$$\pr(\phi(\SL_2(\F_q)^I))
\subset \pr(\Ind_F^{\Q}\overline{\rho}(\G_{\widehat{F}}))\subset
  \pr(\phi(\GL_2(\F_q)^I)),$$
$$ {\rm where }\enspace \phi=(\phi^i)_{i\in I} :\GL_2(\F_q)^I
\hookrightarrow  \GL_2(\F_q)^{J_F} \enspace \text{\rm is given by   }\enspace
 \phi^i(M_i)= (M_i^{\sigma_{i,\tau}})_{\tau\in J_F^i }.$$

Keeping these notations, we introduce the following assumption on the image of 
$\Ind_F^{\Q}\overline{\rho}$

\medskip 
\medskip 
$\mathrm{{\bf (LI}_{Ind \overline{\rho}}{\bf )}}$ the condition 
$\mathrm{\bf (LI_{\overline{\rho}})}$ holds and $\forall$
$i\in I$, $\forall$ $\tau,\tau'\in J_F^i$ $(\tau\neq \tau' \Rightarrow 
\sigma_{i,\tau}\neq \sigma_{i,\tau'})$.

\medskip 
We  now  introduce a genericity  assumption on the weight $k$.
\begin{defin}
We say that the weight $k\in \Z[J_F]$ is non-induced, if there does not exist 
a strict subfield $F'$ of $F$ and a weight $k'\in \Z[J_{F'}]$ such that 
for each $\tau\in J_F$, $k_\tau=k'_{\tau|_{F'}}$.
\end{defin}

\begin{rque}\label{non-IND2}
Define $\widetilde{k}=\sum_{\widetilde{\tau}\in J_{\widetilde{F}}}
k_{\widetilde{\tau}}\widetilde{\tau}\in \Z[J_{\widetilde{F}}]$
 by  putting $k_{\widetilde{\tau}}= k_{\widetilde{\tau}|_F}$, for all 
$\widetilde{\tau}\in J_{\widetilde{F}}$. 
The group   $\G_{\Q}$ acts on $\Z[J_{\widetilde{F}}]$ by 
$\widetilde{k}=\sum_{\widetilde{\tau}\in J_{\widetilde{F}}}
k_{\widetilde{\tau}}\widetilde{\tau}\mapsto 
\widetilde{k}^{\widetilde{\tau}'}=\sum_{\widetilde{\tau}\in J_{\widetilde{F}}}
k_{\widetilde{\tau}\widetilde{\tau}'}\widetilde{\tau}$. It is easy to see
that $k\in \Z[J_F]$ is non-induced, if and only if 
 $\{\widetilde{\tau}'\in \G_{\Q} \enspace |\enspace 
 \widetilde{k}=\widetilde{k}^{\widetilde{\tau}'}\}$ equals $\G_F$.
\end{rque}

\begin{prop}\label{k-not-IND}
Assume that $\mathrm{\bf (LI_{\overline{\rho}})}$ hold and $k$ is non-induced.
Assume moreover that for all $\tau\neq\tau'\in J_F$,
 $p\neq k_\tau+k_{\tau'}-1$. 
 Then $\mathrm{{\bf (LI}_{Ind \overline{\rho}}{\bf )}}$ hold.
\end{prop}

{\bf Proof : } Let $\widetilde{\tau}_1, \widetilde{\tau}_2\in \G_{\Q}$
be such that  for all $g\in \G_{\widehat{F}}$ we have 
$\pr(\overline{\rho}(\widetilde{\tau}_1^{-1}g\widetilde{\tau}_1))=
\pr(\overline{\rho}(\widetilde{\tau}_2^{-1}g\widetilde{\tau}_2))$. 
We have to prove that  $\widetilde{\tau}_1^{-1}\widetilde{\tau}_2\in \G_F$.
We put $\overline{\rho}_i(g)=
\overline{\rho}(\widetilde{\tau}_i^{-1}g\widetilde{\tau}_i)$  ($i=1,2$).

Let $\mathfrak{P}$ be  a prime ideal  of $\widehat{F}$ above 
a  prime ideal  $\gp$ of $F$ dividing  $p$. Denote by   $h_i'$ (resp.
 $h_i$) the residual degree of $\mathfrak{P}^{\widetilde{\tau}_i}$ (resp. of 
$\gp^{\tau_i}$), $i=1,2$. By   Cor.\ref{poids-modere} 
we have 
$\overline{\rho}_i|_{I_{\mathfrak{P}}}^{\enspace\mathrm{s.s.}}=
\varepsilon_{i}\oplus \delta_i,$ 
where $\varepsilon_{i}$ (resp. $\delta_i$)  $I_{\mathfrak{P}}
\rightarrow  I_{\mathfrak{P}^{\widetilde{\tau}_i}} \rightarrow 
\F_{p^{h_i'}}^\times \rightarrow \F_{p^{h_i}}^\times \rightarrow 
\kappa^\times$ is obtained by composing 
the conjugation by  $\widetilde{\tau_i}$, the projection on the tame 
inertia, the norm map, and the character 
$x\mapsto\underset{\tau\in J_{F,\gp^{\tau_i}}}\prod 
\tau(x)^{p_\tau}$ (resp. $x\mapsto\underset{\tau\in 
J_{F,\gp^{\tau_i}}} \prod  \tau(x)^{q_\tau}$),  
where $\{p_\tau,q_\tau\}=\{m_\tau,k_0-1-m_\tau\}$.
In the case {\rm(ORD)} we can even assume that for all
$\tau\in J_F$, we have $p_\tau=m_\tau$ and $q_\tau=k_0-1-m_\tau$. 

Note that $\varepsilon_1\delta_1= \varepsilon_2\delta_2=\omega^{1-k_0}$. 
Since  $I_{\mathfrak{P}}\subset \G_{\widehat{F}}$ and 
$\pr\circ\overline{\rho}_1=\pr\circ\overline{\rho}_2$ on $\G_{\widehat{F}}$,
we may assume that  $\varepsilon_1/\delta_1=\varepsilon_2/\delta_2$. 
By varying   $\mathfrak{P}$ 
we deduce that for all $\widetilde{\tau}\in J_{\widetilde{F}}$, 
$\widetilde{k}_{\widetilde{\tau}}=
\widetilde{k}_{\widetilde{\tau}\widetilde{\tau}_1^{-1}\widetilde{\tau}_2}$   
(here we use   $p>k_0$ and the assumption  $p\neq k_\tau+k_{\tau'}-1$). 
As $k$ is non-induced, it follows from the remark \ref{non-IND2} that 
 $\widetilde{\tau}_1^{-1}\widetilde{\tau}_2\in \G_F$. \hfill $\square$

The following corollary  generalizes a  result of Ribet \cite{ribet3} 
on the image of a Galois representation associated to a 
family of classical modular forms, to the case of the 
family of internal conjugates of a   Hilbert modular form.

\begin{cor} Assume that $\mathrm{\bf (LI_{\overline{\rho}})}$ hold 
and $k$ is non-induced. Assume moreover that  $p>2k_0$ is totally 
split in $F$. Then, 
$$(\GL_2(\F_q)^{J_F})^{\mathcal{D}}
\!\subset\! \Ind_F^{\Q}\overline{\rho}(\G_{\widehat{F}})\!\subset\!
(\overline{\rho}(\G_{\widehat{F}})^{J_F})^{\mathcal{D}} \text{ , where } 
\mathcal{D}=\F_p^{\times(1-k_0)}.$$ 
\end{cor}

$${\rm Put } \enspace
 H(\F_q)=\left(\prod_{i\in I}\GL_2(\F_q)\right)^\mathcal{D}:=
\left\{ (M_i)_{i\in I}\in  \prod_{i\in I}\GL_2(\F_q)\Big{|}\enspace\exists
 \delta\in\mathcal{D},\enspace \forall i, \enspace  \det(M_i)=\delta\right\}.$$

\begin{lemma}\label{borel} Assume $p\!>\!2k_0$. Then,
 
{\rm(i)} for all $\gp$ dividing $p$, $\overline{\rho}( I_{\gp})$
is contained (possibly after conjugation by  an element of $\GL_2(\F_q)$), 
either in the Borel's   subgroup of   $\GL_2(\F_q)$, either in the 
 non-split torus  of $\GL_2(\F_q)$.
The second case cannot occur if  $f$ is ordinary at $p$. 
 
{\rm(ii)} $\Ind_F^{\Q}\overline{\rho}(I_p)\subset  \phi(H(\F_q))$.
\end{lemma}

{\bf Proof : }  (i) Put $h:=|J_{F,\gp}|$. By   
Cor.\ref{poids-modere}  we have 
$\overline{\rho}|_{I_{\gp}}^{\enspace\mathrm{s.s.}}=
\varepsilon_{\gp}\oplus \delta_{\gp}$, 
 where $\varepsilon_{\gp}$ (resp. $\delta_{\gp}$)
$I_{\gp}\rightarrow  \kappa^\times $ are obtained by
composing the tame inertia map $I_{\gp}\rightarrow
\F_{p^h}^\times$ and the character $\varepsilon: x\mapsto
\underset{\tau\in J_{F,\gp}}{\prod}  \tau(x)^{p_\tau}$ (resp. $\delta
:x\mapsto \underset{\tau\in J_{F,\gp}}{\prod}\tau(x)^{q_\tau}$),
where $\{p_\tau,q_\tau\}=\{m_\tau,k_0-1-m_\tau\}$.

Let $x_h$ be a generator of $\F_{p^h}^\times$. As the traces of the  elements 
of $\overline{\rho}(\G_{\widehat{F}})$ are  in $F_q\coprod y\F_q$,
we have  $(\varepsilon(x_h)+\delta(x_h))^2\in\F_q$ and therefore
$\varepsilon(x_h)^2+\delta(x_h)^2\in\F_q$.

If $\varepsilon(x_h)^2,\delta(x_h)^2\in\F_q^\times$, then 
it is easy to see that $\varepsilon(x_h),\delta(x_h)\in\F_q^\times$
(we use $p>k_0$ and $p\neq 2k_\tau-1$). Therefore $I_{\gp}$
fixes a  $\F_q$-rational line, and $\overline{\rho}( I_{\gp})$
 is contained  in a Borel subgroup of  $\GL_2(\F_q)$.

Otherwise $\varepsilon(x_h)^2$ and $\delta(x_h)^2$ are conjugated by the 
non-trivial element  of $\Gal(\F_{q^2}/\F_q)$, hence 
$\varepsilon(x_h)^2=\delta(x_h)^{2q}$. Since  $p\!>\!2k_0$, we have
$\varepsilon(x_h)=\delta(x_h)^q$ and so 
$\varepsilon(x_h)+\delta(x_h)^q\in \F_q^\times$. Hence 
$\tr(\overline{\rho}( I_{\gp}))\subset \F_q$, and therefore
$\overline{\rho}( I_{\gp})\!\subset\! \GL_2(\F_q)$.
In this case  $\overline{\rho}( I_{\gp})$ is contained in a 
non-split torus of $\GL_2(\F_q)$.
If  $f$ is ordinary at $p$, then for all $\tau$,  
$p_\tau=m_\tau<k_0-m_\tau-1=q_\tau$ and therefore
$\varepsilon(x_h)^2\neq \delta(x_h)^{2q}$, and 
so the second case cannot occur.

(ii) The determinant condition ${}^\mathcal{D}$ being satisfied,  all we 
have to check is the following : for all $i\in I$ and $\tau,\tau'\in J_F^i$
the character 
$$ I_p \rightarrow \{\pm 1\} \text{ , } g\mapsto (\overline{\rho}(\widetilde{\sigma}_{i,\tau}^{-1}g
\widetilde{\sigma}_{i,\tau}))^{-1} (\overline{\rho}(\widetilde{\sigma}_{i,\tau'}^{-1}g
\widetilde{\sigma}_{i,\tau'}))$$
is trivial. This follows, as in the proof of the 
Prop.\ref{k-not-IND}, from the assumption $p\!>\!2k_0$. \hfill $\square$


\begin{lemma}\label{big-image}
$\phi(H(\F_q)) \subset \Ind_F^{\Q}\overline{\rho}(\G_{\widehat{F}}) $,
\end{lemma}

{\bf Proof : } We have seen in the beginning of this paragraph that
$\pr(\phi(\SL_2(\F_q)^I))\subset 
\pr(\Ind_F^{\Q}\overline{\rho}(\G_{\widehat{F}}))$. 
By lemma \ref{sl2}, we deduce that 
$\phi(\PSL_2(\F_q)^I)\subset \Ind_F^{\Q}\overline{\rho}(\G_{\widehat{F}})$.

As $\phi(H(\F_q))=\phi(\SL_2(\F_q)^I) \Ind_F^{\Q}\overline{\rho}(I_p)$, we're done.
\hfill $\square$

\begin{prop}\label{LI-Ind}
Assume that $f$ is not a theta series and that 
$\mathrm{{\bf (LI}_{Ind \overline{\rho}}{\bf )}}$  does not hold
for infinitely many primes $p$. Then, there exists 
$\tau\in J_F$, $ \tau\neq \mathrm{id}$ and a finite order Hecke character
$\varepsilon$ of $\widetilde{F}$ of conductor dividing 
$\Nm(\gn)$, such that for all prime 
$v\nmid \Nm(\gn)$ which splits completely 
in $\widetilde{F}$, we have $c(f_{\tau},v)=\varepsilon(v)c(f,v)$.
\end{prop}

{\bf Proof : } Since $f$ is not a theta series, we know that 
$\mathrm{{\bf (LI}_{\overline{\rho}}{\bf )}}$ hold for all
but finitely many primes $\p$. Take such a $p$, and 
assume that $\mathrm{{\bf (LI}_{Ind \overline{\rho}}{\bf )}}$ does not hold.
 Then there exist $\widetilde{\tau}_1,\widetilde{\tau}_2 \in \G_{\Q}$
such that $\tau:=\widetilde{\tau}_2^{-1}\widetilde{\tau}_1|_F\neq\mathrm{id} $
and for all $g\in\G_{\widehat{F}}$, 
$\pr\overline{\rho}(\widetilde{\tau}_2^{-1}g\widetilde{\tau}_1)=
\pr\overline{\rho}(\widetilde{\tau}_2^{-1}g\widetilde{\tau}_2)$.
As $\mathrm{{\bf (LI}_{\overline{\rho}}{\bf )}}$ hold and 
$\G_{\widehat{F}}$ is a normal subgroup of $\G_{\widetilde{F}}$, 
the above relation hold for every $g\in \G_{\widetilde{F}}$. 
Therefore, there exist a character $\varepsilon_{\gal}:\G_{\widetilde{F}}
\rightarrow \kappa^\times$, such that for all $g\in \G_{\widetilde{F}}$, 
$\overline{\rho}_{f_\tau}(g)=\varepsilon_{\gal}(g)\overline{\rho}_{f}(g)$.
Assume that $p>2k_0$. Then, 
the same argument as in the proof of Prop.\ref{k-not-IND} shows that 
$\varepsilon_{\gal}$ is unramified at primes dividing $p$.
By lemma \ref{lifting} $\varepsilon_{\gal}$ can then be 
lifted to a  finite order  Hecke character $\varepsilon$ of  $\widetilde{F}$ 
of conductor dividing $\Nm(\gn)$.  Because of the 
determinant relation $\overline{\psi_{\tau}}=\varepsilon_{\gal}^2
\overline{\psi}$, there finitely many such $\varepsilon$'s.

For every prime $v\nmid \gn p$ which splits completely  in 
$\widetilde{F}$, we have $c(f_{\tau},v)\equiv\varepsilon(v)c(f,v)\pmod{\cP}$.
If $\mathrm{{\bf (LI}_{Ind \overline{\rho}}{\bf )}}$ fails for
infinitely many $\cP$'s, then the congruence above become an 
equality. \hfill $\square$

\begin{cor} Assume that $F$ is a Galois field of odd degree and the central 
character of $\psi$ of $f$ is trivial ($F=\widehat{F}$). Assume moreover that  
$f$ is not a theta series and that
$\mathrm{{\bf (LI}_{Ind \overline{\rho}}{\bf )}}$   does not hold
for infinitely many primes $p$. Then, there exist 
a subfield $F' \subsetneq F$ and a Hilbert modular form 
$f'$ on $F'$, such that the base change of $f'$ to $F$ is  a twist of $f$ 
by a quadratic character of conductor dividing $\Nm(\gn)$.
\end{cor}

{\bf Proof : } As in the proof of Prop.\ref{LI-Ind} there exist a  quadratic 
character $\varepsilon$ of $F$ of conductor dividing $\Nm(\gn)$ and 
$\mathrm{id}\neq \tau\in \Gal(F/\Q)$ such that  we 
have $\rho_{f_\tau}= \varepsilon_{\gal}\otimes\rho$. Let $F'\subset F$ (resp.
$F_i\supset F$)  be the  fixed field of $\tau$ (resp. of 
$\ker(\varepsilon_{\tau^i})$). By assumption $F/F'$ is a cyclic extension 
of {\it odd} degree $h$. Let $F''=\prod_{i=1}^{h} F_i$. Then we have
$$\Gal(F''/F')=\{(u_1,..,u_h)\in\{\pm 1\}^h \enspace | \enspace\prod_{i=1}^{h} 
u_i=1\}\rtimes \{\tau^i\enspace | \enspace 0\leq i\leq h-1\} ,$$
where $\tau$ acts on $(u_1,..,u_h)$ by cyclic permutation. When $h=3$ the group 
$\Gal(F''/F')$ is isomorphic to $A_4$. 

The representation $\rho|_{\G_{F''}}$ is invariant by $\Gal(F''/F')$, 
but Langlangs' Cyclic Descend does not apply directly, because 
the order of $\Gal(F''/F')$ is even. Consider the  quadratic character 
$\delta=\varepsilon\cdot \varepsilon_{\tau^2}\cdot ..\cdot 
\varepsilon_{\tau^{h-1}}$. Then  the $\G_F$-representation  $\delta_{\gal}\rho$ 
is invariant by $\Gal(F/F')$,  so extends  to a representation of $\G_{F'}$.
By applying Langlangs' Cyclic Descend to $\delta\otimes f$
 we obtain $f'$ as desired. \hfill $\square$

\section{Boundary cohomology and congruence criterion.}

We recall that  $f\in S_k(\gn,\psi)$ is a  Hilbert modular newform.

\begin{defin}
We say that a  normalized eigenform $g\in  S_k(\gn,\psi)$ is 
 {\it  congruent } to  $f$ modulo   $\cP$, if  
their respective eigenvalues  for the Hecke operators   
(that is  their Fourier coefficients)  are congruent
 modulo $\cP$.

We say that a   prime $\cP$ is a 
{\it congruence prime} for $f$, if there exists 
a  normalized eigenform $g\in  S_k(\gn,\psi)$ 
 distinct from $f$ and congruent to  $f$ modulo  $\cP$.
\end{defin}

One expects that, as in the elliptic modular case 
(carried out by Hida  \cite{hida63,hida64}
and Ribet \cite{ribet4}),  the congruence  primes
for $f$ are controlled by the special value at $1$ of the 
adjoint $L$-function of $f$. Such  results  have been obtained 
 by Ghate \cite{ghate1} when $k$ is parallel. 

 Following \cite{hida63}, \cite{ghate1} and  using a
vanishing result of  the boundary cohomology,
we obtain a new result in this direction (see Thm.\ref{theo-A} and Thm.\ref{theo-B}(ii)).

\subsection{Vanishing of certain local  components  of 
the boundary  cohomology.}

We  introduce the following condition :

\medskip
\noindent {\bf (MW)}  the middle weight 
 $\frac{|p(J_F)|+|p(\varnothing)|}{2}=\frac{d(k_0-1)}{2}$ 
does not belong to   $\{|p(J)|,J\!\subset\! J_F\}$.

\medskip
This condition is automatically satisfied when  the  motivic  weight
 $d(k_0-1)$ is odd, or when $d=2$ and $k$ is non-parallel.

\begin{lemma} \label{key2}
Let $\rho_0$ be  a  representation of $\G_{\widetilde{F}}$ on a 
finite dimensional $\kappa$-vector space $W$.  Assume that 
for every $g\in \G_{\widetilde{F}}$,  the characteristic polynomial  of 
$(\otimes\Ind_F^{\Q}\overline{\rho})(g)$ annihilates $\rho_0(g)$.

{\rm (i)} If {\bf (I)}, {\bf (II)} and $\mathrm{\bf (LI_{\overline{\rho}})}$
 hold,  then for all $h\in\Z$, the weights  $h$ and $d(k_0-1)-h$ occur
with the  same multiplicity in 
 each  $\G_{\widetilde{F}}$-irreducible  subquotient  of $\rho_0$.

{\rm (ii)} If   {\bf (I)}, $\mathrm{\bf (Irr_{\overline{\rho}})}$ and {\bf (MW)} hold, and  
$p-1>\max(1,\frac{5}{d})\sum_{\tau\in J_F}(k_\tau-1)$, then
 each $\G_{\widetilde{F}}$-irreducible  subquotient  of $\rho_0$ 
 contains at least  two different  weights  for the action of 
the tame inertia at $p$. 
\end{lemma}

{\bf Proof : } We may assume  that $\rho_0$ is  irreducible. 

(i) By the lemmas \ref{borel}(ii) and \ref{big-image} we have 
$\Ind_F^{\Q}\overline{\rho}(I_p)\subset  \phi(H(\F_q))\subset 
\Ind_F^{\Q}\overline{\rho}(\G_{\widehat{F}})$. Let 
$T'$ be the torus of $H(\F_q)$ containing the image of the 
tame inertia, and $N'$ be the normalizer of $T'$ in $H(\F_q)$.
The image by $\phi$ of  $N'/T'\cong \{\pm 1\}^I$ is
the subgroup of  the Weyl group  $N/T=\{\pm 1\}^{J_F}$ of $G$
containing the  elements  which are constant on the partition 
$J_F=\coprod_{i\in I}J_F^i$.  In particular, the longest Weyl 
element $\epsilon_{J_F}$ belongs to the  image of $ N'/T'$.

Let  $x\in W$ be an eigenvector for the action of $T'$. 
By the annihilation condition, there exists a subset $J_x\subset J_F$, 
such that $I_p$ acts on $x$ by the weight $|p(J_x)|$.

Let $g_{J_F}\in \G_{\widetilde{F}}$
be such that $\Ind_F^{\Q}\overline{\rho}(g_{J_F})=\epsilon_{J_F}\mod{T'}$. 
Then $\rho_0(g_{J_F})(x)$ is of weight $|p(J_x\Delta J_F)|=d(k_0-1)-|p(J_x)|$.
Therefore, for each $h\in\Z$, $\rho_0(g_{J_F})$ gives  a  bijection 
between the eigenspaces for the tame inertia of weight $h$ and $d(k_0-1)-h$.

(ii) If $\mathrm {\bf (LI_{\overline{\rho}})}$ hold, then  the statement 
follows from  (i) and {\bf (MW)}. Otherwise, by Prop.\ref{LI} the image 
 $\pr(\overline{\rho}(\G_F))$ is dihedral.
Since $\widetilde{F}$ is totally real, $\pr(\overline{\rho}(\G_{\widetilde{F}}))$
is also dihedral  (see {\S}\ref{dihedral}). 

Denote by $N$ the normalizer of the standard torus $T$ in $G$. Put 
$N'=\Ind_F^{\Q}\overline{\rho}(\G_{\widetilde{F}}) \subset N(\kappa)$
and $T'=N'\cap T(\kappa)$. Then $N'/T'$ is a subgroup of the 
Weyl group $\{\pm 1\}^{J_F}=N/T$ of $G$. 

As we have seen in {\S}\ref{dihedral}, the representation 
 $\Ind_F^{\Q}\overline{\rho}$ is tamely ramified at $p$ and the image 
of the inertia group $I_p$ is contained in $T'$. 

Let  $x\in W$ be an eigenvector for the action of $T'$. 
By the annihilation condition, there exists a subset $J_x\subset J_F$, 
such that $I_p$ acts on $x$ by the weight $|p(J_x)|$.
For every element $\epsilon_J\in N'/T'$, ($J\subset J_F$), 
let $g_J\in \G_{\widetilde{F}}$
be such that $\Ind_F^{\Q}\overline{\rho}(g_J)=\epsilon_J\mod{T'}$. 
Then $\rho_0(g_J)(x)$ is of weight $|p(J_x\Delta J)|$. It remains to 
show that the $|p(J_x\Delta J)|$ are not all equal when 
$\epsilon_J$ runs over the elements of $N'/T'$. Note that,  for all 
$\tau\in J_F$, the $\tau$-projection  $N'/T'\rightarrow \{\pm 1\}$ 
is a surjective homomorphism (because the group 
$\pr(\overline{\rho}_{f_{\tau}}(\G_{\widetilde{F}}))$ is also dihedral). 
Therefore,  we have : 
$$\sum_{\epsilon_J\in N'/T'} |p(J_x\Delta J)|= |N'/T'| \frac{d(k_0-1)}{2}.$$
The statement now follows from the {\bf (MW)} assumption.
\hfill $\square$

\begin{rque}\label{key-diamond} 
The (i) in the previous lemma is a generalization, from the quadratic to 
the arbitrary degree case,   of the key lemma in  \cite{diamond2}.
This  lemma is  false in general under the only assumptions 
{\bf (I)}, {\bf (II)} and $\mathrm{\bf (Irr_{\overline{\rho}})}$, when the degree is $\geq 3$. In fact,
consider the following construction in the cubic case : let $L$ be a  
Galois extension  of $\Q$ of group $A_4$, such that the cubic subfield 
$F$ fixed by the Klein group is totally real; let $K$ be a quadratic extension 
of $F$ in $L$, and consider a theta series $f$ of weight $(2,2,2)$
attached to a Hecke character of $K$; then the tensor
induced representation $\otimes\Ind_F^{\Q}\rho$ has two
irreducible four dimensional subquotients of Hodge-Tate weights 
$(0,2,2,2)$ and $(1,1,1,3)$. 
\end{rque}

Let $\T'\subset \T$ be the subalgebra generated by the Hecke operators
outside a finite set of places containing those dividing $\gn p$.

\begin{theo} \label{boundary}
Assume that {\bf (I)}, $\mathrm{\bf (Irr_{\overline{\rho}})}$ and {\bf (MW)} hold,
and $p-1\!>\!\max(1,\frac{5}{d})\underset{\tau\in J_F}{\sum}(k_\tau\!-\!1)$. 
Denote by $\gm$ the  maximal ideal   of $\T$ corresponding to 
$f$ and  $p$ and put $\gm'=\gm\cap \T'$. Then  

{\rm (i)} the $\gm'$-torsion of the  boundary  cohomology
$\rH^\bullet_{\partial}(Y,\bV_n)(\kappa)[\gm']$ vanishes,

{\rm (ii)}  the Poincar{\'e}   pairing     
$\rH_{!}^d(Y,\bV_n(\cO))'_{\gm'} \times
\rH_{!}^d(Y,\bV_n(\cO))'_{\gm'}\rightarrow \cO$
is  a perfect duality of free $\cO$-modules of finite rank,

$\mathrm{(iii)}$ 
$\rH^\bullet(Y,\bV_n(\cO))_{\gm'}=
\rH_c^\bullet(Y,\bV_n(\cO))_{\gm'}=
\rH_{!}^\bullet(Y,\bV_n(\cO))_{\gm'}$.
\end{theo}

{\bf Proof : } (i) Consider the minimal compactification 
$Y_{\overline{\Q}}\overset{j}{\hookrightarrow} 
Y_{\overline{\Q}}^*\overset{i}{\hookleftarrow}  \partial Y_{\overline{\Q}}^*$.
The Hecke correspondences extend to  $Y_{\overline{\Q}}^*$. 
By  the  Betti-{\'e}tale comparison isomorphism, we  identify  (in a  
Hecke-equivariant way) the following two  long exact cohomology sequences :

\xymatrix{ ...\ar[r] 
& \rH_c^r(Y,\bV_n(\kappa)) \ar[r]\ar@{=}[d]
& \rH^r(Y,\bV_n(\kappa)) \ar[r]\ar@{=}[d]
& \rH^r_{\partial}(Y,\bV_n(\kappa)) \ar[r]\ar@{=}[d] & ...\\ ...\ar[r] 
& \rH^r(Y_{\overline{\Q}}^*,j_{!}\bV_n(\kappa)) \ar[r]
& \rH^r(Y_{\overline{\Q}}^*,j_{*}\bV_n(\kappa)) \ar[r]
& \rH^r(\partial Y_{\overline{\Q}}^*, i^*Rj_{*}\bV_n(\kappa)) \ar[r]& ...}

Consider the $\G_{\Q}$-module $W^r_{\partial}=
\rH^r(\partial Y_{\overline{\Q}}^*, i^*Rj_{*}\bV_n(\kappa))$. 
We have to show that $W^r_{\partial}[\gm']=0$. 

By the Cebotarev Density Theorem and the congruence relations at 
totally split primes of $F$, we can apply lemma \ref{key2} to
$W^r_{\partial}[\gm']$.  Therefore  each $\G_{\widehat{F}}$-irreducible
subquotient  of $W^r_{\partial}[\gm']$  has  at least  two different
weights  for the action of the tame inertia at $p$. So it is enough 
to show that  each $\G_{\Q}$-irreducible subquotient  of  $W^r_{\partial}$ 
is pure (=contains a single weight for the action of the tame inertia at $p$).
Because $\partial M_{\overline{\Q}}^*$ is   zero dimensional, 
the spectral sequence  $\rH^\bullet(\partial Y_{\overline{\Q}}^*,
i^*R^\bullet j_{*}\bV_n(\kappa))\Rightarrow 
\rH^\bullet(\partial Y_{\overline{\Q}}^*,
i^*Rj_{*}\bV_n(\kappa))$ shows that   $W^r_{\partial}=
\rH^0(\partial Y_{\overline{\Q}}^*,i^*R^r j_{*}\bV_n(\kappa))$.

As  $\rH^0(\partial Y_{\overline{\Q}}^*,i^*R^r j_{*}\bV_n(\kappa))$
is a subquotient of 
$\rH^0(\partial Y_{\overline{\Q}}^{1,*},i^*R^r j_{*}\bV_n(\kappa))$
it is enough  to show that  each $\G_{\Q}$-irreducible subquotient  of 
this last is pure.

This will be done using a result of Pink\cite{Pi}. 
We had  to replace $Y$ by $Y^1$, because the group $G$ does not
satisfy the conditions of this reference, while $G^*$ satisfies them.

Consider the decomposition   $T=D_l\times D_h$, 
according to   $\begin{pmatrix} u\epsilon & 0\\ 0 & u^{-1}\end{pmatrix}=
\begin{pmatrix} u& 0\\ 0 & u^{-1}\end{pmatrix}
\begin{pmatrix} \epsilon  & 0\\ 0 & 1\end{pmatrix}$. 
Put $\Gamma^1=\Gamma_1^1(\gc,\gn)$. 
By \cite{Pi}Thm.5.3.1, the restriction  of the {\'e}tale sheaf 
$i^*R^r j_{*}\bV_n(\F_p)$ to a cusp  $\mathcal{C}=\gamma\infty$ of 
$Y_{\overline{\Q}}^{1,*}$, is the  image by the functor of  Pink of the 
$\gamma^{-1}\Gamma^1\gamma\cap B /\gamma^{-1}\Gamma^1\gamma\cap D_lU$-module
$$\underset{a+b=r}\oplus \rH^a(\gamma^{-1}\Gamma^1\gamma\cap D_l ,
\rH^b(\gamma^{-1}\Gamma^1\gamma\cap U,\bV_n(\F_p))).$$

Under the assumption {\bf (II)}, a modulo $p$ version of a theorem
of  Kostant   (see \cite{PoTi}) gives an  isomorphism of $T$-module 
$\rH^b(\gamma^{-1}\Gamma^1\gamma\cap U,\bV_n(\F_p))=\underset{|J|=b}{\oplus}
W_{\epsilon_J(n+t)-t}$. By decomposing $W_{\epsilon_J(n+t)-t}=
W_{\epsilon_J(n+t)-t,l}\otimes W_{\epsilon_J(n+t)-t,h}$ 
according to $T=D_l\times D_h$, we get 
$$ \rH^a(\gamma^{-1}\Gamma^1\gamma\cap D_l,
\rH^b(\gamma^{-1}\Gamma^1\gamma\cap U ,\bV_n(\F_p)))=
\underset{|J|=b}{\oplus} \rH^a(\gamma^{-1}\Gamma^1\gamma\cap D_l,
W_{\epsilon_J(n+t)-t,l}) \otimes W_{\epsilon_J(n+t)-t,h},$$ 
where  Galois acts only on the  second factors of the right hand side.

Therefore $\rH^0(\partial M_{\overline{\Q}}^*,i^*R^r 
j_{*}\bV_n(\F_p))$ is a direct sum of subspaces 
$\rH^0(\mathcal{C},W_{\epsilon_J(n+t)-t,h}(\F_p))$, $|J|\leq r$,
each containing a single Fontaine-Laffaille weight, namely 
the weight $|p(J)|$.

(ii)  As the Poincar{\'e} duality   is perfect over $E$, it is enough 
to show that the $\gm'$-localization of natural map 
$\rH^d(\cO)/\rH_{!}^d(\cO)\rightarrow 
\rH^d(E)/\rH_{!}^d(E)$ is  injective. For this, it is sufficient to
show that   $\rH_{\partial}^d(\cO)_{\gm'}:=
 \rH^d({\partial}M,\bV_n(\cO))_{\gm'}$ is  torsion free,
which is a consequence of the vanishing of 
 $\rH_{\partial}^{d-1}(E/\!\cO)_{\gm'}$. 
By (i) and Nakayama'a lemma  $\rH_{\partial}^{d-1}(\kappa)_{\gm'}=0$ and 
moreover  we have a surjection $\rH_{\partial}^{d-1}(\kappa)_{\gm'}\twoheadrightarrow 
\rH_{\partial}^{d-1}(E/\!\cO)_{\gm'}[\varpi]$, 
where $\varpi$ is an uniformizer of $\cO$.  

(iii) The vanishing of $\rH_{\partial}^\bullet(\kappa)_{\gm'}$
gives the vanishing of $\rH_{\partial}^\bullet(\cO)_{\gm'}=0$.
\hfill $\square$

\subsection{Definition of periods.}\label{periodes}
 By taking the subspace  $\underset{\ga\subset 
\go}{\bigcap} \ker(T_{\ga}-c(f,\ga))$ 
of (\ref{ESH-epsilon})  we obtain  

$$\delta_J : \C f_J \overset{\sim}\longrightarrow \rH_{!}^d
(Y^{\an},\bV_n(\C))[\widehat{\epsilon}_J,f].$$

Fix an isomorphism $\C\cong\overline{\Q}_p$
compatible with $\iota_p$. We recall that   $\rH_{!}^d(Y^{\an},\bV_n(\cO))'$  
denotes the image  of the natural map  $\rH_c^d(Y^{\an},\bV_n(\cO))
\rightarrow \rH^d (Y^{\an},\bV_n(\C))$. 
As $\cO$ is principal, the  $\cO$-module   
$L_{f,J}:=\rH_{!}^d(Y^{\an},\bV_n(\cO))'[\widehat{\epsilon}_J,f]$ 
is free (of rank $1$). We fix a basis $\eta(f,J)$ of $L_{f,J}$. 

\begin{defin}
For each $J\!\subset\!  J_F$ we define the period
$\Omega(f,J)=\frac{\delta_J(f_J)}{\eta(f,J)}\in \C^\times\!/\!\cO^\times.$

\noindent We fix  $J_0\!\subset\!  J_F$ and
put $\Omega_f^+=\Omega(f,J_0)$ and  $\Omega_f^-= \Omega(f,J_F\!\bs\!J_0)$.
\end{defin}

\begin{rque}
The  periods $\Omega_f^{\pm}$  differ from  
the ones originally  introduces by Hida in \cite{hida63}.
The  Hida's  periods put together all the external conjugates of $f$.  
Our  slightly  different definition is  motivated  by the  
congruence criterion  that we want to show  (Thm.\ref{theo-A}). 
As we can prove the perfectness of the twisted Poincar{\'e} pairing 
only for certain local components of the middle degree cohomology
$\rH_{!}^d(Y^{\an},\bV_n(\cO))'$, and that 
in general $f$ and his external conjugates do not belong 
to the same local component,  we have to separate them 
in  the period's definition.
\end{rque}

\subsection{Computation of a  discriminant.}
The aim of this  paragraph is the computation of the discriminant 
$\disc(L_f)$ of the  $\cO$-lattice $L_f:=\rH_{!}^d
(Y^{\an},\bV_n(\cO))'[f] =
\oplus_{J\subset J_F} L_{f,J}$, with respect twisted Poincar{\'e}  pairing
 $[\enspace, \enspace]$. We follow \cite{ghate1} Sect.6. 

We have $\disc(L_f)=\det(([\eta(f,J),\eta(f,J')])_{J,J'\subset J_F})$. 

By \cite{ghate1} (41),  for every $\tau\in J_F$ and $x,y \in \rH_{!}^d
(Y^{\an},\bV_n(\C))$ we have  $[\epsilon_\tau\cdot x,y]=-[x,\epsilon_\tau\cdot y]$.

The  embedding  $\cO\hookrightarrow \C$ that we have  fixed  
gives   an embedding $\tau_0 :F \hookrightarrow \C$. We have 
$$\hspace{-4mm}\disc(L_f)=\!\prod_{\tau_0\in J\!\subset\! J_F}
\begin{vmatrix}  0  & \left[\eta(f,J),\eta(f,J_F\!\bs\!J)\right] \\ 
\left[\eta(f,J_F\!\bs\!J),\eta(f,J)\right] & 0 \end{vmatrix}=\!
\prod_{\tau_0\in J\!\subset\! J_F}-\left(
\frac{[\delta_J(f),\delta_{J_F\!\bs\!J}(f)]}{\Omega(f,J)\Omega(f,J_F\!\bs\!J)}
\right)^2 $$

We have  $[\delta_J(f),\delta_{J_F\!\bs\!J}(f)]\!=\!
2^d \langle \epsilon_{J_F}\delta(f),\iota\cdot\delta(f)\rangle\!=\!
2^d W(f)\langle \epsilon_{J_F}\delta(f),\delta(f^c)\rangle 
\!=\!2^d W(f)(f,f)_{\gn},$

\noindent  where $f^c$ denotes is the complex conjugate 
 of $f$, $\iota$ is the Atkin-Lehner involution  and $W(f)$ is  the 
complex constant  of the functional  equation  of the standard $L$-function 
 of $f$.  Thus we get the following equality in $E^\times\!/\cO^\times$  :
\begin{equation}\label{discriminant}
\disc(L_{f,J_0}\oplus L_{f,J_F\!\bs\!J_0})=
\left(\frac{W(f)(f,f)_{\gn}}{\Omega_f^+\Omega_f^-}\right)^2.
\end{equation}

\subsection{Shimura's  formula for $L(\Ad^0(f),1)$.} \label{L-ad}
For each  prime ideal $v$ of $F$  we define  
$\alpha_v$ and $\beta_v$ by the equations 
$$\alpha_v+\beta_v=c(f,v),\enspace \enspace\alpha_v\beta_v=\begin{cases}
\psi(v)\Nm(v) & \text{, if }v\nmid \gn, \\
0  & \text{, if }v\mid \gn. \end{cases}$$

The naive adjoint $L$-function of $f$ is then given by the Euler product :
\begin{equation}\label{Lnaive}
L^0(\Ad^0(f),s)=\prod_{v\nmid \gn}\left[
(1\!-\!\alpha_v\beta_v^{-1}\Nm(v)^{-s}) (1\!-\!\Nm(v)^{-s})
(1\!-\!\beta_v\alpha_v^{-1}\Nm(v)^{-s}) \right]^{-1}.
\end{equation}

We denote by $f^c$ the (external) complex conjugate of $f$. 
We  introduce  a twisted version  of the $L$-function  
associated to  the tensor product  $f\otimes f^c$  (see \cite{shimura45}) :
$$D(f,f^c,s)=\prod_{v}\left(1\!-\!
\alpha_v\beta_v\overline{\alpha_v\beta_v}\Nm(v)^{-2s}\right)
\left(1\!-\!{\alpha_v}\overline{\alpha_v}\Nm(v)^{-s}\right)^{-1}\cdot$$
$$\cdot\left(1\!-\!{\alpha_v}\overline{\beta_v}\Nm(v)^{-s}\right)^{-1}
\left(1\!-\!{\beta_v}\overline{\alpha_v}\Nm(v)^{-s}\right)^{-1}
\left(1\!-\!{\beta_v}\overline{\beta_v}\Nm(v)^{-s}\right)^{-1}$$

\noindent  as well as a  naive version $D^0(f,f^c,s)$ obtained by 
removing the  factors  for $v|\gn$.

Using that  for all ${v\nmid  \gn}$, 
  $\overline{c(f,v)}= \overline{\psi(v)} c(f,v)$, 
a direct computation, shows that 
\begin{equation}\label{Dnaive}
 \zeta_F^0(2s)D^0(f,f^c,s\!+\!k_0\!-\!1)=
\zeta_F^0(s) L^0(\Ad^0(f),s).
\end{equation}

The advantage to switch to $D(f,f^c,s)$ is that we  have the following 
formula, proved by Shimura  (see \cite{HiTi} lemma 7.2, for a proof that 
the $D(f,f^c,s)$ defined above equals the one  studied by Shimura) :

\begin{theo}\rm{(Shimura \cite{shimura45} (2.31), Prop. 4.13)}
Let $f\in S_k(\gn,\psi)$ be a newform. Then 
$$\Res_{s=1}D(f,f^c,s+k_0\!-\!1)=2^{d-1}(4\pi)^{|k|}\prod_{\tau\in
  J_F}\Gamma(k_\tau)^{-1}R_F[\go_{+}^\times:
\go^{\times 2}]\langle f,f\rangle \text{, where} $$
  $\langle f,f\rangle= \mu(\Gamma\bs \mathfrak{H}_F)^{-1}
(f,f)_{\gn}$, and  $\mu(\Gamma\bs \mathfrak{H}_F)=
\frac{2\Nm(\gd)^{3/2}\zeta_F(2)\Nm(\gn)}
{\pi^d[\go_{+}^\times:\go^{\times 2}]\prod_{v\mid\gn} (1+\Nm(v))^{-1}}.$
\end{theo}

We deduce the formula :
\begin{equation}\label{residu}
\zeta_F^0(2)\Res_{s=1}D(f,f^c,s+k_0\!-\!1)=
\frac{(4\pi)^{|k|}\pi^d\Res_{s=1}\zeta_F^0(s)}{2\Delta h_F\prod_{\tau\in
  J_F}\Gamma(k_\tau) } (f,f)_{\gn}.
\end{equation}

We define the imprimitive adjoint $L$-function    $L^*(\Ad^0(f),s)$
by  completing the naive adjoint $L$-function   $L^0(\Ad^0(f),s)$, defined 
in (\ref{Lnaive}), in order to have  the relation 
$$L^*(\Ad^0(f),s)D^0(f,f^c,s+k_0\!-\!1)=L^0(\Ad^0(f),s)D(f,f^c,s+k_0\!-\!1).$$

An    explicit computation  of  \cite{HiTi} (7.7) gives  
 $L^*(\Ad^0(f),s)=L^0(\Ad^0(f),s)
\prod_{v | \gn} L^*_v(\Ad^0(f),s) $, where for $v | \gn$ 
$$L^*_v(\Ad^0(f),s)=\begin{cases}1\!-\! \Nm(v)^{-s} & \text{, if }f 
\text{ is principal series  and minimal at  } v, \\ 
1\!-\! \Nm(v)^{-s-1} & \text{, if }f 
\text{ is special and minimal at  } v,
\\ 1 & \text{, otherwise. }\end{cases} $$

Following Deligne \cite{deligne33} we  associate to $L^*(\Ad^0(f),s)$
the  Euler factor at infinity    
$$\Gamma(\Ad^0(f),s)=\prod_{\tau\in J_F}
\pi^{-(s+1)/2}\Gamma((s+1)/2)
(2\pi)^{s+k_\tau-1}\Gamma(s+k_\tau-1).$$
Finally from (\ref{Dnaive}) and  (\ref{residu}) we obtain  
\begin{equation}\label{adjoint}
\Lambda^*(\Ad^0(f),1):=\Gamma(\Ad^0(f),s)L^*(\Ad^0(f),s)=
\frac{2^{|k|-1}}{\Delta h_F}(f,f)_{\gn}.
\end{equation}

\begin{rque}
Consider the adjoint $L$-function $L(\Ad^0(\rho),s)$ of the 
three dimensional $\G_F$-representation $\Ad^0(\rho)$ (on
trace zero matrices). By the compatibility between local and 
global Langlands correspondence $L(\Ad^0(\rho),s)$ equals 
the adjoint $L$-function  $L(\Ad^0(f),s)$ associated to the 
automorphic representation attached to $f$. 
Nevertheless $L(\Ad^0(f),s)$ may differ from  
$L^*(\Ad^0(f),s)$ at some  places $v$ dividing $\gn$ (see
\cite{HiTi} (7.3c)).
\end{rque}

\subsection{Construction of  congruences.}

\begin{lemma}\label{C0}
Let   $V_1$ and $V_2$ be two finite  dimensional $E$-vector spaces
and let $L$ be a $\cO$-lattice in  $V=V_1\oplus V_2$. 
For $j=1,2$, put $L_j=L\cap V_j$ and denote $L^j$  the projection of 
 $L$ in $V_j$ following the above direct sum decomposition. Then :

{\rm (i)} $L_j\!\subset\! L^j $ are two lattices of $V_j$, and
 $L_j$ is a direct factor  in $L$.

{\rm (ii)} we have isomorphisms of finite $\cO$-modules :
\begin{equation}\label{fusion}
L^1/L_1\overset{\sim}{\leftarrow} L/L_1\oplus
L_2\overset{\sim}{\rightarrow}L^2/L_2
\end{equation}
This finite $\cO$-module is called the congruence module, and 
is denoted by $C_0(L;V_1,V_2)$.
\end{lemma}

The following proposition follows from  Deligne-Serre 
(Lemma 6.11 of \cite{SeDe}) and will be used to construct congruences :

\begin{prop}\label{relevement}
Keep the notations of the  lemma \ref{C0}. 
Let $\cT$ be a commutative $\cO$-algebra consisting of 
endomorphisms of $V$, preserving  the lattice $L$ and the 
 direct sum decomposition  $V_1\oplus V_2$. Denote, for $j=1,2$, by 
  $\cT_j$   the image  of $\cT$ in $\End(V_j)$.

Assume that  $C_0(L;V_1,V_2)$ is non zero :
$\{\cP\}=\Ass(C_0(L;V_1,V_2))=\Supp(C_0(L;V_1,V_2))$.

Let $\gm_1$ be maximal ideal
$\cT_1$,  of residue  field  $\kappa_1$, such that  
$L^1/L_1\otimes_{\cT_1}\kappa_1$ is non zero, and denote by 
$\overline{\theta_1}:\cT_1\rightarrow \kappa_1$ the corresponding 
character.

Then  there exists a discrete valuation ring
$\cO'$ of maximal ideal   $\cP'$ (with
$\cP'\cap\cO=\cP$), of residue field  $\kappa'\supset \kappa_1$ 
and whose fraction field  $E'$ is a finite extension of $E$, 
and there exists  a 
character $\theta_2:\cT_2\rightarrow \cO'$ such that
for each $T\in\cT$,  $\enspace \overline{\theta_1}(T) \equiv
\theta_2(T) \pmod{\cP'}$.
\end{prop}

{\bf Proof : } Denote by   $\pi_j$ the projection of $\cT$  onto 
$\cT_j$, $j=1,2$. Then  $\gm=\pi_1^{-1}(\gm_1)$
is a maximal ideal of $\cT$ of residue field $\kappa_1$.
Put  $\gm_2=\pi_2(\gm)$.
As the isomorphism (\ref{fusion})  of the lemma $\ref{C0}$  is 
 $\cT$-equivariant, we get 
$$(L^1/L_1)\otimes_{\cT_1}(\cT_1/\gm_1)\cong
(L/(L_1\oplus L_2))\otimes_{\cT}(\cT/\gm) \cong
(L^2/L_2)\otimes_{\cT_2}(\cT_2/\gm_2)$$

By assumption $(L^1/L_1)\otimes_{\cT_1}(\cT_1/\gm_1)$
 is non zero. Therefore  $\gm_2$ is a 
maximal ideal of $\cT_2$ of residue field $\kappa_1$ and 
the corresponding  character  $\overline{\theta_2}:\cT_2\rightarrow \kappa_1$ 
fits in the  following commutative diagram :
$$\xymatrix@R=5pt{    &\cT_1\ar[rd]^{\overline{\theta_1}} & \\
\cT\ar[rd]\ar[ru] &   & \kappa_1    \\
  & \cT_2\ar[ru]_{\overline{\theta_2}}      & }$$

As $\cT_2$ is a (finite) flat $\cO$-algebra, there exist a
prime ideal $\cP_2$, contained in $\gm_2$ and such 
that $\cP_2\cap\cO=0$. The reduction modulo  $\cP_2$ gives 
a character $\theta_2$ of $\cT_2$  as in the statement. 

\begin{theo}{\rm(Theorem A)}\label{theo-A} 
Let  $f$ and  $p$ be such that {\bf (I)}, $\mathrm{\bf (Irr_{\overline{\rho}})}$ and {\bf (MW)} hold,
and $p-1>\max(1,\frac{5}{d})\sum_{\tau\in J_F}(k_\tau-1)$. If  
$\iota_p(\frac{W(f)\Lambda^*(\Ad^0(f),1)}{\Omega_f^+\Omega_f^-})\in\cP$, 
then  $\cP$  is a congruence prime for $f$.
\end{theo}

{\bf Proof : } We put $L=\rH_{!}^{d}(Y^{\an},\bV_n(\cO))_{\gm'}'
[\pm\widehat{\epsilon}_{J_0},\psi]\!\subset\! 
V=\rH_{!}^{d}(Y^{\an},\bV_n(E))_{\gm'} [\pm\widehat{\epsilon}_{J_0},\psi]$. 

We define   $V_1=\rH_{!}^{d}(Y^{\an},\bV_n(E))[\pm\widehat{\epsilon}_{J_0},f]$
(with the notations of {\S}\ref{periodes}
we have  $L_1=L\cap V_1=L_{f,J_0}\oplus L_{f,J_F\!\bs\!J_0}$). 

It follows from (\ref{discriminant}) that 
the bilinear  form  $[\enspace,\enspace ]$ is non-degenerate on $V_1$.
Let  $V_2$ be the orthogonal subspace of $V_1$ in $V$.

By  (\ref{adjoint}) and (\ref{discriminant}) and the assumption on $\cP$, we have
$0\neq \disc(L_1)\in \cP$. By  Thm.\ref{boundary}(ii) the $\cO$-lattice $L$ is 
autodual with respect to the twisted Poincar{\'e} pairing and  therefore, 
$\disc(L_1)=[L^1:L_1]$. Then by lemma \ref{C0}, the congruence module 
$C_0(L;V_1,V_2)$ is non-zero and $\cP$ belongs to its support. By  
Prop.\ref{relevement}  and the duality between   $\T(\C)$ and $S_k(\gn,\psi)$, we 
obtain another normalized eigenform $g\in S_k(\gn,\psi)$  congruent to $f$. 
Hence, $\cP$  is a congruence prime for $f$. \hfill $\square$

\section{Fontaine-Laffaille weights  of the Hilbert modular varieties.}

In this section all the objects are over $\cO$. 
The aim  is to establish a modulo $p$ version 
of Thm.\ref{hodge-pad}   under the   assumptions  that $p$ does not
divide $\Delta$ and $p-1>|n|+d$.

\subsection{The BGG complex over $\cO$.} \label{BGG} ${ }$

\medskip
\noindent{\bf Koszul's complex.}
The  Koszul's complex of the trivial $G$-module   $\cO$ is given by     
$$...\rightarrow  U_{\cO}(\mathfrak{g})\otimes \wedge_{\cO}^2\mathfrak{g}
\rightarrow  U_{\cO}(\mathfrak{g})\otimes\mathfrak{g}
\rightarrow  U_{\cO}(\mathfrak{g})
\rightarrow \cO
\rightarrow 0$$

As  $\mathfrak{g}=\gb\oplus\gu^-$, the 
$\cO[\gb]$-module
$\mathfrak{g}/\gb$ is a direct factor  in 
$\mathfrak{g}$ and  we have a homomorphism of $B$-modules 
$ U_{\cO}(\mathfrak{g})\otimes \wedge_{\cO}^{\bullet}\mathfrak{g}
\rightarrow 
U_{\cO}(\mathfrak{g})\otimes_{U_{\cO}(\gb)}\wedge_{\cO}^{\bullet}(\mathfrak{g}/\gb)$. 
Thus, we   deduce another  complex  
$$ U_{\cO}(\mathfrak{g})\otimes_{U_{\cO}(\gb)}\wedge_{\cO}^{\bullet}
(\mathfrak{g}/\gb) \rightarrow \cO \rightarrow 0,$$
denoted by $S_{\cO}^{\bullet}(\mathfrak{g},\gb)$.

More generally,  for a free  $\cO$-module  $V$ endowed with an 
 action  of $U_{\cO}(\mathfrak{g})$, we consider  the complex  
$S_{\cO}^{\bullet}(\mathfrak{g},\gb)\otimes V$, endowed 
de the diagonal action  of $U_{\cO}(\mathfrak{g})$. 

For every   $U_{\cO}(\gb)$-module  $W$, free over $\cO$,
 we have a canonical isomorphism  of $U_{\cO}(\mathfrak{g})$-modules 
\begin{equation}
\left(U_{\cO}(\mathfrak{g})\otimes_{U_{\cO}(\gb)} W\right)\otimes V
\cong U_{\cO}(\mathfrak{g})\otimes_{U_{\cO}(\gb)} 
\left( W\otimes V|_{\gb}\right),
\end{equation}

So we get another complex 
$$ U_{\cO}(\mathfrak{g})\otimes_{U_{\cO}(\gb)}
\left(\wedge_{\cO}^{\bullet}(\mathfrak{g}/\gb)\otimes  V|_{\gb}\right)
\rightarrow V \rightarrow 0,$$
denoted $S_{\cO}^{\bullet}(\mathfrak{g},\gb,V)$.
In the case where $V=V_n$ we  denote  it 
$S_{\cO}^{\bullet}(\mathfrak{g},\gb,n)$.

\medskip
\noindent{\bf Verma modules.}
For each weight $\mu\in \Z[J_F]$, we define a  $U_{\cO}(\mathfrak{g})$-module 
$V_{\cO}(\mu):=U_{\cO}(\mathfrak{g}) \otimes_{U_{\cO}(\gb)}W_{\mu}(\cO)$, called the 
 {\it Verma  module} of weight $\mu$.

\begin{lemma}
Let $W$ be a $B$-module, free of finite rank  over $\cO$, whose 
weights are smaller than $ (p-1)t$. Then, there exists a filtration  of 
$B$-modules $0=W_0\!\subset\! W_1\!\subset\! ...\!\subset\! W_r=W$ such that 
for each $1\leq i\leq r$, $W_i/W_{i+1}\cong W_{\mu_i}(\cO)$, for some
 $\mu_i\in \Z[J_F]$. Moreover the  $W_{\mu_i}(\cO)$, $1\leq i\leq r$,
are  the irreducible factors  of the $T$-module $W$.

In  particular, if $U$ acts trivially on $W$, then  
$W\cong \oplus_{i=1}^{r}  W_{\mu_i}(\cO)$.
\end{lemma}

{\bf Proof : } Let $\mu_1$ be a maximal weight  of $W$ (for the partial   order  
given by the   positive roots  of $G$) and let  $v\in W$ be a $\cO$-primitive 
vector  of weight  $\mu_1$. Let $W'$ be the $U_{\cO}(\gb)$-submodule
 generated by $v$.  Then  $W'\cong W_{\mu_1}(\cO)$ and 
$W'\otimes \kappa$ is irreducible, because  $\mu_1$ is  smaller than $ (p-1)t$ 
(and  $W'$ is free of rank 1).
As $W$ is free over $\cO$ we have an  exact
 sequence  of $B$-modules  $$0\rightarrow \Tor_1^{\cO}(W/W',\kappa)
\rightarrow W'\otimes \kappa \rightarrow W\otimes \kappa.$$ As $W'\otimes \kappa$ is 
irreducible and $v$ is primitive, the last arrow is injective. Therefore 
$$\Tor_1^{\cO}(W/W',\kappa)=0,$$ that is  $W/W'$ is free over $\cO$.  
The lemma follows then  by induction. 
\hfill $\square$

\begin{lemma}\label{verma} The module  
$S^{i}_{\cO}(\mathfrak{g},\gb,n)$ has  a finite filtration  by    
$U_{ \cO}(\mathfrak{g})$-submodules whose graded pieces are of
 the form $V_{\cO}(\mu)$, $\mu\in \Omega^i(n)$,
where $\Omega^i(n)$ denoted  the set  of  weights of the 
$\gt$-module  $\wedge_{\cO}^{i}(\mathfrak{g}/\gb)\otimes  
V_{n}(\cO)|_{\gb}$.
\end{lemma}

{\bf Proof : } Since $p-1>|n|+d$ the previous lemma applies to 
 $\wedge_{\cO}^{\bullet}(\mathfrak{g}/\gb)\otimes
 V_{n}(\cO)|_{\gb}$. This gives  a filtration  $0\!=\!W_0\!\subset\! W_1\!\subset\! ...\!\subset\! W_r=
\wedge_{\cO}^{i} (\mathfrak{g}/\gb)\otimes  
V_{n}(\cO)|_{\gb}$ whose graded pieces are  $W_{\mu}(\cO)$, $\mu\in \Omega^i(n)$.
 As $U_{\cO}(\mathfrak{g})$ is 
  $U_{\cO}(\gb)$-free, the  functor $U_{\cO}(\mathfrak{g})
\otimes_{U_{\cO}(\gb)}\bullet$  is  exact. \hfill $\square$

\medskip
\noindent{\bf Central characters.}
Let $U_{\cO}(\mathfrak{g}) \rightarrow U_{\cO}(\gt)$
be the  projection coming from the  Poincar{\'e}-Birkoff-Witt decomposition
$U_{\cO}(\mathfrak{g})=U_{\cO}(\gt)\oplus
(\gu^-U_{\cO}(\mathfrak{g})+U_{\cO}(\mathfrak{g})\gu$). 
We take its restriction  to the invariants for the adjoint action  
$\theta : U_{\cO}(\mathfrak{g})^G\rightarrow U_{\cO}(\gt)$.
Note  that $U_{\overline{\F}_p}(\gt)$ 
identifies itself  with the algebra of regular  functions  on
$\Hom_{\cO}(\gt,\overline{\F}_p)\cong \overline{\F}_p[J_F]$
(a Laurent polynomial algebra).   The Weyl group $\{\pm 1\}^{J_F}$ 
of $G$ acts on it by  $(\epsilon_J\cdot P)(\mu)=P(\epsilon_J(\mu+t)-t)$.
The following result  is a analogous to a  the theorem of
Harish-Chandra :

\begin{theo}{\rm(Jantzen \cite{jantzen})}
$\theta_{\overline{\F}_p}$ induces  an algebra isomorphism  
$U_{\overline{\F}_p}(\mathfrak{g})^G\rightarrow 
U_{\overline{\F}_p}(\gt)^{\{\pm 1\}^{J_F}}$.
\end{theo}

For every $\mu\in \Z[J_F]$ and every $\cO$-algebra $R$, we denote by 
$d\mu_R:\gt_R\rightarrow R$ the corresponding
character  and by $\chi_{\mu,R}=d\mu_R\circ \theta_R$ the 
composed map 
$U_{R}(\mathfrak{g})^G\rightarrow U_{R}(\gt)\rightarrow R$. 
This  definition is compatible with the $\cO$-algebra homomorphisms.

If  $V$ is a $U_{R}(\mathfrak{g})$-module generated by a vector $v$ 
of  weight $\mu$, that  is annihilated by $\gu$, then 
$U_{R}(\mathfrak{g})^G$ acts over $V$ by $\chi_{\mu,R}$.
Put $\chi_{\mu,p}=\chi_{\mu,\cO}$
and  $\overline{\chi}_{\mu,p}=\chi_{\mu,\overline{\F}_p}$.

\begin{cor}
Let $\mu\in \Z[J_F]$. If 
$\overline{\chi}_{n,p}=\overline{\chi}_{\mu,p}$,
then  there exists $J\!\subset\! J_F$ such that 
$\mu-(\epsilon_J(n+t)-t)\in p\Z[J_F]$. In particular, 
if   $\mu$ is  smaller than $ (p-1)t$, then  we have $\mu=\epsilon_J(n+t)-t$.
\end{cor}

\begin{prop}\label{chandra}
Let  $\mu\in \Omega^{i}(n)$ (see  lemma \ref{verma}). 
Then  $\overline{\chi}_{n,p}=\overline{\chi}_{\mu,p}$,
if and only if there exists a subset 
$J\!\subset\! J_F$containing $i$ elements  and such that $\mu=\epsilon_J(n+t)-t$.
\end{prop}

 {\bf Proof : }
By  the corollary, it remains to show that for   $J\!\subset\! J_F$,
we have   $\epsilon_J(n+t)-t\in \Omega^{i}(n)$,
if and only if  $|J|=i$. By  the lemma \ref{verma}, we have to show that 
 $W_{\epsilon_J(n+t)-t}(E)$ occurs in 
$\wedge_{E}^{i}(\mathfrak{g}/\gb)\otimes  
V_{n}(E)|_{\gt}$ (with multiplicity one)
 if and only if  $|J|=i$. The weight of 
$\wedge_{E}^{i}(\mathfrak{g}/\gb)\otimes  
V_{n}(E)|_{\gt}$  are of the form 
$\epsilon_{J'}(n+t)-t+\nu$, where $J'\!\subset\! J_F$ is a subset containing 
$i$ elements   and $\nu$ is a weight of  $V_{n}(E)$. 
Therefore  $\epsilon_J(n+t)-t=\epsilon_{J'}(n+t)-t+\nu$ and so 
$n=\epsilon_J(\nu)+\epsilon_J(\epsilon_{J'}(t))-t$.
Because  $n$ is a maximal weight  of $V_{n}(E)$, 
we deduce that $J=J'$. 

\medskip
\noindent{\bf Decomposition with respect to the central characters.}
 By  the lemma \ref{verma} $S^{i}_{\cO}(\mathfrak{g},\gb,n)$
admits  a finite filtration  by    $U_{ \cO}(\mathfrak{g})$-submodules, whose graded are of the
 form $V_{\cO}(\mu)$, $\mu\in \Omega^i(n)$. Therefore $S^{\bullet}_{\cO}(\mathfrak{g},\gb,n)$
is  annihilated by a power of the ideal  $I:=\prod_{\mu\in \Omega^{\bullet}(n)}\ker(\chi_{\mu,p})$ 
of  the commutative  ring   $U_{\cO}(\mathfrak{g})^G$.  As  we will see at the end  of this section as 
a consequence of Prop.\ref{chandra},  $S^{\bullet}_{\cO}(\mathfrak{g},\gb,n)$ is 
annihilated by $I$ itself. We have the following  commutative algebra lemma : 

\begin{lemma}
Let $P_1,..,P_r$  be ideals of the commutative ring   $R$, such that 
 $P_1..P_r=0$ and  for all $i\neq j$,  $P_i+P_j=R$. 
Then   each $R$-module $W$ admits a direct sum decomposition 
 $W=\oplus_{1\leq i\leq r} W^{P_i}$, with  $W^{P_i}=\{m\in W|P_im=0\}$.
\end{lemma}

Consider the maximal ideals   $pR+\ker(\chi_{\mu,p})=
\ker(\overline{\chi}_{\mu,p})$ of $R$, where $\mu\in \Omega^{\bullet}(n)$. 
Let  $\overline{\chi}_1=\overline{\chi}_{n,p}$,
$\overline{\chi}_2$,...,$\overline{\chi}_r$ be the 
set of distinct characters among 
$\overline{\chi}_{\mu,p}$, $\mu\in \Omega^{\bullet}(n)$.
Put $P_i=\prod_{\overline{\chi}_{\mu,p}=\overline{\chi}_i} 
\ker(\chi_{\mu,p})$. By the above lemma we get a decomposition 
\begin{equation}
S^{\bullet}_{\cO}(\mathfrak{g},\gb,n)=
\oplus_{i=1}^{r}
S^{\bullet}_{\cO}(\mathfrak{g},\gb,n)^{P_i}
\end{equation}
which is a direct sum, because  the differentials  are 
$U_{ \cO}(\mathfrak{g})$-equivariant. Moreover,
$V_{\cO}(\mu)_{\overline{\chi}_{n,p}}=
V_{\cO}(\mu)$, if $\overline{\chi}_{\mu,p}=\overline{\chi}_{n,p}$,
and  $V_{\cO}(\mu)_{\overline{\chi}_{n,p}}=0$, otherwise. 
From here and from Prop.\ref{chandra} we get  :

\begin{theo} \label{decomposition-centrale}
The complex  $S^{\bullet}_{\cO}(\mathfrak{g},\gb,
n)_{\overline{\chi}_{n,p}}$ is a direct factor  in 
$S^{\bullet}_{\cO}(\mathfrak{g},\gb,n)$
and we have  $S^{0}_{\cO}(\mathfrak{g},\gb,
n)_{\overline{\chi}_{n,p}}= V_{n}(\cO)$. 
For each $i\geq 1$, 
$S^{i}_{\cO}(\mathfrak{g},\gb,
n)_{\overline{\chi}_{n,p}}$ has a filtration whose 
graded are   given by the $V_{\cO}(\epsilon_J(n+t)-t)$
where $J\!\subset\! J_F$, $|J|=i$ (with multiplicity one).
\end{theo}

\subsection{The  BGG complex for distributions  algebras.}

Let $\cU_{\cO}(G)$ be the distribution $\cO$-algebra  over $G$.
For each $G$-module $V$,  free over  $\cO$,  we define the complex 
$$  0 \leftarrow  V \leftarrow 
\cU_{\cO}(G)\otimes_{\cU_{\cO}(B)}
\left(\wedge_{\cO}^{\bullet}(\mathfrak{g}/\gb)\otimes  
V|_{\gb}\right)  ,$$

\noindent and denote it by  $\mathcal{S}_{\cO}^{\bullet}(G,B,V)$. 
In the case where $V=V_{n}(\cO)$ we denote this complex by  
$\mathcal{S}_{\cO}^{\bullet}(G,B,n)$.

\begin{rque}
The complex $\mathcal{S}_{\cO}^{\bullet}(G,B,V)$ 
is not exact. It will become exact after applying the
Grothendieck linearization functor to the associated complex
of vector bundles over the Hilbert modular variety.
\end{rque}

For all $\mu\in \Z[J_F]$,  we define the Verma module $\cV(\mu)=
\cU_{\cO}(G)\otimes_{\cU_{\cO}(B)} W_{\mu}(\cO)$ (see \S\ref{BGG}).
We recall  that,  since  $p-1> |n|+d$,  
$\Omega^i(n)$ is the set  of  $\mu\in \Z[J_F]$
such that $W_{\mu}(\cO)$ is a irreducible subquotient  of 
$\wedge_{\cO}^{i}(\mathfrak{g}/\gb)\otimes  
V_{n}(\cO)|_{\gb}$. The  lemma \ref{verma} translates as  :

\begin{lemma} The modules
$\mathcal{S}_{\cO}^{\bullet}(G,B,n)$ has a finite filtration by 
 $\cU_{\cO}(G)$-submodules whose  successive quotients are given by
$\cV_{\cO}(\mu)$, with $\mu\in \Omega^i(n)$.
\end{lemma}

As $U_{\cO}(\mathfrak{g})\!\subset\!  
\cU_{\cO}(G)\!\subset\! U_{E}(\mathfrak{g})$,
the  center 
$U_{\cO}(\mathfrak{g})^G$  of $U_{\cO}(\mathfrak{g})$ 
is  contained in the center of $\cU_{\cO}(G)$.
Consider the central  characters  
$\chi_{\mu,p}=\chi_{\mu,\cO}$
and  $\overline{\chi}_{\mu,p}=\chi_{\mu,\overline{\F}_p}$ (see \S\ref{BGG}).

If  $W$ is a $\cU_{\cO}(G)$-module  generated by a vector 
$v$ of weight $\mu$, that  is  annihilated by $\gu$,
then  $U_{\cO}(\mathfrak{g})^G$ acts  on $W$ by the character $\chi_{\mu,p}$. 
Put  $I=\prod_{\mu\in \Omega^{\bullet}(n)}\ker(\chi_{\mu,p})$. 
By  the  last lemma  the finite ${\cO}$-module  
$\mathcal{S}_{\cO}^{\bullet}(G,B,n)$ is  a  
$R:=U_{\cO}(\mathfrak{g})^G/I$-module. 
Let  $\overline{\chi}_1=\overline{\chi}_{n,p},\overline{\chi}_2,...
\overline{\chi}_r$ be the distinct  algebra homomorphisms
from  $R$ in $\overline{\F}_p$. For $1\leq j\leq r$, we put
$$\mathcal{S}_{\cO}^{\bullet}(G,B,n)_{\overline{\chi}_j}=
\left\{x\in\mathcal{S}_{\cO}^{\bullet}(G,B,n)
\Big{|} \left(\prod_{\mu\in \Omega^{\bullet}(n), \overline{\chi}_{\mu,p} 
=\overline{\chi}_j}  \ker(\chi_{\mu,p})\right)x=0\right\}.$$

The same way as in Thm.\ref{decomposition-centrale}  we obtain a decomposition 
\begin{equation}\label{decomposition}
\mathcal{S}_{\cO}^{\bullet}(G,B,n)= \oplus_{j=1}^{r}
\mathcal{S}_{\cO}^{\bullet}(G,B,n)_{\overline{\chi}_j}.
\end{equation}

The main theorem of this section states then :

\begin{theo}
$\displaystyle \mathcal{S}_{\cO}^{i}(G,B,n)_{\overline{\chi}_{n,p}}
\cong \bigoplus_{J\!\subset\! J_F, |J|=i}\cV_{\cO} (\epsilon_J(n+t)-t).$
\end{theo}

{\bf Proof : } We start with the case $n=0$.

As $\gu$ is abelian, $U$ acts trivially on  
$\wedge_{\cO}^{i}(\mathfrak{g}/\gb)$ and therefore, 
by  the lemma \ref{verma}, we have
$$\wedge_{\cO}^{i}(\mathfrak{g}/\gb)\cong
\bigoplus_{J\!\subset\! J_F, |J|=i} W_{\epsilon_J(t)-t}(\cO). $$ 

As  $\cU_{\cO}(G)$ is free over $\cU_{\cO}(B)$ we obtain :
\begin{equation}
\mathcal{S}_{\cO}^{i}(G,B,0)=
\mathcal{S}_{\cO}^{i}(G,B,0)_{\overline{\chi}_{0,p}}
\cong \bigoplus_{J\!\subset\! J_F, |J|=i}\cV_{\cO}
(\epsilon_J(t)-t).
\end{equation}

If $n\geq 0$,  we deduce from the $n=0$  case a  decomposition
\begin{equation}
\mathcal{S}_{\cO}^{i}(G,B,n)\cong \bigoplus_{J\!\subset\! J_F, |J|=i} 
\cU_{\cO}(G)\otimes_{\cU_{\cO}(B)}
\left( W_{\epsilon_J(t)-t}(\cO)\otimes V_{n}(\cO) \right).
\end{equation}

Using  (\ref{decomposition}) the theorem is a consequence of the 
following lemma, whose proof is a direct application of the proof 
 of Prop.\ref{chandra}.

\begin{lemma} 
$\left(\cU_{\cO}(G)\otimes_{\cU_{\cO}(B)}
\left( W_{\epsilon_J(t)-t}(\cO)\otimes
V_{n}(\cO) \right)\right)_{\overline{\chi}_{n,p}}
\cong \cV_{\cO}(\epsilon_J(n+t)-t).$
\end{lemma}

\subsection{BGG complex for crystals.}
Our reference is the section 4 of  \cite{MoTi}.
For every integer $r\geq 0$ we put  $S_r=\Spec(\Z/p^{r+1})$.
For a $\Z[\frac{1}{\Delta}]$-scheme $X$,   we put  $X_r=X \times S_r$. 

We have  an equivalence of categories between the category of crystals over 
$(\overline{X}_0 /S_r)_{\log}^{\crys}$ and the category of
$\cO_{\overline{X}_r}$-modules $\cM$ 
which are locally free and   endowed with integrable, quasi-unipotent
connection with logarithmic poles  $\nabla : \cM \rightarrow \cM
\otimes_{\cO_{\overline{X}_r}} \Omega^1_{\overline{X}_r/S_r}(\dlog(\infty_X))$.

\medskip
We have  a functor $\mathrm{L}$, called the  {\it linearization functor},  
from  the category of sheaves of  
 $\cO_{\overline{X}_r}$-modules to  the category of 
 crystals on  $(\overline{X}_0 /S_r)_{\log}^{\crys}$.

By the log-crystalline Poincar{\'e} lemma, we have a resolution :
\begin{equation}\label{poincare-crys}
0 \rightarrow \cM \rightarrow
\mathrm{L}(\cM\otimes_{\cO_{\overline{X}_r}} 
\Omega^{\bullet}_{\overline{X}_r/S_r}(\dlog \infty)).
\end{equation}

Let  $W_1$ and  $W_2$ be two $B$-modules with weights smaller than $ (p-1)t$. 
Put $\overline{\cW}_i=\cF_{B}(W_i)$, $i=1,2$ (see \S\ref{hodge-tate}).
By   \cite{MoTi}{\S}5.2.4   we have a homomorphism 
\begin{equation}
\Hom_{\cU_{\cO}(G)} (\cU_{\cO}(G) \otimes_{\cU_{\cO}(B)}W_1), \cU_{\cO}(G)
\otimes_{\cU_{\cO}(B)}W_2)\rightarrow 
\mathrm{Diff.Op.}(\overline{\cW}_{2,r},\overline{\cW}_{1,r}),
\end{equation}
which becomes  an isomorphism  after tensoring with $E$ (see (\ref{op-diff})).

\medskip
We apply now the above construction to the toroidal compactification 
over of  the Hilbert modular variety $\overline{M}{}'$ and 
the vector bundle $\overline{\cV}_n$. For every  $r\geq 0$ we have an injective 
 homomorphism  of complexes  of vector bundles  over $\overline{M}_r'$   
\begin{equation}\label{iota}
 \mathcal{K}^{\bullet}_n:=
\underset{J\!\subset\! J_F}{\oplus}\overline{\cW}_{\epsilon_J(n+t)-t}
\hookrightarrow \overline{\cV}_n \otimes_{\cO_{\overline{M}_r'}}
\Omega^{\bullet}_{\overline{M}_r'/S_r}(\dlog \infty).
\end{equation}

\begin{prop}\label{qis}
The map $(\ref{iota})$ is a strict  injective homomorphism  of filtered complexes.
\end{prop}

By the last  proposition   $\mathrm{L}( \mathcal{K}^{\bullet}_n)$
is a direct factor in $\mathrm{L}(\overline{\cV}_n \otimes_{\cO_{\overline{M}_r'}} 
\Omega^{\bullet}_{\overline{M}_r'/S_r}(\dlog \infty))$, which is 
 exact by  the Poincar{\'e}'s crystalline lemma. Therefore
 $\mathrm{L}( \mathcal{K}^{\bullet}_n)$ is also exact. 
As the functor $\mathrm{L}$ is exact, we deduce filtered  isomorphisms
$\cH_{\ldR}^j(\overline{M}_r'/S_r,  \overline{\cV}_n)
\cong \cH^j(\overline{M}_r'/S_r, \mathcal{K}^{\bullet}_n)$.

Recall that $p$ does not divide $\Delta$ and $p-1>|n|+d$. Under this 
assumptions we have

\begin{theo}\label{bgg-modp}
The spectral sequence  given by the  Hodge filtration 
$$\rE_1^{i,j}=\bigoplus_{J \!\subset\! J_F, |p(J)|=i} 
\rH^{i+j-|J|} (\overline{M}_r, \overline{\cW}_{\epsilon_J(n+t)-t,n_0}) 
\Rightarrow\cH_{\ldR}^{i+j}(\overline{M}_r,\overline{\cV}_n)$$
degenerates at $\rE_1$  :
\begin{equation}\label{bgg-modp1}
\gr^i \cH_{\ldR}^r(\overline{M}_r, \overline{\cV}_n) 
=\bigoplus_{J \subset J_F, |J|\leq r, |p(J)|=i} \rH^{r-|J|}
(\overline{M}_r,  \overline{\cW}_{\epsilon_J(n+t)-t}).
\end{equation}
\end{theo}

{\bf Proof : } The proof  is formally  the same as the one 
of Thm.\ref{hodge-pad}(ii), once we have Prop.\ref{qis}.
The degeneration of the spectral sequence follows from a 
result of Illusie \cite{illusie2} Prop.4.13. applied to the
semi-stable morphism $\overline{\pi_s}:\overline{\cA^s}\rightarrow
\overline{M^1}$ of smooth $\Z_p$-schemes.
\hfill $\square$

\begin{rque} \label{hecke-equiv}
(i) It follows from the same arguments as in 
Cor.\ref{hodge-complex}(i), that the above decomposition is 
 Hecke equivariant, except for the $T_{\gp}$ operators, when 
 $\gp$ divides $p$. When $p$ is totally split in $F$,  
we could use Wedhorn's results \cite{Wed}  to write
$T_{\gp}$ as a sum of correspondences and try adapt to this case the method of 
\cite{FaJo}. Unfortunately, this approach is not available 
when $p$ is not  totally split in $F$. 

In the  proof of Thm.\ref{theo-B}, we will use  different method to 
prove the $T_{\gp}$-equivariance 
of the  above decomposition after a localization outside $p$. 

(ii) The commutativity of the  Hecke operators outside $p$ follows  from the  
degeneration at $\rE_1$ as in the proof of Cor.\ref{hodge-complex}(i). The last  
graded piece  $\rH^0(\overline{Y},\overline{\cW}_{\epsilon_{J_F}(n+t)-t,n_0})$ 
of the  filtration is   independent of the toroidal compactification
 by the Koecher principle (\ref{koecher}).
\end{rque}

\section{Integral cohomology  over  certain local  components  of 
the Hecke algebra.}

\subsection{The  key lemma.}
Let $q=p^r$ and denote by  $\sigma_1$,..,$\sigma_r$ the elements of
$\Gal(\F_q/\F_p)$.

\begin{theo} {\rm(Brauer-Nesbitt, Steinberg \cite{St})}
The group $\SL_2(\F_q)$ has exactly  $q$ irreducible representations 
 on finite dimensional  $\F_q$-vector spaces, namely  the
$\otimes_{j=1}^r(\Sym^{a_j})^{\sigma_j}$, for $0\leq a_j\leq p-1$.
\end{theo}

\begin{cor}\label{bn}
For every finite set  $I$, the group 
$ \prod_{i\in I}\SL_2(\F_q)$ has exactly  $q^{|I|}$ irreducible 
representations on finite dimensional  $\F_q$-vector spaces, namely the
$ \otimes_{i\in I}
\left(\otimes_{j=1}^r(\Sym_i^{a_{i,j}})^{\sigma_j}\right)$, 
for  $0\leq a_{i,j}\leq p-1$.
\end{cor}

In  \cite{Ma}  Mazur  states the following :

\begin{lemma} Let $\Phi$ be a group 
and $\rho_0$ be  a  representation of $\Phi$
on a  finite dimensional $\F_q$-vector space $W$. 
Let  $\rho: \Phi \rightarrow \GL_2(\F_q)$ be an absolutely irreducible 
representation such that  for all $g\in \Phi$, the characteristic polynomial
 of  $\rho(g)$ annihilates $\rho_0(g)$. Then,
$\rho_0^{\mathrm{s.s.}}=\rho\oplus..\oplus\rho$ and in particular  
$\rho\!\subset\!\rho_0$.
\end{lemma}

The  corresponding statement 
for a another group than $\GL_2$ is false  in  general.
Here is an example for $\GL_3$ : take 
$\rho=\Sym^2: \GL_2(\F_q)\rightarrow \GL_3(\F_q)$ and 
$\rho_0=\det:\GL_2(\F_q)\rightarrow \GL_1(\F_q)$.
Nevertheless, a we have a generalization for the special group :
$$H(\F_q)=\left(\prod_{i\in I}\GL_2(\F_q)\right)^\mathcal{D}:=
\left\{ (M_i)_{i\in I}\in  \prod_{i\in I}\GL_2(\F_q)\Big{|}\enspace\exists
 \delta\in\mathcal{D},\enspace \forall i\in I, \enspace  \det(M_i)=\delta\right\}$$

and the particular  representation 
$$\rho_1= \bigotimes_{i\in I,\tau\in J_F^i}
\mathrm{St}_i^{\sigma_{i,\tau}} : H(\F_q) \rightarrow
  \GL_{2^d}(\F_q) \text{ , } (M_i)_{i\in I}\mapsto 
\bigotimes_{i\in I,\tau\in J_F^i} M_i^{\sigma_{i,\tau}},$$ 
where $(J_F^i)_{i\in I}$ is a partition of $J_F$ and for all  $i\in I$,
 $(\sigma_{i,\tau})_{\tau\in J_F^i}$ are  two by two distinct elements 
of $\Gal(\F_q/\F_p)$ ($\St=\Sym^1$ is    
  the standard two-dimensional representation  of $\GL_2$).   
\begin{lemma}\label{variante-mazur}
Let $\rho_0$ be  a  representation of $H(\F_q)$
on a  finite dimensional $\F_q$-vector space $W$, such that 
 for all $g\in H(\F_q)$ the characteristic polynomial  
of $\rho_1(g)$ annihilates $\rho_0(g)$. Then 
$\rho_0^{\mathrm{s.s.}}={\rho_1}\oplus..\oplus{\rho_1}$ (each irreducible 
subquotient  of $\rho_0$ is isomorphic to $\rho_1$).
\end{lemma}

 {\bf Proof : } We can assume that $\rho_0$ is absolutely irreducible. 
Consider the exact 
sequence  $1\!\rightarrow\! H_1(\F_q)=\prod_{i\in I}\SL_2(\F_q)\!\rightarrow\!
 H(\F_q)\! \overset{\nu}{\rightarrow} \!\mathcal{D} \!\rightarrow\! 1$.
  By Cor.\ref{bn}, 
we know that  each irreducible  subquotient  of ${\rho_0}_{|H_1(\F_q)}$ is of 
the form   $\otimes_{i\in I}\left(\otimes_{j=1}^r(\Sym_i^{a_{i,j}})^{\sigma_j}
\right) \text{, with } 0\leq a_{i,j}\leq p-1$.

The subspace corresponding to the highest weight of the representation
${\rho_0}_{|H_1(\F_q)}$ is preserved by the  standard torus of $H(\F_q)$, and 
therefore contains an eigenvector $x$
for the action of this torus. Because $\rho_0$ is irreducible, it is 
generated by $x$, and therefore $\rho_0$ isomorphic to a  twist of 
$\otimes_{i\in I}\left(\otimes_{j=1}^r(\Sym_i^{a_{i,j}})^{\sigma_j}\right)$
 by some power  of the character $\nu$
(in particular, ${\rho_0}_{|H_1(\F_q)}$ is also irreducible).

 As the characteristic polynomial   of $\rho_1$ annihilates $\rho_0$, 
the set  of the  weights of $\rho_0$ is a subset of the set  of the  weights of 
$\rho_1$, and therefore $\rho_0=\rho_1$. \hfill$\square$

In {\S}\ref{ind-image} we proved  under the assumption
 $\mathrm{\bf (LI_{\Ind\overline{\rho}})}$, that 
$\Ind_F^{\Q}\overline{\rho}(\G_{\widehat{F}})$ contains the image of the map 
$\phi=(\phi^i)_{i\in I} : H(\F_q) \hookrightarrow  \GL_2(\F_q)^{J_F}$.

Denote by $\widehat{F}'$ the fixed field of 
$\overline{\rho}^{-1}(\phi(H(\F_q)))$.

\begin{lemma} {\rm(Key Lemma)} \label{key}
Let $\rho_0$ be  a  representation of $\G_{\widehat{F}'}$ on a 
finite dimensional $\kappa$-vector space $W$. 
Assume $\mathrm{\bf (LI_{\Ind\overline{\rho}})}$ and assume that, 
for every $g\in \G_{\widehat{F}'}$,  the characteristic polynomial  of 
$(\otimes\Ind_F^{\Q}\overline{\rho})(g)$ annihilates $\rho_0(g)$.
Then  each $\G_{\widehat{F}'}$-irreducible subquotient   of $\rho_0$ 
is isomorphic to  $\otimes\Ind_F^{\Q}\overline{\rho}$.
\end{lemma}

 {\bf Proof : } It is enough to treat the case where $\rho_0$ is irreducible. 
The idea is show that the action of $\G_{\widehat{F}'}$ on $W$ is through 
the  algebraic group $H(\F_q)$, and then  use the lemma 
$\ref{variante-mazur}$.

Put  $\overline{\rho}'= (\Ind_F^{\Q}\overline{\rho})_{|\G_{\widehat{F}'}}$. 
Because of the annihilation assumption, the group $\rho_0(\ker(\overline{\rho}'))$
is an unipotent $p$-group   and  therefore $W^{\ker(\overline{\rho}')}$
is non-zero.  Moreover  the subspace $W^{\ker(\overline{\rho}')}$ is 
preserved  by $\G_{\widehat{F}'}$. Because $W$ is irreducible, we get 
$W^{\ker(\overline{\rho}')}=W$ and therefore the action of $\G_{\widehat{F}'}$
 on  $W$ is through $H(\F_q)$.  Thus we get a homomorphism  $\rho_0'$  
fitting in the following commutative diagram :
$$\xymatrix@C=50pt{
 \G_{\Q}\ar[r]^{\otimes\Ind_F^{\Q}\overline{\rho}}&  \GL_{2^d}(\kappa)\\
  \G_{\widetilde{F}}\ar@{_{(}->}[u] \ar[r]^{\Ind_F^{\Q}\overline{\rho}}& 
 \GL_2(\kappa)^{J_F} \ar[u]^{\otimes} \\
  \G_{\widehat{F}'} \ar[d]_{\rho_0}
\ar@{_{(}->}[u]\ar@{->>}[r]^-{\phi^{-1}\circ\overline{\rho}'}& 
H(\F_q)\ar@{_{(}->}^{\phi}[u]
\ar@{-->}[dl]_{\rho_0'}\ar@{_{(}->}[u] \ar@/_3pc/[uu]_{\rho_1} \\
\GL(W)  & }$$
The characteristic polynomial  of   $\rho_1$ annihilates   the representation 
$\rho_0'$.  By the lemma \ref{variante-mazur} each $H$-irreducible 
subquotient of $W$ is  isomorphic to $\rho_1$, that is  to say
$W^{\mathrm{s.s.}}= \oplus\rho_1$ as  $H(\F_q)$-modules. As the action of 
$\G_{\widehat{F}'}$ on both sides is  through $H(\F_q)$, we're done. 
\hfill$\square$

\subsection{Localized cohomology of the  Hilbert modular variety.}${}$

Let $\T'\subset \T$ be the subalgebra generated by the Hecke operators
outside a finite set of places containing those dividing $\gn p$. Put
$\gm'=\gm\cap \T'$.

\begin{theo}\label{coh-loc}
Assume  $f$ and  $p$ satisfy $\mathrm{{\bf (I)}}$, $\mathrm{{\bf (II)}}$ and 
$\mathrm{\bf (LI_{\Ind\overline{\rho}})}$. Then 

{\rm (i)}  $\rH^\bullet(Y,\bV_n(\kappa))_{\gm'}=
\rH^d(Y,\bV_n(\kappa))_{\gm'}$,

{\rm (ii)} 
$\rH^\bullet(Y,\bV_n(\cO)_{\gm'}
=\rH^d(Y,\bV_n(\cO))_{\gm'}$ is a free $\cO$-module 
of  finite rank and the $\cO$-module  $\rH^\bullet(Y,\bV_n(E/\!\cO))_{\gm'}
=\rH^d(Y,\bV_n(E/\!\cO))_{\gm'}$ is  divisible   of finite corank .

$\mathrm{(iii)}$ 
$\rH^d(Y,\bV_n(\cO))_{\gm'} 
\times \rH^d(Y,\bV_n(E/\!\cO))_{\gm'}\rightarrow \cO$ is a
perfect Pontryagin pairing.
\end{theo}

{\bf Proof : } 
(i) By Faltings' Comparison Theorem  \cite{Fa-jami} and Thm.\ref{bgg-modp}(i) the integer
$|p(J)|$ is not a Fontaine-Laffaille  weight 
of  $\rH^r(\kappa)$, when $r<d$.  Wedhorn \cite{Wed} has established the 
congruence  relations for  all totally split   primes  of  $F$.
By  the Cebotarev Density Theorem  the assumptions of the key lemma 
\ref{key} are fulfilled. We deduce  that  $\rH^r(\kappa)[\gm']=0$, and therefore 
$\rH^r(\kappa)_{\gm'}=0$
by Nakayama's lemma. The case  $n>d$ follows by  Poincar{\'e} duality.

(ii)(iii)  By  the long exact cohomology sequence
$$...\rightarrow \rH^{r-1}(\kappa)\longrightarrow
\rH^r(\cO)\overset{\varpi}{\longrightarrow}
\rH^r(\cO)\longrightarrow \rH^r(\kappa)\rightarrow  ... ,$$

\noindent and by the vanishing  of $\rH^r(\kappa)_{\gm'}$, 
for  $r \neq d$, we  deduce that (for $r \neq d$)
the  multiplication by an uniformizer $\varpi$ is a surjective endomorphism  of 
$\rH^r(\cO)_{\gm'}$, so this last vanishes.

The  same way, by the long exact sequence
$$...\rightarrow \rH^r(\varpi^{-1}\cO/\cO)\longrightarrow
\rH^r(E/\!\cO)\overset{\varpi}{\longrightarrow}
\rH^r(E/\!\cO)\longrightarrow
\rH^{r+1}(\varpi^{-1}\cO/\cO) \rightarrow  ... ,$$

\noindent we deduce a surjection $\rH^r(\kappa)_{\gm'}
\twoheadrightarrow  \rH^r(E/\!\cO)_{\gm'}[\varpi]$, when $r\neq d$. 
Since  $\rH^r(E/\!\cO)_{\gm'}$ is a torsion $\cO$-module, it
vanishes (for $r \neq d$).

\medskip 

The  localization at $\gm'$ of the  long exact  sequence  of $\cO$-modules :
$$...\rightarrow \rH^{r-1}(E/\!\cO)\longrightarrow
\rH^r(\cO)\longrightarrow \rH^r(E)\longrightarrow \rH^r(E/\!\cO)
\rightarrow... ,$$
 is concentrated  at the three  terms of degree  $r=d$.
From this we deduce the freeness. \hfill $\square$

\subsection{On the Gorensteiness of Hecke the algebra.}

\begin{theo}{\rm(Theorem B)}\label{theo-B}
Let  $f$ and  $p$ be such that $\mathrm{{\bf (I)}}$, $\mathrm{{\bf (II)}}$ 
and $\mathrm{\bf (LI_{\Ind\overline{\rho}})}$  hold. Then

{\rm(i)} $\rH^\bullet(Y,\bV_n(\kappa))[\gm]=
\rH^d(Y,\bV_n(\kappa))[\gm]$ is a  $\kappa$-vector 
space  of dimension $2^d$. 

{\rm(ii)} $\rH^\bullet(Y,\bV_n(\cO))_{\gm}=
\rH^d(Y,\bV_n(\cO))_{\gm}$ 
is free of rank  $2^d$ over  $\T_{\gm}$.  

{\rm(iii)} $\T_{\gm}$ is Gorenstein.
\end{theo}

{\bf Proof : } In this proof we put 
$W=\rH^d(Y_{\overline{\Q}},\bV_n(\kappa))_{\gm}$. 
By using an auxiliary level structure as in \cite{diamond2}, we
can assume that the condition {\bf (NT)} of \S\ref{hmv} is fullfilled. 

(i) As in the proof of Thm.\ref{coh-loc}(i), by  lemma \ref{key} 
we have an isomorphism of $\G_{\widehat{F}'}$-modules
$$ W[{\gm}]^{\mathrm{s.s.}}=
(\otimes \Ind_F^{\Q}\overline{\rho})^{\oplus r}.$$
It is crucial   to observe that  $I_p\subset \G_{\widehat{F}'}$.
By Thm.\ref{BL-thm} we have $r\geq 1$.
In order to show that  $r=1$  we consider  the restriction of these 
representations  to $I_p$. 
The multiplicity of the 
maximal Fontaine-Laffaille weight $|p(J_F)|$ in the right hand side  is $r$, 
by Thm.\ref{BL-thm}, Cor.\ref{hodge-complex}(ii)
and Fontaine-Laffaille's theory. 

On the other hand, the multiplicity of $|p(J_F)|$ in the left hand side 
is equal, by Thm.\ref{bgg-modp},  to the dimension of 
$\rH^0(\overline{Y}\otimes\kappa,\overline{\cW}_{\epsilon_{J_F}(n+t)-t,n_0})[\gm]$ :
In fact, by the remark \ref{hecke-equiv} it is sufficient to check 
the $T_{\gp}$-equivariance of the $\gm'$-localization of the projective limit over 
$r$ of (\ref{bgg-modp1}). By Thm.\ref{coh-loc}(ii) it is question of 
checking the $T_{\gp}$-equivariance of an isomorphism of {\it free} $\cO$-modules.
 Therefore, it is enough to be checked after  extending  the  scalars to 
$\C$. Then the  Strong Multiplicity One Theorem applies (recall
that $p$ is prime to the level $\gn$). We owe  this  idea to  Diamond  
(\cite{diamond2} proof of Prop.1).

\smallskip
We will now see that  $\dim_\kappa \rH^0(\overline{Y}\otimes \kappa, 
\overline{\cW}_{\epsilon_{J_F}(n+t)-t,n_0})[\gm]=1$.
We have $\overline{\cW}_{\epsilon_{J_F}(n+t)-t,n_0 }=
\underline{\omega}^k\otimes\underline{\nu}^{n_0t/2}$.
So we are led to show that two normalized Hilbert modular forms of weight $k$, 
level $\gn$ and coefficients in  $\kappa=\T_{\gm}/\gm$ 
 having  the same  eigenvalues for all the Hecke operators  are equal.
One should be careful to observe  that the Hecke operators 
permute the  connected components $M_1(\gc,\gn)$
 of the Shimura variety     $Y=Y_1(\gn)$ (here 
the ideal $\gc$ runs over a set of representatives of  $\Cl_F^+$).
We use then the Hecke relations  between Fourier coefficients 
and  eigenvalues  for the Hecke operators   and 
the   $q$-expansion principle (see \ref{qdev})  at the  $\infty$ cusp  
of each connected  component  $M_1(\gc,\gn)$.

\medskip

Even if we do not know the degeneration at 
$\rE_1$ of the   Hodge to  de Rham spectral sequence, we 
obtain by the same  arguments that $r\leq 1$ (instead of  $r=1$), 
because  we have always $\rH^0(\overline{Y}\otimes \kappa, 
\overline{\cW}_{\epsilon_{J_F}(n+t)-t,n_0})[\gm]
\supset  \gr^{|p(J_F)|}W[\gm]$. But 
$\rH^\bullet(Y,\bV_n(\cO))_{\gm}$ is non zero 
 as $\rH^\bullet(Y,\bV_n(\cO))_{\gm}\otimes \Q$
is  free of rank $2^d$ over $\T_{\gm}\otimes \Q$, and therefore $r=1$.

(ii)(iii) Mazur's argument in the elliptic modular case remains valid.
By the theorem A, the twisted Poincar{\'e} pairing (\ref{accoupl-twisted}) on 
$\rH^d(Y,\bV_n(\cO))_{\gm}= \rH_c^d(Y,\bV_n(\cO))_{\gm}$ 
 is  a perfect duality of $\T_{\gm}$-modules,  
so it would be enough to show (ii).

Again using  the perfectness of the twisted Poincar{\'e} pairing
$W\times W\rightarrow \kappa$ 
we obtain  $W \cong \Hom_{\T_{\gm}} (W,\kappa)$, and so 
$W\otimes_{\T_{\gm}} \kappa  =W/\gm W \cong \Hom(W[\gm],\kappa)$,
and therefore 
$$\dim_\kappa(W\otimes_{\T_{\gm}} k)= \dim_\kappa(W[\gm]),$$
which equals  $2^d$, by (i). Then (ii) follows from the following

\begin{lemma}\label{freeness}
Let $\cT$ a torsion free local $\cO$-algebra 
($\cT\hookrightarrow \cT\otimes_{\cO}E$) of  maximal ideal 
 $\gm$ and  residue field $\kappa=\cT/\gm$.

Let $M$ be a finitely generated $\cT$-module such that  $M\otimes_{\cO}E$
is  free of rank $r$ over  $\cT\otimes_{\cO}E$. If
 $M\otimes_{\cT} \kappa$ is a $\kappa $-vector space  of dimension $\leq r$, 
then  $M$ is free of rank $r$ over $\cT$. 
\end{lemma}

{\bf Proof : } 
Since  $\cM\otimes_{\cT} k$ is of dimension $\leq r$, the  Nakayama's  lemma 
gives  a surjective homomorphism  of $\cT$-modules
$\cT^r\twoheadrightarrow M$. Denote by  $I$ its kernel. We have an 
 exact  sequence of $\cO$-modules 
$$0\rightarrow I \rightarrow \cT^r \rightarrow M \rightarrow 0.$$
By  tensoring it by  $\otimes_{\cO}E$ (or equivalently by  
$\otimes_{\cT}(\cT\otimes_{\cO}E)$) 
we obtain another  exact sequence   $$0\rightarrow I \otimes_{\cO}E
 \rightarrow (\cT\otimes_{\cO}E)^r \rightarrow M\otimes_{\cO}E\rightarrow 0.$$
By comparing the  dimensions over $E$ we get  $I\otimes_{\cO}E=0$. 
Since $I$ is torsion free, $I=0$. \hfill$\square$

\subsection{An application to $p$-adic ordinary families.}

For $r\geq 1$, consider the following open compact subgroups of $G(\A_f)$
$$K_0(p^r)=\left\{u\in K_1(\gn)| u\equiv 
\begin{pmatrix} *& * \\ 0 & *\end{pmatrix}\pmod{p^r}\right\},$$
$$K_{11}(p^r)=\left\{u\in K_1(\gn)| u\equiv 
\begin{pmatrix} 1 & * \\ 0 & 1\end{pmatrix}\pmod{p^r}\right\}.$$

Let $Y_0(p^r)$ (resp. $Y_{11}(p^r)$) be the Hilbert modular variety of 
level $K_0(p^r)$ (resp. $K_{11}(p^r)$).

\medskip
The cohomology group $\rH^{\bullet}(Y_{11}(p^r),\bV_n(E/\cO))^*$ has 
a natural action of $K_0(p^r)/K_{11}(p^r)\simeq
(\go/p^r)^\times\times (\go/p^r)^\times$
(we denote by ${}^*$ the Pontryagin dual). Therefore the group 
$T(\Z_p)/\overline{\go^\times}$ acts on 
the inductive  limit $\rH^{\bullet}(Y_{11}(p^{\infty}),\bV_n(E/\cO))^*:=
\limind \rH^{\bullet}(Y_{11}(p^r),\bV_n(E/\cO))^*$.

By Hida's stabilization lemma, the ordinary part of 
$\rH^{\bullet}(Y_{11}(p^{\infty}),\bV_n(E/\cO))^*$ (that is the part
where the Hecke operators $T_{0,\gp}$ of Def.\ref{ord-autom} 
are  invertible for all  $\gp$ dividing $p$) is independent on $n$. We
denote it by  $\cH_{\ord}^{\bullet}:=
\rH^{\bullet}_{\ord}(Y_{11}(p^{\infty}),E/\cO)^*$.

By the above discussion    $\cH_{\ord}^{\bullet}$ has a structure of a 
 $\Lambda:=\cO[[T(\Z_p)/\overline{\go^\times}]]$-module.
 It is of finite type, by a theorem of Hida. 

\medskip
We also define the $p$-adic ordinary Hecke $\Lambda$-algebra 
$\T_{k,\ord}^{\infty}:=\limproj \T_{k,\ord}(Y_{11}(p^r))$. 
As $\T_{k,\ord}^{\infty}$ is independant of $k$, we denote it by 
$\T^{\infty}_{\ord}$. 
Then  $\cH_{\ord}^{\bullet}$ is  a  $\T^{\infty}_{\ord}$-module.

\medskip
An  arithmetic  character of $T(\Z_p)/\overline{\go^\times}$ is by definition 
a character whose restriction to an open subgroup is given by an 
algebraic character. It is immediate that such a character is a 
product of an algebraic character and a finite order character. An algebraic 
character of $T(\Z_p)\simeq D(\Z_p)\times D(\Z_p)$ trivial on $\go^\times$
is necessarily of the form $(u,\epsilon)\mapsto u^n\epsilon^{-m}$, 
where $m,n\in \Z[J_F]$ and
$m+2n\in \Z t$. Hence, the general form of an arithmetic character $\psi$
of  $T(\Z_p)/\overline{\go^\times}$ is   $(u,\epsilon)\mapsto 
u^n\epsilon^{-m}\psi_1(u)\psi_2(\epsilon)$, where $\psi_1,\psi_2$ are finite 
order characters.  Every such $\psi$ induces an $\cO$-algebra homomorphism 
$\Lambda\rightarrow \cO$, whose kernel is denoted by $P_\psi$. 

\medskip
Let $\gm$ be a maximal ordinary ideal of $\T=\T_k(\gn)$ and  $\gm_{\infty}$ be 
a maximal ideal of $\T^{\infty}_{\ord}$ above  $\gm$.
We denote by $\T^{\infty}_{\gm_{\infty}}$
(resp. $\cH_{\gm_{\infty}}^{\bullet}$) the localization of 
$\T^{\infty}_{\ord}$ (resp. of  $\cH_{\ord}^{\bullet}$) at
$\gm_{\infty}$.

\begin{prop} \label{control} Let  $\gm$ be such  that $\mathrm{{\bf (I)}}$, $\mathrm{{\bf (II)}}$ and
 $\mathrm{\bf (LI_{\Ind\overline{\rho}})}$ hold. Then 

\rm{(i)} $\cH_{\gm_{\infty}}^d$ is free of finite rank over $\Lambda$ and
we have exact control :
 $$\cH_{\gm_{\infty}}^d/P_\psi\cH_{\gm_{\infty}}^d\simeq
\rH^{\bullet}(Y_{11}(p^r),\bV_\psi(E/\cO))^*_{\gm_r},$$

\rm{(ii)} $\cH_{\gm_{\infty}}^d$ is free of  rank $2^d$ over 
$\T^{\infty}_{\gm_{\infty}}$, and 

\rm{(iii)}  Hida's control theorem for the  Hecke algebra holds, that is 
$\T^{\infty}_{\gm_{\infty}}$ is a free $\Lambda$-algebra  of finite rank
and  for every $\psi$ we have
$\T^{\infty}_{\gm_{\infty}}/P_\psi\T^{\infty}_{\gm_{\infty}}\simeq 
\T_\psi(Y_{11}(p^r))_{\gm_r}$.
\end{prop}

{\bf Proof : } 
(i) The proof is very similar to the one of \cite{MoTi}Thm.9. It uses 
that a $\Lambda$-module is free, if it is free of constant rank over $\cO$ for
infinitely many specializations. In our case, it is enough to 
specialize at the weight of the form $k+(p-1)k'$, and use 
the exact control criterion and Thm.\ref{coh-loc}. 
We omit the details. Note that (i) follows from  (ii) and (iii).

(ii) Consider   $\Lambda\rightarrow \T^{\infty}_{\gm_{\infty}}
\rightarrow \End_{\cO}(\cH_{\gm_{\infty}}^d)$. The specialization 
 at $\psi=\psi_k$ gives 
$$\cO \rightarrow \T^{\infty}_{\gm_{\infty}}/P_k \T^{\infty}_{\gm_{\infty}}
\rightarrow \End_{\cO}(\cH_{\gm_{\infty}}^d/P_k \cH_{\gm_{\infty}}^d).$$

By Thm.\ref{coh-loc} we have 
$\rH^d(Y_0(p),\bV_n(E/\cO))^*_{\gm}\simeq \rH^d(Y,\bV_n(E/\cO))^*_{\gm}\simeq
\rH^d(Y,\bV_n(\cO))_{\gm}$ and an exact control : $\cH_{\gm_{\infty}}^d/P_k 
\cH_{\gm_{\infty}}^d\simeq  \rH^d(Y,\bV_n(\cO))_{\gm}$.

From here and from  Thm.B we obtain that $\cH_{\gm_{\infty}}^d\otimes_{\T^{\infty}_{\gm_{\infty}}}
(\T^{\infty}_{\gm_{\infty}}/P_\psi\T^{\infty}_{\gm_{\infty}})\simeq
\cH_{\gm_{\infty}}^d\otimes_{\Lambda} \Lambda/P_k$ 
is free  of rank $2^d$ over  $\T_{\gm}$. Hence 
$\cH_{\gm_{\infty}}^d\otimes_{\T^{\infty}_{\gm_{\infty}}}\kappa$
is free  of rank $2^d$ over  $\T_{\gm}
\otimes_{\T^{\infty}_{\gm_{\infty}}/P_\psi\T^{\infty}_{\gm_{\infty}}}\kappa=\kappa$.
Then lemma \ref{freeness}  applies to the $\T^{\infty}_{\gm_{\infty}}$-module
$\cH_{\gm_{\infty}}^d$ which is finitely generated  over the local algebra $\Lambda$.

(iii) As  $\cH_{\gm_{\infty}}^d$ is a free  $\Lambda$-module, it 
admits a direct sum decomposition with respect to the Weyl group
action on the Betti cohomology :
$$\cH_{\gm_{\infty}}^d=\bigoplus_{J\subset J_F} 
\cH_{\gm_{\infty}}^d[\widehat{\epsilon_J}].$$

Every $\cH_{\gm_{\infty}}^d[\widehat{\epsilon_J}]$ is free 
of rank 1 over $\T^{\infty}_{\gm_{\infty}}$ and free over $\Lambda$.
Therefore   $\T^{\infty}_{\gm_{\infty}}$ is  free over $\Lambda$ and exact control 
holds.   \hfill$\square$

\begin{cor}  \label{thm-ord}
 Let $f\in S_{k+(p-1)k'}(Y_0(p^r))$ be a newform and $p$ be a prime not
dividing $\N_{F/\Q}(\gd)$, such that $p-1>\sum (k_\tau-1)$ and 
$\mathrm{\bf (LI_{\Ind\overline{\rho}})}$ holds. Then theorems A and B hold. 
\end{cor}

\newpage
\section*{List of symbols} 

\begin{multicols}{3}
\noindent $A$ \dotfill HBAV, Def.1.6

\noindent $A^t$ \dotfill dual HBAV 

\noindent $\cA$ \dotfill universal HBAV \S1.4

\noindent $\ga$ \dotfill ideal of $\go$

\noindent $B$ \dotfill standard Borel of $G$ 

\noindent $\gb$ \dotfill Lie algebra of  $B$

\noindent $\gc,\gc_+$ \dotfill  \S1.3

\noindent $c(f,\ga)$ \dotfill \S1.10

\noindent $\Cl_F$ \dotfill class group \S1.2

\noindent $\Cl_F^+$ \dotfill narrow class group \S1.1

\noindent $d$ \dotfill degree of $F$

\noindent $D$ \dotfill $\Res^{F}_{\Q} \Gm$ 

\noindent $D_{\gp}$ \dotfill decomposition group\S2.6

\noindent $\gd$ \dotfill different of $F$

\noindent $\mathcal{D}$ \dotfill \S3.4 

\noindent $E$ \dotfill large $p$-adic field

\noindent $f$ \dotfill Hilbert modular newform

\noindent $F$ \dotfill totally real number field

\noindent $\widetilde{F}$ \dotfill Galois closure of $F$

\noindent $\widetilde{F}'$ \dotfill before lemma 6.5

\noindent $\widehat{F}$ \dotfill \S3.4 

\noindent $\cF_D,\cF_B,\cF_G$ \dotfill functors \S2.3 

\noindent $g$ \dotfill Hilbert modular form

\noindent $G$ \dotfill $\Res^{F}_{\Q} \GL_2$ 

\noindent $\mathfrak{g}$ \dotfill Lie algebra of  $G$

\noindent $G^*$ \dotfill $G\times_{D}\Gm$ 

\noindent $g_J$ \dotfill \S1.11

\noindent $g_\tau$ \dotfill internal conjugate \S2.5

\noindent $\G_L$ \dotfill Galois group of $L$

\noindent $\mathfrak{G}$ \dotfill \S1.6

\noindent $H$ \dotfill above lemma 3.15

\noindent $\mathfrak{H}_F,\mathfrak{H}$ \dotfill  \S1.1

\noindent $h^+$ \dotfill  \S1.1

\noindent $\cH^1_{\dR}$ \dotfill  \S1.4

\noindent $I_{\gp}$ \dotfill inertia group at $\gp$ \S2.6

\noindent $J$ \dotfill subset of $J_F$

\noindent $J_F$ \dotfill set of infinite places of $F$

\noindent $J_F^i$ \dotfill \S3.5

\noindent $J_{F,\gp}$ \dotfill  \S2.6

\noindent $k=n+2t$ \dotfill  Def.\ref{weight} 

\noindent $k_0,n_0$ \dotfill Def.\ref{weight} 

\noindent $K$ \dotfill  CM field

\noindent $K_\infty,K_\infty^+$ \dotfill \S1.1

\noindent $K_1(\gn),K^1_1(\gn)$ \dotfill \S1.1

\noindent $K^\bullet_n$ \dotfill \S2.2

\noindent $\mathcal{K}^\bullet_n$ \dotfill \S2.3

\noindent $L$ \dotfill  field

\noindent $m=(k_0t-k)/2$ \dotfill  Def.\ref{weight} 

\noindent $M,M^1$ \dotfill connected HMV \S\ref{hmv}

\noindent $M'$ \dotfill \S1.4

\noindent $\overline{M},\overline{M^1}$ \dotfill \S1.6

\noindent $M^*,M^{1*}$ \dotfill \S1.8

\noindent $\gm$ \dotfill maximal ideal of $\T$ \S0.3

\noindent $\gm'$ \dotfill maximal ideal of $\T'$ \S6.2

\noindent $n$ \dotfill weight of $G$

\noindent $N$ \dotfill  normalizer of $T$ 

\noindent $\gn$ \dotfill level ideal

\noindent $\cO$ \dotfill integer ring of $E$

\noindent $\go$ \dotfill integer ring of $F$

\noindent $\go'$ \dotfill \S1.4

\noindent $\go_+^\times,\go_{\gn,1}^\times$ \dotfill \S0.4

\noindent $p$ \dotfill prime number 

\noindent $p(J)$ \dotfill \S2

\noindent $\cP$ \dotfill maximal ideal of $\cO$

\noindent $\gp$ \dotfill prime of $F$ dividing $p$

\noindent $q$ \dotfill power of $p$

\noindent $R$ \dotfill ring

\noindent $s$ \dotfill \S2.1

\noindent $S$ \dotfill base scheme

\noindent $S_k(\gn,\psi)$ \dotfill Def.1.3

\noindent $S_{\ga},T_{\ga}$ \dotfill Hecke operators \S1.10

\noindent $t=\sum_{\tau\in J_F}\tau$ \dotfill  Def.\ref{weight}

\noindent $T$ \dotfill standard torus of $G$ 

\noindent $\gt$ \dotfill Lie algebra of  $T$

\noindent $\T$ \dotfill Hecke algebra \S0.3

\noindent $\T'$ \dotfill reduced Hecke algebra 

\noindent $U$ \dotfill standard unipotent of $G$ 

\noindent $U(\gb),U(\mathfrak{g})$ \dotfill \S2.2
 
\noindent $\gu$ \dotfill Lie algebra of  $U$

\noindent $v$ \dotfill finite place of $F$

\noindent $V,V_n$ \dotfill  $G$-modules \S1.12

\noindent $\bV_n$ \dotfill local system \S1.12, \S2.1

\noindent $\overline{\cV},\overline{\cW}$ \dotfill \S2.3

\noindent $W,W_\mu$ \dotfill  $B$-modules 

\noindent $W_f$ \dotfill   \S2.1

\noindent $Y,Y^1$ \dotfill HMV \S\ref{hmv} 

\noindent $\alpha$  \dotfill $\mu_{\gn}$-level structure, Def.1.8

\noindent $\beta$  \dotfill \S2.4

\noindent $\gamma$  \dotfill element of $G(\R)$

\noindent $\Gamma_1(\gc,\gn)$  \dotfill \S1.1

\noindent $\Gamma^1_1(\gc,\gn)$ \dotfill \S1.1

\noindent $\delta,\delta_J$  \dotfill  \S1.12

\noindent $\Delta$  \dotfill $\Nm(\gn\gd)$

\noindent $\delta_{\gp},\varepsilon_{\gp}$  \dotfill tame characters of $I_{\gp}$ 

\noindent $\varepsilon$  \dotfill quadratic character

\noindent $\epsilon$  \dotfill unit of $F$

\noindent $\epsilon_J$  \dotfill \S1.11

\noindent $\widehat{\epsilon}_J$  \dotfill \S1.12


\noindent $\eta$  \dotfill id{\`e}le


\noindent $\iota$  \dotfill  Def.1.6

\noindent $\iota_p$  \dotfill embedding of $\overline{\Q}$ in  $\overline{\Q}_p$

\noindent $\kappa$  \dotfill residue field of $\cO$

\noindent $\lambda$  \dotfill $\gc$-polarization, Def.1.7

\noindent $\overline{\lambda}$  \dotfill $\gc$-polarization class,Def.1.7

\noindent $\mu$  \dotfill weight of $B$

\noindent $\nu$  \dotfill reduced norm  $G \rightarrow D$ \S\ref{hmfv}

\noindent $\underline{\nu}$  \dotfill \S1.5

\noindent $\xi$  \dotfill element of $F$

\noindent $\pi$  \dotfill \S1.4

\noindent $\varpi_v$  \dotfill uniformizer of $F_v$

\noindent $\overline{\pi}$  \dotfill \S1.6

\noindent $\rho=\rho_{f,p}$  \dotfill $p$-adic repr. \S0.1

\noindent $\overline{\rho}=\overline{\rho}_{f,p}$ \dotfill mod $p$ repr. \S0.1 

\noindent $\rho_1$  \dotfill \S6.1

\noindent $\sigma$  \dotfill cone \S1.6

\noindent $\Sigma$  \dotfill fan \S1.6

\noindent $\sigma_{i,\tau}$  \dotfill \S3.5

\noindent $\tau$  \dotfill  infinite place of $F$

\noindent $\phi$  \dotfill \S3.5

\noindent $\varphi$  \dotfill Hecke character

\noindent $\chi_n$  \dotfill \S2.2
  
\noindent $\psi$  \dotfill Hecke character, Def.1.3

\noindent $\omega$  \dotfill  mod $p$ cycl. char. 

\noindent $\Omega_f^{\pm}$  \dotfill  periods, Def.4.5

\noindent $\underline{\omega}$  \dotfill \S1.4

\noindent $\underline{\omega}^k$  \dotfill \S1.5

\noindent $(\enspace,\enspace)_{\gn}$  \dotfill  \S1.4
\end{multicols}

\newpage
\section*{Appendix}  
\renewcommand\thesection{A}
\renewcommand\thesubsection{A.\arabic{subsection}}
\setcounter{subsection}{0} 
\setcounter{theo}{0} 

The notations in this section are independent from the rest of the text.

Let $K$ be a local field  of characteristic zero  with a perfect
residue field  $\kappa$ of characteristic $p$.  Let
$W$ be the ring   of Witt vectors of $\kappa$ and $K_0$ 
be the fraction field of $W$. We  denote by $\sigma$ the 
Frobenius of $\kappa$, $W$ and $K_0$. 
Let $E$ be  a finite   extension  of $\Q_p$ in $\overline{\Q}_p$. 


\subsection{ $p$-adic representations.} 

A $p$-adic representation of $\G_K$ is a finite dimensional 
$\Q_p$-vector space, endowed with continuous action of $\G_K$.
The $p$-adic representations form an abelian category, denoted
 $\Rep(\G_K)$. We denote by $\Rep_E(\G_K)$ the subcategory 
of $\Rep(\G_K)$ consisting of $E$-linear representations.

There are several interesting subcategories of $\Rep(\G_K)$, 
as the one of Hodge-Tate  representations $\Rep_{\mathrm{HT}}(\G_K)$,
the one of de Rham representations $\Rep_{\dR}(\G_K)$,
the one of semi-stables  representations $\Rep_{\mathrm{st}}(\G_K)$,
and the one of  crystalline representations $\Rep_{\crys}(\G_K)$
$$\Rep(\G_K)\supset \Rep_{\mathrm{HT}}(\G_K) \supset 
\Rep_{\dR}(\G_K)\supset  \Rep_{\mathrm{st}}(\G_K) \supset 
\Rep_{\crys}(\G_K).$$

Let  $\Rep_{\overline{\Q}_p}(\G_K)$ be the category of continuous
representations of $\G_K$ on  finite dimensional $\overline{\Q}_p$-vector spaces.
All these representations are obtained by scalar extension from $\Rep_E(\G_K)$ (for some $E$).
By an abuse of language will still  call them   $p$-adic representations.

\subsection{ Hodge-Tate weights.}
Denote by $C$ the $p$-adic  completion of the algebraic closure
$\overline{K}$ of $K$. The $\G_K$-action on $\overline{K}$ 
extends by continuity to an action on $C$. 
Put $B_{\mathrm{HT}}=\oplus_i C(i)$, where the $\G_K$-action on 
$C(i)$ is twisted by the $i$-th power  of the cyclotomic character .

Let $V\in \Rep(\G_K)$. Then, by definition,  
$V\in \Rep_{\mathrm{HT}}(\G_K)$, if and only, if 
$\dim_K(V\otimes_{\Q_p}B_{\mathrm{HT}})^{\G_K}=\dim_{\Q_p}V$.
For $V\in \Rep_{\mathrm{HT}}(\G_K)$, we say  that $i$ is a  
Hodge-Tate weight of $V$, if  $V_i:=(V\otimes_{\Q_p}C(i))^{\G_K}\neq 0$ 
and we call  $h^i=\dim_K V_i$ its  multiplicity. 
We have a equality of $\G_K$-modules 
$V\otimes_{\Q_p}C=\bigoplus_i V_i\otimes_K C(-i)$.

If $V\in \Rep_E(\G_K)$, then   for all $i\in \Z$ 
$V_i=(V\otimes_{\Q_p}C(i))^{\G_K}$ is a $E\otimes_{\Q_p} K$-module 
in a natural way. It is not free in general. By decomposing the 
 $\Q_p$-algebra  $E\otimes_{\Q_p} K$ as a product of fields $\prod_j L_j$
(endowed with injections  $\sigma : E\hookrightarrow L_j$, 
$\tau : K\hookrightarrow L_j$), we obtain :
$$(V\otimes_{\Q_p}C(i))^{\G_K}
\underset{E\otimes K}{\otimes} L_j=
(V\otimes_{E}C(i))^{\G_{L_j}} $$

There is another way to index the  Hodge-Tate weights 
that is more appropriate to  the modular case. Consider the 
functor $\Rep_E(\G_K)\rightarrow \Rep_{\overline{\Q}_p}(\G_K)$
sending $V$ to $V_{\overline{\Q}_p}:=V\otimes_E \overline{\Q}_p$.

\begin{defin}\label{poids-ht} For all 
 $\tau: K\hookrightarrow \overline{\Q}_p$ 
we put $h_{\tau,i}=\dim_{\overline{\Q}_p}
(V_{\overline{\Q}_p}\otimes_{\tau,K} C(i))^{\G_K}$. The integer 
$h_{\tau,i}$ is called the multiplicity of $i$ as Hodge-Tate weight 
of $V$. For all $\tau$, we have $\sum_{i\in\Z} h_{\tau,i}= \dim_E V$. 
\end{defin}

\begin{ex}
Assume that $E=\Q_p$. Then   $V_i=(V\otimes_{\Q_p}C(i))^{\G_K}$
 is a $K$-vector space  and the $i$-th  Hodge-Tate  number is given by
 $h^i=\dim_K(V_i)$, $i\in \Z$. 

If we change the action of $\G_K$ on $C$ by an automorphism
$g\mapsto \tau^{-1} g\tau$, with $\tau\in \G_{\Q_p}$, 
then  we send $V_i$ onto $V_i^\tau$ by $v\otimes a\mapsto
v\otimes \tau(a)$ that does not  change  $h^i$ 
(since  the cyclotomic character  is invariant by  $g\mapsto \tau^{-1} g\tau$). 
\end{ex}

\subsection{Crystalline representations and filtered modules.}
The category  $\Rep_{\crys}(\G_K)$ of  crystalline representations 
is the $p$-adic  analogue of  the unramified  $l$-adic representations.

\begin{defin}
(i) A filtered $\phi$-module  over $K$ is  a $K_0$-vector space 
 $D$ of finite dimension, endowed with a $\sigma$-linear bijective map  
$\phi:D\rightarrow D$  and a  filtration $\Fil^i D$ of 
$D_K=D\otimes_{K_0} K$ which is decreasing ($\Fil^i D\supset \Fil^{i+1} D$), 
exhaustive ($\cup \Fil^i D=D_K$) and separated ($\cap \Fil^i D=0$).
We denote by $\MF_K$ the additive category  filtered $\phi$-module  over $K$.

(ii) A filtered $\phi$-module $D$   over $K$ is called  weakly admissible, if
it contains a $W$-lattice  $M$, such that 
$\sum p^{-i}\phi(\Fil^i D\cap M)=M$.
Such a lattice is called {\it strongly divisible} (or adapted to $D$). 
The weakly admissible filtered $\phi$-module  over $K$ form 
a full subcategory of $\MF_K$, denoted by $\MF_K^f$.
\end{defin}

\begin{rque} 
To an object  $D\in \MF_K$ one can associate  Newton and Hodge
polygons and the notion of be  weakly admissible can be 
expressed in terms of inequalities between these two polygons.
\end{rque}

Fontaine's theory gives  an equivalence of categories  $D_{\crys}$ between 
$\Rep_{\crys}(\G_K)$ and a certain full subcategory of 
admissible objects $\MF_K^a$  of $\MF_K$. The Hodge-Tate weights of
$V\in \Rep_{\crys}(\G_K)$ are given by the jumps of the filtration on 
$D_{\crys}(V)$.

It is known that  admissible implies weakly admissible, and
recently Colmez and Fontaine proved the converse, in the 
more general semi-stable case. When $K$ is an  unramified extension of
$\Q_p$ and the lenght of the filtration is $\leq p-1$, this has been 
established earlier by   Fontaine and Laffaille \cite{FoLa}.

\subsection{Crystalline representations modulo $p$.} 
In the sequel, we assume  $K$ to be unramified ($K=K_0$).
 Fontaine and Laffaille \cite{FoLa} have introduced

\begin{defin}
(i) A filtered  $F$-module  over $W$ is defined by the following data 

$\bullet$  a  $W$-module $M$, 

$\bullet$  a  filtration $\Fil^iM$   by $W$-submodules which is decreasing
exhaustive and separated,

$\bullet$  $\sigma$-linear maps $\varphi_i: \Fil^iM\rightarrow M$ satisfying 
$\varphi_i|_{\Fil^{i+1}M}=p\varphi_{i+1}$.
 
We denote by $\MF_W$ the $\Z_p$-linear additive category of  filtered 
 $F$-modules  over $W$. 

(ii) We define two full abelian subcategories   
$\MF_{W,lf}\!\subset\! \MF_{W,tf}\!\subset\! \MF_W$, by  $M\in \MF_W$

$$M\in \MF_{W,tf}\iff \begin{cases} M\text{ is of finite  type  over } W,\\
\Fil^iM\text{ are  direct factors in }M, \\
\sum \varphi_i(\Fil^iM)=M.\end{cases}, $$

$$M\in \MF_{W,lf}\iff \begin{cases} M\text{ is of finite length  over } W,\\
\sum \varphi_i(\Fil^iM)=M.\end{cases}. $$
\end{defin}

\begin{rque}
Let $M$ be a strongly divisible lattice  for  $D\in \MF_K^f$, and put
$\Fil^iM=\Fil^iD\cap M$ and $\varphi_i=\phi/p^i$. Then we  have 
$M\in \MF_{W,tf}$. If moreover $M$ is a free $W$-module of  finite type,
then   we have $M/p^rM\in  \MF_{W,lf}$, for all $r \in \n$. 
\end{rque}

\subsection{The tame inertia.}
We choose to uniformize the Local Class Field Theory isomorphism by sending 
a uniformizer to the geometric Frobenius. 

 The tame inertia
$I_K^t$ is the quotient of the inertia group $I_K$ by its  maximal pro-$p$ subgroup
(called the wild inertia). There is a canonical isomorphism 
 $I_K^t\simeq \limproj \F_{p^h}^\times$ (see \cite{serre2}). 
A tame character of level $h\in \n$ is a character of
 $I_K^t$  factoring through  $\F_{p^h}^\times$. 

By the Local Class Field Theory, the inertia group of the maximal 
abelian extension of $K$ is isomorphic to $W^\times$. 
Thus we get a homomorphism $I_K\rightarrow W^\times \rightarrow
\kappa^\times$  equal to the  tame 
character $I_K\rightarrow I_K^t \rightarrow \kappa^\times$ of level $[\kappa:\F_p]$
(see \cite{serre2}).

\subsection{Fontaine-Laffaille's theory  and a theorem of Wintenberger.}
In  \cite{FoLa} Fontaine and Laffaille have introduced a contravariant functor
$$V_{\mathrm{FL}} :\MF_{W,tf}\longrightarrow \Rep_{\Z_p}(\G_K)\text{ , such that}$$

$\bullet$ the restriction to $\MF_{W,lf}^{<p-1}$ is exact and 
fully faithful,

$\bullet$ if $M\in \MF_{W,lf}^{<p-1}$, then $\mathrm{length}_W M=
\mathrm{length}_{\Z_p} V_{\mathrm{FL}}(M)$,

$\bullet$ if $M\in \MF_{W,tf}^{<p-1}$ is free, 
then  $V_{\mathrm{FL}}(M)$ is free and 
$\mathrm{rank}_W M= \mathrm{rank}_{\Z_p} V_{\mathrm{FL}}(M)$.

Let $X$ be the abelian  group of periodic map  $\xi:\Z\rightarrow \Z$.
By a result of Wintenberger \cite{Wint} we can decompose a 
filtered $F$-module  $M$ of  finite type over $W$  as a sum of 
isotypic components indexed by $X$,  $M=\bigoplus_{\xi\in X} M_{\xi}$.

\subsection{Hodge-Tate weights and Fontaine-Laffaille weights.}\label{ht-modere}
 The aim of this paragraph is to explain how the  theory of  Fontaine 
and Laffaille relates  the  Hodge-Tate weights of a crystalline 
 representation  to the weights of the tame inertia acting on the 
semi-simplification of its reduction modulo $p$. 
This formulation is due to Wintenberger \cite{Wint}.

Let $V$ be a  $p$-adic crystalline representation of   Hodge-Tate weights
 between $0$ and $p-1$ 
($V\in \Rep_{\mathrm{crys}}^{<p-1}(\G_K)$).
We can associate to it  $D=D_{\mathrm{crys}}(V)\in \MF_K^f$.
The  multiplicity of  $i\in \Z$ as  a  Hodge-Tate weight of $V$ 
is equal to $h^i=\dim_K(\Fil^i D)-\dim_K(\Fil^{i+1} D)$. 

Let $M\in \MF_{W,tf}^{<p-1}$ be a  $W$-lattice adapted to $D$.
By a theorem of Wintenberger we have two natural decompositions of  $M$ :
$$M=\bigoplus_{i\in \Z} M_i, \enspace 
M=\bigoplus_{\xi\in X} M_{\xi}.$$

Let's define $D_i=M_i\otimes_W K_0$ and $D_\xi=M_\xi\otimes_W K_0$.
Then we have $h^i=\dim D_i$ and 
$$D_i=\bigoplus_{\xi\in X, \xi(0)=i} D_{\xi}.$$

For each $r\in \n$, we have $M/p^rM\in \MF_{W,lf}^{<p-1}$. 
Put $L_r=V_{\mathrm{FL}}(M/p^rM)$ and $L=\limproj L_r$. 
Then  $L$ is a lattice of $V$ and  by construction we have $L/pL=L_1$. 

Moreover $M/pM=\oplus_{\xi\in X} M_{\xi}/pM_{\xi}$. 
But  $M_{\xi}/pM_{\xi}$ is a sum of copies  of a  simple object of
 $ \MF_{W,lf}^{<p-1}$ and  $V_{\mathrm{FL}}(M_{\xi}/pM_{\xi})$ 
is equal  to the sum of the same number of copies of the 
tame character $\theta(\xi)$. We deduce that the tame characters  occurring 
in $L/pL$ correspond exactly  to the $\xi$ occurring  in $M$. 
By the theorem of  Brauer-Nesbitt the semi-simplification of $L/pL$
does not depend  of the particular choice of a lattice.

\begin{theo} {\rm (Fontaine-Laffaille)} 
Assume  $V\in \Rep_{\crys}^{<p-1}(\G_K)$
has  Hodge-Tate weights $i_1,...,i_r$ (with multiplicities). 
Let $L$ be a stable lattice  of $V$ and consider  the tame inertia action
on  $(L/pL)$. Take a decomposition  
$(L/pL)^{\mathrm{ss}}=\oplus \overline{L}_j$, such that 
the tame inertia acts on $\overline{L}_j$ by a certain tame
character of level $h^j=\dim_{\F_p}\overline{L}_j$ and  weights
  $i_1^j,..,i_{h^j}^j$ ($r=\sum h^j$).  Then  the multisets $\{i_1,...,i_r\}$  and 
$\{i_k^j\enspace |\enspace 1\leq k\leq h^j\}$ are equal.
\end{theo}

The result remains  valid for  $V\in \Rep_{\crys,E}^{<p-1}(\G_K)$
(see  Def.\ref{poids-ht}).

\bibliographystyle{siam}

\bigskip

Universit{\'e} Paris 7,  UFR de Math{\'e}matiques,  Case 7012, 
2 place Jussieu,  75251 PARIS

\medskip
 Email : \texttt{dimitrov@math.jussieu.fr}

\end{document}